\documentclass[11pt]{amsart}
\usepackage{amssymb}
\usepackage{lscape}
\usepackage[all]{xy}
\usepackage{array}
\usepackage{epic,eepic}

\newtheorem{Thm}{Theorem}[section]
\newtheorem{Lem}[Thm]{Lemma}
\newtheorem{Cor}[Thm]{Corollary}
\newtheorem{Prop}[Thm]{Proposition}
\newtheorem{Conj}[Thm]{Conjecture}
\theoremstyle{definition} 
\newtheorem{Rem}[Thm]{Remark}

\newcommand{\A}{\mathbb{A}}
\newcommand{\D}{\mathbb{D}}
\newcommand{\E}{\mathbb{E}}

\newcommand{\Z}{\mathbb{Z}}

\newcommand{\N}{\mathbb{N}}
\newcommand{\C}{\mathbb{C}}

\newcommand{\df}{\colon}

\newcommand{\id}{{\rm id}}

\newcommand{\F}{{\mathcal F}}
\newcommand{\U}{{\mathcal U}}
\newcommand{\G}{{\mathcal G}}
\newcommand{\T}{{\mathcal T}}
\newcommand{\X}{{\mathcal X}}
\newcommand{\Y}{{\mathcal Y}}

\newcommand{\RR}{{\mathcal R}}

\newcommand{\MM}{{\mathcal M}}
\newcommand{\cA}{{\mathcal A}}
\newcommand{\cS}{{\mathcal S}}
\newcommand{\cP}{{\mathcal P}}

\newcommand{\CC}{{\mathcal C}}
\newcommand{\stCC}{{\underline{\mathcal C}}}

\newcommand{\g}{\mathfrak{g}}
\newcommand{\n}{\mathfrak{n}}
\newcommand{\h}{\mathfrak{h}}

\newcommand{\mm}{\mathbf{m}}
\newcommand{\nn}{\mathbf{n}}
\newcommand{\ee}{\mathbf{e}}
\newcommand{\aaa}{\mathbf{a}}
\newcommand{\bb}{\mathbf{b}}

\newcommand{\sss}{\mathbf{s}}
\newcommand{\ttt}{\mathbf{t}}

\newcommand{\ii}{{\mathbf i}}
\newcommand{\jj}{{\mathbf j}}
\newcommand{\kk}{{\mathbf k}}
\newcommand{\yy}{{\mathbf y}}

\newcommand{\hI}{\widehat{I}}

\newcommand{\LL}{\Lambda}
\newcommand{\GG}{\Gamma}

\newcommand{\vpi}{\varpi}
\newcommand{\vph}{\varphi}
\newcommand{\la}{\lambda}

\newcommand{\ad}{\operatorname{ad}}

\newcommand{\md}{\operatorname{mod}}
\newcommand{\nil}{\operatorname{nil}}
\newcommand{\rep}{\operatorname{rep}}

\newcommand{\add}{\operatorname{add}}
\newcommand{\Gen}{\operatorname{Fac}}
\newcommand{\Cogen}{\operatorname{Sub}}

\newcommand{\gldim}{\operatorname{gl.dim}}
\newcommand{\pdim}{\operatorname{proj.dim}}
\newcommand{\soc}{\operatorname{soc}}

\newcommand{\Hom}{\operatorname{Hom}}
\newcommand{\Ext}{\operatorname{Ext}}
\newcommand{\End}{\operatorname{End}}

\newcommand{\Ima}{\operatorname{Im}}
\newcommand{\Ker}{\operatorname{Ker}}
\newcommand{\Coker}{\operatorname{Coker}}
\newcommand{\tp}{\operatorname{top}}

\newcommand{\dm}{{\rm dim}\,}
\newcommand{\dimv}{\underline{\dim}}

\newcommand{\Span}{\operatorname{Span}}

\newcommand{\bil}[1]{\langle #1\rangle}
\newcommand{\abs}[1]{\left| #1\right|}
\newcommand{\ebrace}[1]{\langle #1 \rangle}

\newcommand{\op}{{\rm op}}

\newcommand{\irr}{\operatorname{Irr}}

\newcommand{\bsm}{\begin{smallmatrix}}
\newcommand{\esm}{\end{smallmatrix}}

\newcommand{\bbsm}{\left[\begin{smallmatrix}}
\newcommand{\besm}{\end{smallmatrix}\right]}

\newcommand{\bbm}{\begin{matrix}}
\newcommand{\ebm}{\end{matrix}}

\newcommand{\GL}{\operatorname{GL}}

\newcommand{\Gmin}{{G^{\rm min}}}
\newcommand{\Nmin}{N^{\rm min}}
\newcommand{\Norm}{\operatorname{Norm}}
\newcommand{\Lie}{\operatorname{Lie}}

\newcommand{\cE}{\mathcal{E}}

\newcommand{\ba}{\mathbf{a}}
\newcommand{\bd}{\mathbf{d}}

\newcommand{\ff}{\mathbf{f}}
\newcommand{\bi}{\mathbf{i}}

\newcommand{\ra}{\rightarrow}

\def\shuff#1#2{\mathbin{
\hbox{\vbox{ \hbox{\vrule \hskip#2 \vrule height#1 width 0pt}%
\hrule}%
\vbox{ \hbox{\vrule \hskip#2 \vrule height#1 width 0pt
\vrule}%
\hrule}%
}}}

\def\SHUF{{\mathchoice{\shuff{7pt}{3.5pt}}%
{\shuff{6pt}{3pt}}%
{\shuff{4pt}{2pt}}%
{\shuff{3pt}{1.5pt}}}}%

\def\shuffle{\,\SHUF\,\,}

\makeindex

\begin{document}

\date{June 10, 2011}

\bigskip
\title[Kac-Moody groups and cluster algebras]
{Kac-Moody groups and cluster algebras}

\author{Christof Gei{\ss}}
\address{Christof Gei{ss}\newline
Instituto de Matem\'aticas\newline
Universidad Nacional Aut{\'o}noma de M{\'e}xico\newline
Ciudad Universitaria\newline
04510 M{\'e}xico D.F.\newline
Mexico}
\email{christof@math.unam.mx}

\author{Bernard Leclerc}
\address{Bernard Leclerc\newline
LMNO, Universit\'e de Caen\newline
CNRS UMR 6139\newline
F-14032 Caen Cedex\newline
France}
\email{bernard.leclerc@unicaen.fr}

\author{Jan Schr\"oer}
\address{Jan Schr\"oer\newline
Mathematisches Institut\newline
Universit\"at Bonn\newline
Endenicher Allee 60\newline
D-53115 Bonn\newline
Germany}
\email{schroer@math.uni-bonn.de}



\begin{abstract}
Let $Q$ be a finite quiver without oriented cycles, let
$\LL$ be the associated preprojective algebra,
let $\g$ be the associated Kac-Moody Lie algebra with
Weyl group $W$, and let $\n$ be the positive part of $\g$.
For each Weyl group element $w$, a
subcategory $\CC_w$ of $\md(\LL)$ was introduced by Buan, Iyama,
Reiten and Scott.
It is known that $\CC_w$ is a Frobenius category
and that
its stable category $\underline{\CC}_w$ is a 
Calabi-Yau category of dimension two.
We show that $\CC_w$ yields a cluster algebra structure
on the coordinate ring $\C[N(w)]$ of
the unipotent group $N(w) := N \cap (w^{-1}N_-w)$.
Here $N$ is the pro-unipotent pro-group with Lie algebra
the completion $\widehat{\n}$ of $\n$.
One can identify $\C[N(w)]$ with a subalgebra of $U(\n)_{\rm gr}^*$, 
the graded dual of the universal
enveloping algebra $U(\n)$ of $\n$.
Let $\cS^*$ be the dual of Lusztig's semicanonical
basis $\cS$ of $U(\n)$.
We show that all cluster monomials of $\C[N(w)]$
belong to $\cS^*$, and that
$\cS^* \cap \C[N(w)]$ is a $\C$-basis 
of $\C[N(w)]$.
Moreover, we show that the cluster algebra obtained from 
$\C[N(w)]$ by formally inverting the generators
of the coefficient ring is isomorphic to the algebra 
$\C[N^w]$ of regular
functions on the unipotent cell $N^w$ 
of the Kac-Moody group with Lie algebra $\g$.
We obtain a corresponding dual semicanonical 
basis of $\C[N^w]$.
As one application we obtain a basis for each
acyclic cluster algebra, which contains all cluster
monomials in a natural way.
\end{abstract}

\maketitle

\setcounter{tocdepth}{1}

\tableofcontents

\parskip2mm



\section{Introduction}


\subsection{}
This is the continuation of an extensive project to obtain a 
better understanding
of the relations between the following topics:
\begin{itemize}

\item[(i)]
Representation theory of quivers,

\item[(ii)]
Representation theory of preprojective algebras,

\item[(iii)]
Lusztig's (semi)canonical basis of universal enveloping algebras,

\item[(iv)]
Fomin and Zelevinsky's theory of cluster algebras,

\item[(v)]
Frobenius categories and 2-Calabi-Yau categories, 

\item[(vi)]
Cluster algebra structures on coordinate algebras of 
unipotent groups, Bruhat cells and flag varieties.

\end{itemize}
The topics (i) and (iii) are closely related.
The numerous connections have been studied by many authors.
Let us just mention Lusztig's work on canonical bases of quantum
groups, and Ringel's Hall algebra approach to quantum groups.
An important link between (ii) and (iii), due to Lusztig
\cite{Lu1,Lu2} and Kashiwara and Saito \cite{KSa}
is that the elements of the (semi)canonical basis are naturally
parametrized by the irreducible components of the varieties of
nilpotent representations of a preprojective algebra. 

Cluster algebras were invented by Fomin and Zelevinsky
\cite{BFZ,FZ1,FZ2},
with the aim of providing a new algebraic and combinatorial setting
for canonical bases and total positivity.
One important breakthrough was the insight that 
the class of acyclic cluster algebras with a skew-symmetric exchange matrix
can be categorified using the so-called cluster categories.
Cluster categories were introduced by Buan, Marsh, Reineke, Reiten and
Todorov \cite{BMRRT}, see also \cite{Ke}.
In a series of papers by some of these authors and also by Caldero and
Keller \cite{CK,CK2}, it was established that cluster categories have
all necessary properties to provide the mentioned categorification.
We refer to the nice overview article \cite{BM} for more details on
the development of this beautiful theory which established a
strong connection between the topics (i), (iv) and (v).
More recently, 
a different and more general type of categorification using representations of
quivers with potentials was developed by Derksen, Weyman and Zelevinsky
\cite{DWZ1,DWZ2}.
This provides another strong link between topics (i) and (iv).

In \cite{GLSRigid} we showed that the representation theory of
preprojective algebras $\LL$ of Dynkin type 
(\textit{i.e.}~type $\A$, $\D$ or $\E$)
is also closely related to cluster algebras.
We proved that $\md(\LL)$ can be regarded as a categorification
of a natural (upper) cluster structure on the polynomial algebra $\C[N]$.
Here $N$ is a maximal unipotent subgroup of a complex Lie group 
of the same type as $\LL$. Let $\n$ be its Lie algebra,
and let $U(\n)$ be the universal enveloping algebra of $\n$.
The graded dual $U(\n)^*_{\rm gr}$ can be identified with the coordinate
algebra $\C[N]$.
By means of our categorification, we were able to prove that 
all the cluster monomials of $\C[N]$ belong to the dual
of Lusztig's semicanonical basis of $U(\n)$.
Note that the cluster algebra $\C[N]$ is in general not acyclic.

The aim of this article is a vast generalization of these results to the more
general setting of symmetric Kac-Moody groups and their unipotent cells.
We also provide additional tools for studying the associated categories
and cluster structures.
For many cluster algebras
we construct a basis (called {\it dual semicanonical basis}) which
contains all cluster monomials in a natural way.
In particular, we obtain such a basis for all acyclic cluster algebras.
Also, we construct a dual PBW-basis of the cluster algebras involved.
This provides another close link between Lie theory and the representation
theory of preprojective algebras.
We show that the coordinate rings $\C[N(w)]$ and $\C[N^w]$ 
are genuine cluster algebras in a natural way, and not just upper cluster 
algebras in the sense of \cite{BFZ}.

Let us give some more details.
We consider preprojective algebras 
$\LL=\LL_Q$ attached to quivers $Q$ which are not necessarily of Dynkin type. 
These algebras are therefore infinite-dimensional in general.
The category $\nil(\LL)$ of all finite-dimensional nilpotent
representations of $\LL$ is then too large to be related to a 
cluster algebra of finite rank.
Moreover, it does not have projective or injective objects,
and it lacks an Auslander-Reiten translation. 
However, Buan, Iyama, Reiten and Scott \cite{BIRS} have attached
to each element $w$ of the Weyl group $W = W_Q$ of $Q$
a subcategory $\CC_w$ of $\nil(\LL)$.
They show that
the categories $\CC_w$ are
Frobenius categories
and the corresponding stable categories ${\stCC}_w$ are 
Calabi-Yau categories of dimension two.
(These results were also discovered and proved independently
in \cite{GLSUni1} in the special case when $w$ is an
{\it adaptable} element of $W$.)
Each subcategory $\CC_w$ contains a distinguished maximal rigid $\LL$-module
$V_\ii$ associated to each reduced expression $\ii = (i_r,\ldots,i_1)$ of $w$.
(A module $X$ is called {\it rigid} if $\Ext_\LL^1(X,X) = 0$.)

Special attention is given to the algebra $B_\ii := \End_\LL(V_\ii)^\op$,
which turns out to be quasi-hereditary.
There is an equivalence between $\CC_w$ 
and the category
of $\Delta$-filtered $B_\ii$-modules.
This allows us to describe mutations of maximal
rigid $\LL$-modules in $\CC_w$ in terms of the $\Delta$-dimension
vectors of the corresponding $B_\ii$-modules.

To the subcategory $\CC_w$ we associate
a cluster algebra $\cA(\CC_w)$ which in general is
not acyclic,
and we show that 
$\CC_w$ can be seen as a categorification of
the cluster algebra $\cA(\CC_w)$. 
Each of the modules $V_\ii$ provides an initial seed
of this cluster algebra.
(As a very special case, we also obtain in this way a new 
categorification of every acyclic cluster algebra with a
skew-symmetric exchange matrix and a certain choice of coefficients.)
The proof relies on the fact that
the algebra $\cA(\CC_w)$ has a natural realization
as a certain subalgebra of the graded dual $U(\n)^*_{\rm gr}$, where 
$\n$ is now the positive part of the symmetric Kac-Moody Lie algebra 
$\g = \n_- \oplus \h \oplus \n$ 
of the same type as $\LL$. We show that again all the
cluster monomials belong to the dual
of Lusztig's semicanonical basis of $U(\n)$. 

Next, we prove that $\cA(\CC_w)$ has a 
simple monomial basis
coming from the objects of the additive
closure $\add(M_\ii)$, where $M_\ii = M_1 \oplus \cdots \oplus M_r$
is another $\LL$-module in $\CC_w$ associated to a reduced expression
$\ii$ of $w$.
The modules $M_k$ are rigid, but $M_\ii$ is not rigid, except in some
trivial cases.
We call it the {\it dual PBW-basis} of $\cA(\CC_w)$,
and regard it as a generalization (in the dual setting) of the bases 
of $U(\n)$ constructed by Ringel in terms of quiver representations,
when $\g$ is finite-dimensional \cite {Ri6}. 
We use this to prove that $\cA(\CC_w)$ is spanned
by a subset of the dual semicanonical basis of $U(\n)^*_{\rm gr}$.
Thus, we obtain a natural basis of $\cA(\CC_w)$
containing all the cluster monomials.
We call it the {\it dual  semicanonical basis} of $\cA(\CC_w)$.
We prove that $\cA(\CC_w)$ is isomorphic
to the coordinate ring of the finite-dimensional unipotent
subgroup $N(w)$ of the symmetric Kac-Moody group attached to $\g$.
Moreover, we show that the cluster algebra obtained from 
$\cA(\CC_w)$ by formally inverting the generators
of the coefficient ring is isomorphic to the algebra of regular
functions on the unipotent cell $N^w$ of the Kac-Moody group.    
This solves Conjecture IV.3.1 of \cite{BIRS}.

Note also that in the Dynkin case the unipotent cells $N^w$
are closely related to the double Bruhat cells of type $(e,w)$,
whose coordinate ring is known to be an upper cluster algebra
by a result of [BFZ]. However, our proof is different and
shows that $\C[N^w]$ is not only an upper cluster algebra but a
genuine cluster algebra.

Finally, we explain how the results of this
paper are related to those of \cite{GLSFlag}, in which
a cluster algebra structure on the coordinate ring of the unipotent radical
$N_K$ of a parabolic subgroup of a complex simple algebraic group of
type $\A, \D, \E$ was introduced. 
We give a proof of Conjecture~9.6 of \cite{GLSFlag}.

\subsection{Remark}
Our preprint \cite{GLSUni1} contains special cases of the main
results of this article:
When $w$ is an {\it adaptable} Weyl group element,
we constructed and studied the subcategories $\CC_w$ 
independently of \cite{BIRS}, using different methods.
For this case, \cite{GLSUni1} contains a 
proof of \cite[Conjecture IV.3.1]{BIRS}.
Since \cite{GLSUni1} is already cited in several published
articles, we decided to keep it on the arXiv as
a convenient reference, but it will not be published
in a journal.

\subsection{Notation}
Throughout let $K$ be an algebraically closed field.
For a $K$-algebra $A$ let $\md(A)$ be the category of
finite-dimensional left $A$-modules.
By a {\it module} we always mean a finite-dimensional left module.
Often we do not distinguish between a module and its isomorphism 
class.
Let $D := \Hom_K(-,K)\df \md(A) \to \md(A^\op)$ be the usual
duality. 

For a quiver $Q$ let $\rep(Q)$ be the category of finite-dimensional
representations of $Q$ over $K$.
It is well known that
we can identify $\rep(Q)$ and $\md(KQ)$.

By a {\it subcategory} we always mean a full subcategory.
For an $A$-module $M$ let $\add(M)$ be the subcategory of all
$A$-modules
which are isomorphic to finite direct sums of direct summands of $M$.
A subcategory $\U$ of $\md(A)$ 
is an {\it additive subcategory} if any finite direct
sum of modules in $\U$ is again in $\U$.
By $\Gen(M)$ 
(resp. $\Cogen(M)$) 
we denote the subcategory of
all $A$-modules $X$ such that there exists some $t \ge 1$ and some
epimorphism $M^t \to X$ (resp. monomorphism $X \to M^t$).

For an $A$-module $M$ let $\Sigma(M)$ be the number of isomorphism classes
of indecomposable direct summands of $M$.
An $A$-module is called 
{\it basic} 
if it can be written as a direct sum
of pairwise non-isomorphic indecomposable modules.

For an $A$-module $M$ and a simple $A$-module $S$ let
$[M:S]$ be the Jordan-H\"older multiplicity of $S$
in a composition series of $M$.
Let $\dimv(M) := \dimv_A(M) := ([M:S])_S$ be the {\it dimension vector}
of $M$, where $S$ runs through all isomorphism classes of simple
$A$-modules.

For a set $U$ we denote its cardinality by $|U|$.
If $f\df X \to Y$ and $g\df Y \to Z$ are maps, then the composition
is denoted by $gf = g \circ f\df X \to Z$.

If $U$ is a subset of a $K$-vector space $V$, then let
$\Span_K\ebrace{U}$ be the subspace of $V$ generated by
$U$.

By $K(X_1,\ldots,X_r)$ (resp. $K[X_1,\ldots,X_r]$)
we denote the field of rational functions (resp. the polynomial ring)
in the variables $X_1,\ldots,X_r$ with coefficients in $K$.

Let $\C$ be the field of complex numbers, and let $\N = \{0,1,2,\ldots\}$
be the natural numbers, including $0$.
Set $\N_1 := \N \setminus \{ 0 \}$.

Recommended introductions to representation theory of finite-dimensional
algebras and Auslander-Reiten theory are the books \cite{ARS,ASS,GR,Ri1}.


\section{Definitions and known results}


\subsection{Preprojective algebras and nilpotent varieties}\label{prenil}
Let $Q = (Q_0,Q_1,s,t)$ be a finite quiver without 
oriented cycles.
(As usual, $Q_0$ is the set of vertices, $Q_1$ is the set of arrows,
an arrow $a \in Q_1$ starts in a vertex $s(a)$ and terminates in
$t(a)$.)
Let
$$
\LL = \LL_Q = K\overline{Q}/(c)
$$ 
be the associated 
{\it preprojective algebra}.
We assume that $Q$ is connected and has vertices
$Q_0 = \{1,\ldots,n\}$.
Here $K$ is an algebraically closed field,  
$K\overline{Q}$ is the path algebra of the 
{\it double quiver } $\overline{Q}$ 
of $Q$ which is obtained
from $Q$ by adding to each arrow $a\df i \to j$ in $Q$ an arrow
$a^*\df j \to i$ pointing in the opposite direction, and $(c)$ is
the ideal generated by the element
$$
c = \sum_{a \in Q_1} (a^*a - aa^*).
$$
Clearly, the path algebra $KQ$ is a subalgebra of $\LL$. 
Let
$\pi_Q\df \md(\LL) \to \md(KQ)$
be the corresponding restriction functor.

A $\LL$-module $M$ is called {\it nilpotent} 
if a composition
series of $M$ contains only the simple modules $S_1,\ldots,S_n$ 
associated to the vertices of $Q$.
Let $\nil(\LL)$ 
be the abelian category of finite-dimensional
nilpotent $\LL$-modules.

Let $d = (d_1, \ldots, d_n) \in \N^n$.
By 
$$
\rep(Q,d) = \prod_{a \in Q_1} 
\Hom_{K}(K^{d_{s(a)}},
K^{d_{t(a)}}) 
$$
we denote the affine space of representations of $Q$ with dimension
vector $d$.
Furthermore, let $\md(\LL,d)$ 
be the affine variety of elements
$$
(f_a,f_{a^*})_{a \in Q_1} \in
\prod_{a \in Q_1} \left( \Hom_K(K^{d_{s(a)}},
K^{d_{t(a)}}) \times \Hom_K(K^{d_{t(a)}},
K^{d_{s(a)}}) \right) 
$$
such that the following holds:
\begin{itemize}

\item[(i)]
For all $i \in Q_0$ we have
$$ 
\sum_{a \in Q_1: s(a) = i} f_{a^*}f_a =
\sum_{a \in Q_1: t(a) = i} f_a f_{a^*}.
$$

\end{itemize}
By $\LL_d$ 
we denote the variety
of all 
$(f_a,f_{a^*})_{a \in Q_1} \in \md(\LL,d)$ such that
the following condition holds:
\begin{itemize}

\item[(ii)]
There exists some $N$ such that for each path $a_1a_2 \cdots a_N$ of
length $N$ in the double quiver $\overline{Q}$ of $Q$
we have $f_{a_1}f_{a_2} \cdots f_{a_N} = 0$.

\end{itemize}
(It is not difficult to check that $\LL_d$ is indeed an affine variety,
namely for a fixed $d$ we can choose $N=d_1+\cdots +d_n$ in condition
(ii) above.)
If $Q$ is a Dynkin quiver, then (ii) follows already from condition
(i).
One can regard (ii) as a nilpotency condition, which explains why
the varieties $\LL_d$ are often called {\it nilpotent varieties}.
Note that $\rep(Q,d)$ can be considered as a subvariety of $\LL_d$.
In fact $\rep(Q,d)$ forms an irreducible component of $\LL_d$.
Lusztig \cite[Section 12]{Lu1} proved that 
all irreducible components of $\LL_d$ have the same dimension, namely
$$
\dm \rep(Q,d) = \sum_{a \in Q_1} d_{s(a)}d_{t(a)}.
$$
One can interpret $\LL_d$ as the variety of nilpotent
$\LL$-modules with dimension vector $d$.
The group 
$$
\GL_d = \prod_{i=1}^n \GL_{d_i}(K)
$$ 
acts on $\md(\LL,d)$, 
$\LL_d$ and $\rep(Q,d)$ by conjugation.
Namely, for 
$g = (g_1,\ldots,g_n) \in \GL_d$ and
$x = (f_a,f_{a^*})_{a \in Q_1} \in \md(\LL,d)$ define
$$
g.x := (g_{t(a)}f_a g_{s(a)}^{-1},
g_{s(a)}f_{a^*} g_{t(a)}^{-1})_{a \in Q_1}.
$$
The action on $\LL_d$ and $\rep(Q,d)$ is obtained via restriction.
The isomorphism classes of $\LL$-modules
in $\md(\LL,d)$ and $\LL_d$, and $KQ$-modules in $\rep(Q,d)$,
respectively, correspond to the orbits of these actions.
For a module $M$ in $\md(\LL,d)$, (resp. in $\LL_d$ or in $\rep(Q,d)$), let
$\GL_d.M$ denote its $\GL_d$-orbit.

There is a bilinear form $\bil{-,-} = \bil{-,-}_Q\df \Z^n \times \Z^n \to \Z$
associated to $Q$ defined by
$$
\bil{d,e} := \bil{d,e}_Q := \sum_{i \in Q_0} d_ie_i -
\sum_{a \in Q_1} d_{s(a)}e_{t(a)}.
$$
The dimension vector of a $KQ$-module $M$ is denoted by 
$\dimv(M) = \dimv_Q(M)$.
(Note that $\dimv_Q(M) = \dimv_\LL(M)$, since we can consider $M$ also as
a $\LL$-module.)
For $KQ$-modules $M$ and $N$ set
$$
\bil{M,N} := \bil{M,N}_Q :=
\dm \Hom_{KQ}(M,N) - \dm \Ext_{KQ}^1(M,N).
$$
It is known that $\bil{M,N} = \bil{\dimv(M),\dimv(N)}$.
Let $(-,-) = (-,-)_Q\df \Z^n \times \Z^n \to \Z$
be the symmetrization of the bilinear form $\bil{-,-}$,
\textit{i.e.} $(d,e) := \bil{d,e} + \bil{e,d}$.
For $\LL$-modules $X$ and $Y$ set
$$
(X,Y)_Q := \bil{\pi_Q(X),\pi_Q(Y)}_Q + 
\bil{\pi_Q(Y),\pi_Q(X)}_Q.
$$

\begin{Lem}[{\cite[Lemma 1]{CB}}]\label{extandhom}
For any $\LL$-modules $X$ and $Y$ we have
$$
\dm \Ext_\LL^1(X,Y) = \dm \Hom_\LL(X,Y) + \dm \Hom_\LL(Y,X) - 
(X,Y)_Q.
$$
In particular,
$\dm \Ext_\LL^1(X,X)$ is even, and 
$\dm \Ext_\LL^1(X,Y) = \dm \Ext_\LL^1(Y,X)$.
\end{Lem}

\begin{Cor}
For a nilpotent $\LL$-module $X$ with dimension vector
$d$ the following are equivalent:
\begin{itemize}

\item
The closure $\overline{\GL_d.X}$ of $\GL_d.X$ is an irreducible
component of $\LL_d$;

\item
The orbit $\GL_d.X$ is open in $\LL_d$;

\item
$\Ext_\LL^1(X,X) = 0$.

\end{itemize}
\end{Cor}

\subsection{Semicanonical bases}\label{section_semican}
We recall the definition of the dual semicanonical basis and its
multiplicative properties, following \cite{Lu1, Lu2, GLSSemi1,GLSSemi2}.
From now on, assume that $K = \C$.

For each dimension vector $d = (d_1,\ldots,d_n)$
we defined the affine variety $\LL_d$.
A subset $C$ of $\LL_d$ is said to be constructible if it is
a finite union of locally closed subsets.
For a $\C$-vector space $V$, a function 
$$
f\df \LL_d \to V
$$ 
is {\it constructible} 
if 
the image $f(\LL_d)$ is finite and 
$f^{-1}(m)$ is a constructible
subset of $\LL_d$ for all $m \in V$.

The set of constructible functions  $\LL_d\to\C$ is denoted by 
$M(\LL_d)$. This is a $\C$-vector space. 
Recall that the group $\GL_d$ acts on $\LL_d$ by conjugation.
By $M(\LL_d)^{\GL_d}$ 
we denote the subspace of
$M(\LL_d)$ consisting of the constructible functions which are
constant on the $\GL_d$-orbits in $\LL_d$.
Set
$$
\widetilde{\MM} := \bigoplus_{d \in \N^n} M(\LL_d)^{\GL_d}.
$$
For $f' \in M(\LL_{d'})^{\GL_{d'}}$, $f'' \in M(\LL_{d''})^{\GL_{d''}}$
and $d = d'+d''$
we define a constructible function 
$$
f := f' \star f''\df \LL_d \to \C
$$ 
in $M(\LL_d)^{\GL_d}$ by
$$
f(X) := \sum_{m \in \C} m \, \chi_{\rm c}\left(\left\{ 
U \subseteq X \mid f'(U)f''(X/U) = m
\right\} \right)
$$
for all $X \in \LL_d$, where $U$ runs over 
the points of the Grassmannian of all submodules
of $X$ with $\dimv(U) = d'$.
Here, for a constructible subset $V$ of a complex variety 
we denote by $\chi_{\rm c}(V)$ its
(topological) Euler characteristic with respect to cohomology with compact 
support.
This turns $\widetilde{\MM}$ into an associative $\C$-algebra.

\begin{Rem}\label{newconv}
Note that the product $\star$ defined here is opposite to 
the convolution product we have used in 
\cite{GLSSemi1,GLSVerma,GLSSemi2}.
This new convention turns out to be better adapted to our choice
of categorifying $\C[N(w)]$ and $\C[N^w]$ by categories
closed under factor modules.
It is also compatible with our choice in \cite{GLSFlag}
of categorifying coordinate rings of partial flag varieties 
by categories closed under submodules.
\end{Rem}

For the canonical basis vector $e_i := \dimv(S_i)$ we know
that $\LL_{e_i}$ is just a point, which (as a $\LL$-module)
is isomorphic to the simple module $S_i$.
Define ${\bf 1}_i\df \LL_{e_i} \to \C$ by ${\bf 1}_i(S_i) := 1$.
By $\MM$ we denote the subalgebra of $\widetilde{\MM}$ generated
by the functions ${\bf 1}_i$ where $1 \le i \le n$.
Set $\MM_d := \MM \cap M(\LL_d)^{\GL_d}$.
It follows that
$$
\MM = \bigoplus_{d \in \N^n} \MM_d
$$
is an $\N^n$-graded $\C$-algebra.
Let $U(\n)$ be the enveloping algebra of the positive
part $\n$ of the Kac-Moody Lie algebra $\g$ associated to $Q$,
see Sections~\ref{KMalg} and \ref{universal}.

\begin{Thm}[{Lusztig \cite{Lu2}}]\label{lusztig1}
There is an isomorphism of $\N^n$-graded $\C$-algebras
$$
U(\n) \to \MM
$$
defined by $E_i \mapsto {\bf 1}_i$ for $1 \le i \le n$.
\end{Thm}

Let $\irr(\LL_d)$ be the set of irreducible components
of $\LL_d$.

\begin{Thm}[{Lusztig \cite{Lu2}}]\label{lusztig2}
For each $Z \in \irr(\LL_d)$ there
is a unique $f_Z\df \LL_d \to \C$ in $\MM_d$
such that $f_Z$ takes value 1 on some dense open subset of
$Z$ and value 0 on some dense open subset of any other irreducible
component $Z'$ of $\LL_d$.
Furthermore, the set 
$$
\cS := \left\{ f_Z \mid Z \in \irr(\LL_d), d \in \N^n \right\}
$$
is a $\C$-basis of $\MM$.
\end{Thm}

The basis $\cS$ is called the {\it semicanonical basis} 
of $\MM$.
By Theorem~\ref{lusztig1} we just identify $\MM$ and $U(\n)$ and
consider $\cS$ also as a basis of $U(\n)$.
Since $U(\n)$ is a cocommutative Hopf algebra, its graded dual 
$$
U(\n)^*_{\rm gr} = \bigoplus_{d \in \N^n} U_d^*
$$
is a commutative  $\C$-algebra.
Let $\MM_d^*$ be the dual space of $\MM_d$, and set
$$
\MM^* := \bigoplus_{d \in \N^n} \MM_d^*.
$$
Again we identify $\MM^*$ and $U(\n)^*_{\rm gr}$.

For $X \in \LL_d$ define an {\it evaluation function}
$$
\delta_X\df \MM_d \to \C
$$
by $\delta_X(f) := f(X)$.

It is not difficult to show that
the map $X \mapsto \delta_X$ from $\LL_d$ to $\MM_d^*$ is constructible
in the above sense.
So on every irreducible component 
$Z \in \irr(\LL_d)$ there is a Zariski open set on which
this map is constant.
Define $\rho_Z := \delta_X$ for any $X$ in this open set.
The $\C$-vector space $\MM_d^*$ 
is spanned by the functions
$\delta_X$ with $X \in \LL_d$.
Then by construction
$$
\cS^* := \left\{ \rho_Z \mid Z \in \irr(\LL_d), d \in \N^n \right\}
$$
is the basis of $\MM^* \equiv U(\n)^*_{\rm gr}$ 
dual to Lusztig's semicanonical basis 
$\cS$
of $U(\n)$.

In case $X$ is a rigid $\LL$-module, the orbit of $X$ in $\LL_d$ is open,
its closure is an irreducible component $Z$, and $\delta_X=\rho_Z$ belongs to
$\cS^*$.

For a module $X \in \LL_d$ and an $m$-tuple 
$\ii = (i_1,\ldots,i_m)$ with $1 \le i_j \le n$ for
all $j$, let
$\F_{\ii,X}$ denote the projective variety of composition series of type
$\ii$ of $X$.
Thus an element in 
$\F_{\ii,X}$ 
is a flag
$$
(0 = X_0 \subset X_1 \subset \cdots \subset X_m = X)
$$
of submodules $X_j$ of $X$
such that for all $1 \le j \le m$ the subfactor 
$X_j/X_{j-1}$ is isomorphic to the simple $\LL$-module $S_{i_j}$
associated to the vertex $i_j$ of $Q$.
Let 
$$
d_\ii\df \LL_d \to \C
$$
be the map which
sends $X \in \LL_d$ to $\chi_{\rm c}(\F_{\ii,X})$.
It follows from the definition of $\star$ that
$d_\ii = {\bf 1}_{i_1} \star \cdots \star {\bf 1}_{i_m}$.
The $\C$-vector space
$\MM_d$ is spanned by the maps
$d_\ii$.
We have $\delta_X(d_\ii) = \chi_{\rm c}(\F_{\ii,X})$.

\begin{Thm}[{\cite{GLSSemi1}}]\label{2.6i}
For $X,Y \in \nil(\LL)$ we have
$\delta_X\delta_Y = \delta_{X\oplus Y}$. 
\end{Thm}

In \cite{GLSSemi2} a more complicated formula than the one in
Theorem~\ref{2.6i}
is given, expressing
$\delta_X\delta_Y$ as a linear combination of $\delta_Z$ where $Z$ runs over 
all possible middle terms of non-split short exact sequences with end terms 
$X$ and $Y$.
The formula is especially useful when 
$\dm\Ext_\LL^1(X,Y)=1$.
In this case, the following hold:

\begin{Thm}[{\cite[Theorem 2]{GLSSemi2}}]\label{2.6ii}
Let $X,Y \in \nil(\LL)$.
If $\dm \Ext_\LL^1(X,Y) = 1$ with
\[
0 \to X\to E'\to Y\to 0 \text{\;\;\; and \;\;\;} 0\to Y\to E''\to X\to 0
\]
the corresponding non-split short exact sequences, then
$$
\delta_X\delta_Y = \delta_{E'} + \delta_{E''}.
$$
\end{Thm}

\subsection{Frobenius categories}
Let $A$ be a $K$-algebra.
Let $\CC$ be a subcategory of a module category $\md(A)$ 
which is closed under extensions.
Clearly, we have
$$
\Ext_{\CC}^1(X,Y) = \Ext_A^1(X,Y)
$$ 
for all modules $X$ and $Y$ in $\CC$.
An $A$-module $C$ in $\CC$ is called $\CC$-{\it projective}  
(resp. $\CC$-{\it injective})
if $\Ext_A^1(C,X) = 0$ (resp. $\Ext_A^1(X,C) = 0$) for all $X \in \CC$.
If $C$ is $\CC$-projective and $\CC$-injective, then $C$ is also called
$\CC$-{\it projective-injective}.
We say that $\CC$ has {\it enough projectives} 
(resp. {\it enough injectives})
if for each $X \in \CC$ there exists a short exact sequence
$0 \to Y \to C \to X \to 0$
(resp. $0 \to X \to C \to Y \to 0$)
where $C$ is $\CC$-projective (resp. $\CC$-injective)
and $Y \in \CC$.
If $\CC$ has enough projectives and enough injectives, and
if these coincide (\textit{i.e.}~an object is $\CC$-projective if and only if
it is $\CC$-injective), then
$\CC$ is called a {\it Frobenius subcategory} 
of $\md(A)$.
In particular, $\CC$ is a Frobenius category in the sense of
Happel \cite{H2}.
Of course,
for $A = \LL$, an $A$-module $C$ in $\CC$ is
$\CC$-projective if and only if it is 
$\CC$-injective, see Lemma~\ref{extandhom}.

By definition the objects in the {\it stable category} ${\stCC}$ are the 
same as the objects in $\CC$, and the morphism spaces 
$\Hom_{\stCC}(X,Y)$ are the morphism
spaces in $\CC$ modulo morphisms factoring through 
$\CC$-projective-injective objects.
The category ${\stCC}$ is a triangulated 
category in a natural way \cite{H2}, where
the shift is given by the relative inverse syzygy functor
$$
\Omega^{-1}\df {\stCC} \to {\stCC}.
$$
For all $X$ and $Y$ in $\CC$ there
is a functorial isomorphism
$$
\Hom_{\stCC}(X,\Omega^{-1}(Y)) \cong \Ext_\CC^1(X,Y).
$$
The category $\stCC$ is a {\it 2-Calabi-Yau category}, if
for all $X,Y \in \stCC$ there is
a functorial isomorphism
$$
\Ext_\CC^1(X,Y) \cong D\Ext_\CC^1(Y,X).
$$

\subsection{Frobenius categories associated to Weyl group elements}
\label{terminal}
By $\hI_1,\ldots,\hI_n$ we denote the indecomposable injective $\LL$-modules 
with socle $S_1,\ldots,S_n$, respectively.
Here $S_1,\ldots,S_n$ are the 1-dimensional simple $\LL$-modules
corresponding to the vertices of the quiver $Q$.
(The modules $\hI_i$ are infinite-dimensional if $Q$ is not
a Dynkin quiver.)

For a $\LL$-module $X$ and a simple module $S_j$ let
$\soc_{(j)}(X) := \soc_{S_j}(X)$ be the sum of all submodules
$U$ of $X$ with $U \cong S_j$.
(In this definition, we do not assume that $X$ is finite-dimensional.)
For a sequence $(j_1,\ldots,j_t)$ of indices with
$1 \le j_p \le n$ for all $p$, there is a unique chain
$$
0 = X_0 \subseteq X_1 \subseteq \cdots \subseteq X_t \subseteq X
$$
of submodules of $X$ such that $X_p/X_{p-1} = \soc_{(j_p)}(X/X_{p-1})$.
Define $\soc_{(j_1,\ldots,j_t)}(X) := X_t$.

For the rest of this section, 
let $\ii = (i_r,\ldots,i_1)$ be a reduced expression of an
element $w$ of the Weyl group $W = W_Q$ of $Q$.
(By definition, this is the Weyl group of the Kac-Moody Lie
algebra $\g$ associated to $Q$, see Section~\ref{KMalg}.)
For $1 \le k \le r$ let 
$$
V_k := V_{\ii,k} := \soc_{(i_k,\ldots,i_1)}\left(\hI_{i_k}\right), 
$$
and set
$V_\ii := V_1 \oplus \cdots \oplus V_r$.
(The module $V_\ii$ is dual to the cluster-tilting object
constructed in \cite[Section III.2]{BIRS}.)
Define
$$
\CC_\ii := \Gen(V_\ii) \subseteq \nil(\LL).
$$
For $1 \le j \le n$ let $k_j := \max\{ 1 \le k \le r \mid i_k = j \}$.
Define 
$I_{\ii,j} := V_{\ii,k_j}$
and set
$$
I_\ii := I_{\ii,1} \oplus \cdots \oplus I_{\ii,n}.
$$
The category $\CC_\ii$ and the module $I_\ii$ depend only on $w$, and not
on the chosen reduced expression $\ii$ of $w$.
Therefore, we define 
$$
\CC_w := \CC_\ii \text{ and } I_w := I_\ii.
$$
(If $Q$ is a Dynkin quiver, and $w = w_0$ is the longest Weyl group
element, then $\CC_w = \nil(\LL) = \md(\LL)$.)
Without loss of generality, we assume that for each $1 \le j \le n$
there is some $1 \le k \le r$ with $i_k = j$.
Otherwise, we could just replace $Q$ by a quiver with fewer vertices.
Note also that $\CC_w = \add(I_w)$ if and only if $i_k \not= i_s$
for all $k \not= s$.
In this case, most of our theory becomes trivial.

The following three theorems are proved in \cite{BIRS}.
They were also obtained independently and by different methods
in \cite{GLSUni1} in the case when $w$ is adaptable.

\begin{Thm}\label{main1}
For any Weyl group element $w$
the following hold:
\begin{itemize}

\item[(i)]
$\CC_w$ is a Frobenius category;

\item[(ii)]
The stable category ${\stCC}_w$ is a 2-Calabi-Yau category;

\item[(iii)]
$\CC_w$ has $n$ indecomposable $\CC_w$-projective-injective modules,
namely the indecomposable direct summands of $I_w$;

\item[(iv)]
$\CC_w = \Gen(I_w)$.

\end{itemize}
\end{Thm}

We denote the relative inverse syzygy functor of $\stCC_w$ by
$\Omega_w^{-1}$.

Recall that a $\LL$-module $T$  is 
{\it rigid} 
if $\Ext_\LL^1(T,T) = 0$.
Let $\CC$ be a subcategory of $\md(\LL)$, and
let $T \in \CC$ be rigid.
Recall that for all $X,Y \in \md(\LL)$ we have
$\dm \Ext_\LL^1(X,Y) = \dm \Ext_\LL^1(Y,X)$.
We need the following definitions:
\begin{itemize}

\item
$T$ is $\CC$-{\it maximal rigid} if 
$
\Ext_\LL^1(T \oplus X,X) = 0
$
with $X \in \CC$ implies $X \in \add(T)$;

\item
$T$ is a $\CC$-{\it cluster-tilting module} if 
$
\Ext_\LL^1(T,X) = 0
$
with $X \in \CC$
implies $X \in \add(T)$.

\end{itemize}

\begin{Thm}\label{main9}
For a rigid $\LL$-module $T$ in $\CC_w$ the following are equivalent:
\begin{itemize}

\item[(i)]
$\Sigma(T) = {\rm length}(w)$;

\item[(ii)]
$T$ is $\CC_w$-maximal rigid;

\item[(iii)]
$T$ is a $\CC_w$-cluster-tilting module.

\end{itemize}
\end{Thm}

For $1 \le k \le r$ let 
\begin{align*}
k^- &:= \max\{0,1 \le s \le k-1 \mid i_s = i_k \},\\
k^+ &:= \min\{ k+1 \le s \le r,r+1 \mid i_s = i_k \}.
\end{align*}
For $1 \le i,j \le n$ let $q_{ij}$ be the number of edges between 
the vertices $i$ and $j$ of the underlying graph of our quiver $Q$.

Following Berenstein, Fomin and Zelevinsky we define a quiver
$\GG_\ii$ as follows:
The vertices of $\GG_\ii$ are just the numbers $1,\ldots,r$.
For $1 \le s,t \le r$ there are $q_{i_s,i_t}$ arrows from $s$ to $t$ provided
$t^+ \ge s^+ > t > s$.
These are called the {\it ordinary arrows} of $\GG_\ii$.
Furthermore, for each $1 \le s \le r$ there is an arrow $s \to s^-$
provided $s^- > 0$.
These are the {\it horizontal arrows} of $\GG_\ii$. 

On the other hand, let $A$ be a $K$-algebra, and
let $X = X_1^{n_1} \oplus \cdots \oplus X_t^{n_t}$ 
be a finite-dimensional $A$-module, where the $X_i$
are pairwise non-isomorphic indecomposable
modules and $n_i \ge 1$.
Let $S_i = S_{X_i}$ be the simple $\End_A(X)^\op$-module
corresponding to $X_i$.
Then
$\Hom_A(X,X_i)$ is the indecomposable projective $\End_A(X)^\op$-module with
top $S_i$.
The basic facts on the quiver $\GG_X$ of the endomorphism algebra
$\End_A(X)^\op$ are collected in \cite[Section 3.2]{GLSRigid}.
In particular, we have a 1-1 correspondence between the vertices
of $\GG_X$ and the modules $X_1,\ldots,X_t$.

\begin{Thm}\label{main6}
The module $V_\ii$ is $\CC_w$-maximal rigid, 
and we have $\GG_{V_\ii} = \GG_\ii$.
\end{Thm}

For example, 
let $Q$ be a quiver with underlying graph 
$\xymatrix@-0.5pc{1 \ar@<0.5ex>@{-}[r]\ar@<-0.5ex>@{-}[r] & 2 \ar@{-}[r]& 3}$.
Then $\ii := (i_7,\ldots,i_1) := (3,1,2,3,1,2,1)$ is a reduced expression
of a Weyl group element $w \in W_Q$.
The quiver $\GG_\ii$ looks as follows:
$$
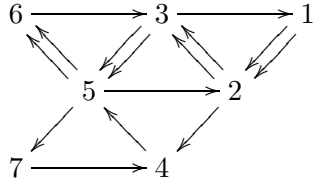
\xymatrix@-0.7pc{ 
6 \ar[rr] &&
3 \ar[rr] \ar@<0.5ex>[dl]\ar@<-0.5ex>[dl] && 1
\ar@<0.5ex>[dl]\ar@<-0.5ex>[dl]
\\
&5 \ar[rr] \ar[dl]\ar@<0.5ex>[ul]\ar@<-0.5ex>[ul] && 
2 \ar[dl]\ar@<0.5ex>[ul]\ar@<-0.5ex>[ul]
\\
7 \ar[rr] && 4 \ar[ul]
}
$$

We often try to visualize $\LL$-modules.
For example, let $Q$ be the quiver
$$
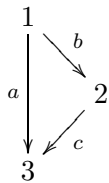
\xymatrix@-0.7pc{
1 \ar[dd]_a\ar[dr]^b\\
& 2 \ar[dl]^c\\
3
}
$$
and let $\ii := (i_6,\ldots,i_1) := (3,2,1,3,2,1)$.
Then the $\LL$-module $V_\ii = V_1 \oplus \cdots \oplus V_6$
looks as follows:
\begin{align*}
V_1 &= {\bsm1\esm} &
V_2 &= {\bsm&1\\2\esm} &
V_3 &= {\bsm&&&1\\1&&2\\&3\esm} \\
V_4 &= {\bsm&&&&1\\&1&&2\\2&&3\\&1\esm} &
V_5 &= {\bsm&&&&&&1\\&&&1&&2\\1&&2&&3\\&3&&1\\&&2\esm} &
V_6 &= {\bsm&&&&&&&1\\&&&&1&&2\\&1&&2&&3\\2&&3&&1\\&1&&2\\&&3\esm}
\end{align*}
The numbers can be interpreted as basis vectors or
as composition factors.
For example, the module $V_5$ is
a 9-dimensional $\LL$-module
with dimension vector $\dimv_\LL(V_5) = (d_1,d_2,d_3) = (4,3,2)$.
More precisely, one could display $V_5$ as follows:
$$
\xymatrix@-0.7pc{
&&&&&&1\ar[dl]^b\\
&&&1\ar[dl]_b\ar[dr]^a&&2\ar[dl]^c\\
1\ar[dr]_a && 2\ar[dr]^{b^*}\ar[dl]_c&&3\ar[dl]^{a^*}\\
&3\ar[dr]_{c^*}&&1\ar[dl]^b\\
&&2
}
$$
This picture shows how the different arrows
of the quiver $\overline{Q}$ of $\LL$ act
on the 9 basis vectors of $V_5$.
For example, one can see immediately that
the socle of $X$ is isomorphic to $S_2$, and the
top is isomorphic to $S_1 \oplus S_1 \oplus S_1$.

\subsection{Relative homology for $\CC_w$}
We recall some notions from relative homology theory
which, for Artin algebras, was developed by Auslander and Solberg
\cite{AS1,AS2}.

Let $A$ be a $K$-algebra, and let $X,Y,Z,T \in \md(A)$.
Set 
$$ 
F_T := \Hom_A(T,-)\df \md(A) \to \md(\End_A(T)^\op).
$$
A short exact sequence
$$
0 \to Z \to Y \to X \to 0
$$
is $F_T$-{\it exact} if
$0 \to F_T(Z) \to F_T(Y) \to F_T(X) \to 0$
is exact.
By $F_T(X,Z)$ we denote the set of equivalence classes of
$F_T$-exact sequences with end terms $X$ and $Z$ as above.

Let $\Y_T$ be the subcategory of all
$X \in \md(A)$ such that
there exists an exact sequence
\begin{equation}\label{seq7}
\cdots \xrightarrow{f_3} T_3 \xrightarrow{f_2} T_1 
\xrightarrow{f_1} T_0 \xrightarrow{f_0} X \to 0
\end{equation}
where $T_i \in \add(T)$ for all $i$ and the short exact sequences
$$
0 \to \Ker(f_i) \to T_i \to \Ima(f_i) \to 0
$$
are $F_T$-exact for all $i \ge 0$.
We call sequence~(\ref{seq7}) an $\add(T)$-{\it resolution} of $X$.
We say that (\ref{seq7}) has {\it length} at most $d$ if $T_j = 0$
for all $j > d$.
Note that
$$
\add(T) \subseteq \Y_T.
$$
Dually, one defines $\add(T)$-{\it coresolutions} 
$$
0 \to X \xrightarrow{g_0} T_0 \xrightarrow{g_1} T_1 
\xrightarrow{g_2} T_2 \xrightarrow{g_3} \cdots
$$
where we require now that the sequences
$$
0 \to \Ima(g_i) \to T_i \to \Coker(g_i) \to 0
$$
are $F^T$-exact, where $F^T$ is the contravariant functor $\Hom_A(-,T)$.

For $X \in \Y_T$ and $Z \in \md(A)$ let
$\Ext_{F_T}^i(X,Z)$, $i \ge 0$ 
be the cohomology groups of the cocomplex obtained by
applying the functor $\Hom_A(-,Z)$ to the sequence
$$
\cdots \xrightarrow{f_3} T_2 
\xrightarrow{f_2} T_1 \xrightarrow{f_1} T_0.
$$

\begin{Lem}[{\cite{AS1}}]
For $X \in \Y_T$ and $Z \in \md(A)$ there is a functorial
isomorphism 
$$
\Ext_{F_T}^1(X,Z) \cong F_T(X,Z).
$$
\end{Lem}

\begin{Prop}[{\cite[Proposition 3.7]{AS2}}]\label{seslift}
For $X \in \Y_T$ and $Z \in \md(A)$ 
there is a functorial isomorphism
$$
\Ext_{F_T}^i(X,Z) \to \Ext_{\End_A(T)^\op}^i(\Hom_A(T,X),\Hom_A(T,Z)) 
$$
for all $i \ge 0$.
\end{Prop}

\begin{Cor}\label{fullyfaithful}
The functor 
$$
\Hom_A(T,-)\df \Y_T \to \md(\End_A(T)^\op)
$$
is fully faithful.
In particular, $\Hom_A(T,-)$ has the following properties:
\begin{itemize}

\item[(i)]
If $X \in \Y_T$ is indecomposable, then $\Hom_A(T,X)$ is indecomposable;

\item[(ii)]
If $\Hom_A(T,X) \cong \Hom_A(T,Y)$ for some $X,Y \in \Y_T$, 
then $X \cong Y$.

\end{itemize}
\end{Cor}

Note that Corollary~\ref{fullyfaithful} follows already from
\cite[Section 3]{A}, see also 
\cite[Lemma 1.3 (b)]{APR}.

\begin{Cor}\label{fullyfaithful2}
Let $T \in \md(A)$, and let $\CC$ be an extension closed subcategory
of $\Y_T$.
If 
$$
\psi\df 0 \to \Hom_A(T,X) \xrightarrow{\Hom_A(T,f)} \Hom_A(T,Y)
\xrightarrow{\Hom_A(T,g)} \Hom_A(T,Z) \to 0
$$
is a short exact sequence of $\End_A(T)^\op$-modules with $X,Y,Z \in \CC$,
then 
$$
\eta\df 0 \to X \xrightarrow{f} Y \xrightarrow{g} Z \to 0 
$$
is a short exact sequence in $\md(A)$.
\end{Cor}

Now we apply the above ideas to the category $\CC_w$.
The following proposition is proved in \cite{GLSUni1} for adaptable $w$
and in \cite{BIRS} for arbitrary $w$.
In a more general framework it is proved in \cite{KR}.

\begin{Prop}\label{cormutation1}
Let $T$ be a $\CC_w$-maximal rigid module, and let $X \in \CC_w$.
Then there exists an $\add(T)$-resolution of the form
$$
0 \to T_1 \to T_0 \to X \to 0
$$
and an $\add(T)$-coresolution of the form
$$
0 \to X \to T_0' \to T_1' \to 0.
$$
\end{Prop}

\begin{Cor}\label{fullyfaithful3}
For each $\CC_w$-maximal rigid module $T$ we have
$\CC_w \subseteq \Y_T$.
\end{Cor}

\begin{Cor}\label{pdimone}
For each $X \in \CC_w$
the projective dimension of the $\End_\LL(T)^\op$-module
$\Hom_\LL(T,X)$ is at most one.
\end{Cor}

\begin{Cor}[{\cite[Theorem 5.3.2]{Iy}}]\label{tilt}
If $T$ and $R$ are $\CC_w$-maximal rigid $\LL$-modules,
then the $\End_\LL(T)^\op$-module $\Hom_\LL(T,R)$ is a
classical tilting module, and
$$
\End_{\End_\LL(T)^\op}(\Hom_\LL(T,R)) \cong \End_\LL(R).
$$
\end{Cor}

\subsection{The cluster algebra $\cA(\CC_w,T)$}\label{clustintro}
We refer to \cite{FZSurv} for an excellent survey on cluster algebras.  
Here we only recall the main definitions and introduce 
a cluster algebra $\cA(\CC_w,T)$ associated to
a Weyl group element $w$ and a $\CC_w$-maximal rigid $\LL$-module $T$.

If $\widetilde{B} = (b_{ij})$ is any $r \times (r-n)$-matrix
with integer entries, then
the {\it principal part} $B$ of 
$\widetilde{B}$ is obtained from 
$\widetilde{B}$ by deleting the last $n$ rows.
Given some $1 \le k \le r-n$ define a new $r \times (r-n)$-matrix 
$\mu_k(\widetilde{B}) = (b_{ij}')$ by
$$
b_{ij}' =
\begin{cases}
-b_{ij} & \text{if $i=k$ or $j=k$},\\
b_{ij} + \dfrac{|b_{ik}|b_{kj} + b_{ik}|b_{kj}|}{2} & \text{otherwise},
\end{cases}
$$
where $1 \le i \le r$ and $1 \le j \le r-n$.
One calls $\mu_k(\widetilde{B})$ a {\it mutation} 
of $\widetilde{B}$.
If $\widetilde{B}$ is an integer matrix whose principal part is
skew-symmetric, then it is 
easy to check that $\mu_k(\widetilde{B})$ is also an integer matrix 
with skew-symmetric principal part.
In this case, Fomin and Zelevinsky define a cluster algebra
$\cA(\widetilde{B})$ as follows.
Let $\F = \C(y_1,\ldots,y_r)$ be the field of rational
functions in $r$ commuting variables $y_1,\ldots,y_r$.
Define
$\yy := (y_1,\ldots,y_r)$. 
One calls $(\yy,\widetilde{B})$ the {\it initial seed} 
of
$\cA(\widetilde{B})$.
For $1 \le k \le r-n$ define 
$$
y_k^* := 
\frac{\prod_{b_{ik}> 0} y_i^{b_{ik}} + \prod_{b_{ik}< 0} y_i^{-b_{ik}}}{y_k}.
$$
The pair 
$
(\mu_k(\yy),\mu_k(\widetilde{B}))
$, 
where
$\mu_k(\yy)$ is obtained from $\yy$ by replacing $y_k$ by 
$y_k^*$,
is the 
{\it mutation in direction} 
$k$ of the seed
$(\yy,\widetilde{B})$. 

Now one can iterate this process of mutation and obtain inductively
a set of seeds.
Thus
each seed consists of an $r$-tuple of algebraically independent 
elements of $\F$
called a 
{\it cluster} 
and of a matrix called the 
{\it exchange matrix}.
The elements of a cluster are its 
{\it  cluster variables}.
Given a cluster $(f_1,\ldots,f_r)$, the monomials 
$f_1^{m_1}f_2^{m_2} \cdots f_r^{m_r}$ where $m_k \ge 0$ for all
$k$ are called {\it cluster monomials}.
A seed has $r-n$ neighbours obtained by mutation in direction
$1 \le k \le r-n$.
One does not mutate the last $n$ elements of a cluster, they serve
as ''coefficients'' and belong to every cluster.
The {\it cluster algebra}
$
\cA(\widetilde{B})
$
is by definition the subalgebra of $\F$ generated by the
set of all cluster variables appearing in all seeds obtained 
by iterated mutation starting with the initial seed. 

It is often convenient to define a cluster algebra using an
oriented graph, as follows.
Let $\GG$ be a quiver without loops or 2-cycles
with vertices $\{ 1,\ldots, r\}$.
We can define an $r \times r$-matrix $B(\GG) = (b_{ij})$ 
by setting
$$
b_{ij} = 
(\text{number of arrows $j \to i$ in $\GG$}) 
-(\text{number of arrows $i \to j$ in $\GG$}).
$$
Let $B(\GG)^\circ$ 
be the $r \times (r-n)$-matrix obtained by deleting
the last $n$ columns of $B(\GG)$.
The principal part of $B(\GG)^\circ$ is skew-symmetric, hence
this yields a cluster algebra 
$\cA(B(\GG)^\circ)$.

We apply this procedure to our subcategory $\CC_w$.
Let $T = T_1 \oplus \cdots \oplus T_r$ be a basic 
$\CC_w$-maximal rigid $\LL$-module
with $T_k$ indecomposable for all $k$.
Without loss of generality assume that $T_{r-n+1}, \ldots, T_r$ are
$\CC_w$-projective-injective. 
By $\GG_T$ 
we denote the quiver of the endomorphism algebra $\End_\LL(T)^\op$.
We then define the cluster algebra
$$
\cA(\CC_w,T) := \cA(B(\GG_T)^\circ).
$$
In particular, we denote by
$\cA(\CC_w)$ the cluster algebra 
$\cA(\CC_w,V_\ii)$ attached to the $\CC_w$-maximal rigid
module $V_\ii$ of Section~\ref{terminal}. 
Thus $\cA(\CC_w) := \cA(B(\GG_\ii)^\circ)$.
(Up to isomorphism of cluster algebras, this definition
does not depend on the choice of $\ii$, see Section~\ref{reach}.)

\subsection{Mutation of  rigid modules}\label{intro1.6}
The results of this section are straightforward generalizations of 
results in \cite{GLSRigid}, see
\cite[Sections 12,13,14]{GLSUni1} and \cite{BIRS}.

Let $A$ be a $K$-algebra, and $M$ be an $A$-module.
A homomorphism $f\df X \to M'$ in $\md(A)$ is 
a 
{\it left $\add(M)$-approximation} 
of $X$ if $M' \in \add(M)$ 
and the induced map
$$
\Hom_A(f,M)\df \Hom_A(M',M) \to \Hom_A(X,M) 
$$
is surjective.
A morphism $f\df V \to W$ is called
{\it left minimal} 
if every morphism $g\df W \to W$ with $gf = f$ is an 
isomorphism.
Dually, one defines right $\add(M)$-approximations and
right minimal morphisms.
Some well known basic properties of approximations can be found in 
\cite[Section 3.1]{GLSRigid}.

\begin{Prop}\label{corquivershape}
Let $T$ be a basic $\CC_w$-maximal rigid $\LL$-module, and let $X$ be
an indecomposable direct summand of $T$ which is not $\CC_w$-projective-injective.
Then there are short exact sequences
$$
0 \to X \xrightarrow{f'} T' \xrightarrow{g'} Y \to 0
$$
and
$$
0 \to Y \xrightarrow{f''} T'' \xrightarrow{g''} X \to 0
$$
such that the following hold:
\begin{itemize}

\item[(i)]
$f'$ and $f''$ are minimal left $\add(T/X)$-approximations, and
$g'$ and $g''$ are minimal right $\add(T/X)$-approximations;

\item[(ii)]
$Y \oplus T/X$ is a basic $\CC_w$-maximal rigid $\LL$-module
(in particular $Y$ is indecomposable), and
$X \not\cong Y$;

\item[(iii)]
$\dm \Ext_\LL^1(Y,X) = \dm \Ext_\LL^1(X,Y) = 1$;

\item[(iv)]
We have $\add(T') \cap \add(T'') = 0$;

\item[(v)]
The quiver $\GG_T$ of $\End_\LL(T)^\op$ has no loops and no 2-cycles;

\item[(vi)]
We have
$$
\gldim(\End_\LL(T)^\op) =
\begin{cases}
3 & \text{$\CC_w \not= \add(I_w)$},\\ 
1 & \text{$\CC_w = \add(I_w)$ and $n > 1$},\\
0 & \text{$\CC_w = \add(I_w)$ and $n = 1$}.
\end{cases}
$$

\end{itemize}
\end{Prop}

In the situation of the above proposition, we call $\{ X,Y \}$ an
{\it exchange pair associated to} $T/X$, and we write
$$
\mu_X(T) = Y \oplus T/X.
$$
We say that $Y \oplus T/X$ is the mutation of $T$ in direction $X$.
The short exact sequence
\[
0 \to X \xrightarrow{f'} T' \xrightarrow{g'} Y \to 0
\]
is the {\it exchange sequence} starting in $X$ and ending in $Y$.
Thus, we have
$$
\mu_Y(\mu_X(T)) = T.
$$

Let $T = T_1 \oplus \cdots \oplus T_r$ be a basic $\CC_w$-maximal rigid
$\LL$-module with $T_k$ indecomposable for all $k$.
Without loss of generality we assume that 
$T_{r-n+1}, \ldots, T_r$ are $\CC_w$-projective-injective.
As in Section~\ref{clustintro}
write $B(T) := B(\GG_T) = (t_{ij})_{1 \le i,j \le r}$, and let
$B(T)^\circ = (t_{ij})$ 
be the $r \times (r-n)$-matrix obtained from
$B(T)$ by deleting the last $n$ columns.

For $1 \le k \le r-n$ let
$$
0 \to T_k \to T' \to T_k^* \to 0
$$
and
$$
0 \to T_k^* \to T'' \to T_k \to 0
$$
be exchange sequences associated to the direct summand $T_k$
of $T$.
It follows that
$$
T'  = \bigoplus_{t_{ik} < 0} T_i^{-t_{ik}}
\;\;\;\; \text{ and }  \;\;\;\;
T'' = \bigoplus_{t_{ik} > 0} T_i^{t_{ik}}.
$$
Set 
$$
T^* = \mu_{T_k}(T) = T_k^* \oplus T/T_k.
$$
The quivers of the endomorphism
algebras $\End_\LL(T)^\op$ and $\End_\LL(\mu_{T_k}(T))^\op$ are related
via Fomin and Zelevinsky's mutation rule:

\begin{Thm}\label{main4}
Let $w$ be a Weyl group element.
For a basic $\CC_w$-maximal rigid $\LL$-module $T$ as above and 
$1 \le k \le r-n$
we have
$$
B(\mu_{T_k}(T))^\circ = \mu_k(B(T)^\circ). 
$$
\end{Thm}

\subsection{Categorification} \label{ssec:categorif}
An (additive) \emph{categorification} of a cluster algebra $\cA(\widetilde{B})$
as in section~\ref{clustintro} is given by the following:
\begin{itemize}
\item[(A)] 
A $\C$-linear, Hom-finite Frobenius category $\cE$ with a 
cluster structure in the sense of~\cite[II.1]{BIRS} on the basic 
$\cE$-maximal rigid objects.
\item[(B)]
A basic $\cE$-maximal rigid object $T$ such that 
$B(\Gamma_T)^\circ=\widetilde{B}$.
\item[(C)]
A cluster character $\chi_?\df\operatorname{obj}(\cE)\ra \C(y_1,\ldots,y_r)$
in the sense of Palu~\cite[Definition 1.2]{P}, with triangles
replaced by short exact sequences.
\item[(D)]
The cluster character $\chi_?$ induces a bijection between basic,
$T$-reachable $\cE$-maximal rigid objects and clusters in 
$\cA(\widetilde{B})$.
\end{itemize}

\begin{Rem} (1) Conditions (A)-(C) imply obviously that each cluster monomial
in $\cA(\widetilde{B})$ is of the form $\chi_R$ for some $\cE$-rigid object $R$.
Thus condition $(D)$ is a kind of injectivity requirement for $\chi_?$.

(2) By the results in Section~\ref{intro1.6}  we have a cluster structure
on $\CC_w$. We can take $T=V_\bi$, for which we know $\Gamma_{V_\bi}$ by
Theorem~\ref{main6}. By Theorems~\ref{2.6i} and~\ref{2.6ii} our $\delta_?$
is a good candidate for a cluster character. In fact, by 
Theorems~\ref{main7} and~\ref{basesthm} below, we know that
$\delta_X\in\cA(B(\Gamma_{V_\bi})^\circ)$ for \emph{all} $X\in\CC_w$.
(By Theorem~\ref{main7} the algebra $\cA(B(\Gamma_{V_\bi})^\circ)$ is
up to isomorphism a subalgebra of $\MM^*$.) 
Property~(D) holds in our situation because of the construction
of the dual semicanonical basis. For this reason we call
$(\CC_w,V_\bi)$ a categorification of $\cA(B(\Gamma_{V_\bi})^\circ)$.
\end{Rem}


\section{Main results}


In this section, let $K = \C$ be the field of complex numbers.

\subsection{The cluster algebra $\cA(\CC_w)$ as a subalgebra of 
$\MM^*\equiv U(\n)^*_{\rm gr}$}\label{reach}
For a reduced expression $\ii = (i_r,\ldots,i_1)$ of a Weyl group
element $w$
let $\T(\CC_w)$ 
be the graph with vertices the isomorphism
classes of basic $\CC_w$-maximal rigid $\LL$-modules and with edges given 
by mutations. 
Let $T = T_1\oplus \cdots \oplus T_r$ be a vertex of $\T(\CC_w)$, and
let $\T(\CC_w,T)$ 
denote the connected 
component of $\T(\CC_w)$ containing $T$.
Two modules in $\T(\CC_w)$ are {\it mutation equivalent}
if they belong to the same connected component.
A $\LL$-module $X$ is called $T$-{\it reachable} 
if $X \in \add(R)$ for some vertex $R$ of $\T(\CC_w,T)$. 
Denote by $\RR(\CC_w,T)$ the subalgebra of $\MM^*$
generated by the $\delta_{R_i}\ (1\le i\le r)$ where 
$R = R_1 \oplus \cdots \oplus R_r$ runs
over all vertices of $\T(\CC_w,T)$. 
The following theorem is our first main result.
The proof is given in Section~\ref{clustercattoclusteralg}.

\begin{Thm}\label{main7}
Let $w$ be a Weyl group element.
Then the following hold:
\begin{itemize}

\item[(i)]
There is a unique isomorphism $\iota \df \cA(\CC_w,T)\to
\RR(\CC_w,T)$ such that 
\[
\iota(y_i)=\delta_{T_i} \qquad (1\le i\le r);
\]

\item[(ii)]
If we identify the two algebras $\cA(\CC_w,T)$ and
$\RR(\CC_w,T)$ via $\iota$, then the clusters of
$\cA(\CC_w,T)$ are identified with the $r$-tuples
$\delta(R) = (\delta_{R_1},\ldots,\delta_{R_r})$, where $R$
runs over the vertices of the graph $\T(\CC_w,T)$.
In particular, $\{ \delta_X \mid X \text{ is $T$-reachable} \}$ is
the set of cluster monomials in $\RR(\CC_w,T)$,
and all cluster monomials belong to 
the dual semicanonical basis 
$\cS^*$ of $\MM^*\equiv U(\n)^*_{\rm gr}$.

\end{itemize}
\end{Thm}

The proof of Theorem~\ref{main7} 
relies on Theorem~\ref{main4} and the multiplication formula
in Theorem~\ref{2.6ii}.

Write $\RR(\CC_w) := \RR(\CC_w,V_\ii)$.
(The algebra $\RR(\CC_w)$ and its cluster algebra structure
do not depend on $\ii$, since all $\CC_w$-maximal rigid
modules of the form $V_\ii$ are mutation equivalent,
see \cite[Proposition III.4.3]{BIRS}.)
Theorem~\ref{main7} shows that the cluster algebra 
$\cA(\CC_w)$ is canonically isomorphic to the subalgebra
$\RR(\CC_w)$ of $U(\n)^*_{\rm gr}$.

As an application, our theory provides an algorithm which computes
the Euler characteristics $\chi_{\rm c}(\F_{\kk,R})$
for all cluster monomials $\delta_R$ in $\RR(\CC_w)$
and all composition series types $\kk = (k_1,\ldots,k_s)$,
see Section~\ref{Euler2}. This is quite remarkable, since starting
from the definitions this seems to be an impossible task in almost
all cases.

\subsection{Dual PBW-bases and dual semicanonical bases}
Let $\ii = (i_r,\ldots,i_1)$ be a reduced expression of a Weyl
group element $w$.
Let $V_\ii = V_1 \oplus \cdots \oplus V_r$ be defined
as before.
For each $1 \le k \le r$ 
there is a canonical embedding 
$$
\iota_k\df V_{k^-} \to V_k.
$$
Here we set $V_0 := 0$.
Let $M_k$ be the cokernel of $\iota_k$, and define
$$
M_\ii := M_1 \oplus \cdots \oplus M_r.
$$
These modules play an important role in our theory.
(In case $w$ is adaptable and $\ii$ is $Q^\op$-adapted, 
the module $M_\ii$ is a terminal
$KQ$-module in the sense of \cite{GLSUni1}.)

In the spirit of Ringel's 
construction of PBW-bases for quantum groups \cite{Ri6},
we construct dual PBW-bases for our cluster algebras 
$\cA(\CC_w)$.
The following theorem is our second main result.
The proof will be given in
Section~\ref{PBWsection}.

\begin{Thm}\label{basesthm}
Let $\ii = (i_r,\ldots,i_1)$ be a reduced expression of a
Weyl group element $w$, and let $M_\ii = M_1 \oplus \cdots \oplus M_r$ 
be defined as above.
\begin{itemize}

\item[(i)]
The cluster algebra $\RR(\CC_w)$ is a polynomial ring in
$r$ variables.
More precisely, we have
$$
\RR(\CC_w) = \C[\delta_{M_1},\ldots,\delta_{M_r}] 
= \Span_\C\ebrace{\delta_X \mid X \in \CC_w};
$$

\item[(ii)]
The set 
$
\left\{ \delta_M \mid M \in \add(M_\ii) \right\}
$ 
is a $\C$-basis of $\RR(\CC_w)$;

\item[(iii)]
The subset $\cS_w^* := \cS^* \cap \RR(\CC_w)$
of the dual semicanonical basis
is a $\C$-basis of $\RR(\CC_w)$ containing all cluster
monomials.

\end{itemize}
\end{Thm}

Let 
$\widetilde{\RR}(\CC_w)$
be the algebra
obtained from $\RR(\CC_w)$ by formally inverting
the elements $\delta_P$ for all $\CC_w$-projective-injectives $P$.
In other words, $\widetilde{\RR}(\CC_w)$ is the cluster algebra
obtained from $\RR(\CC_w)$ by inverting the generators of
its coefficient ring. 
Similarly, let $\underline{\RR}(\CC_w)$ be
the cluster algebra obtained from $\RR(\CC_w)$ by specializing
the elements $\delta_P$ to 1.
For both cluster algebras  $\widetilde{\RR}(\CC_w)$
and $\underline{\RR}(\CC_w)$ we get a $\C$-basis which is
easily obtained from the dual semicanonical basis $\cS^*_w$
and again contains all cluster monomials, see Sections~\ref{generalbase} and
\ref{sectinvertspecial}.

\subsection{The shift functor in $\stCC_w$}    
As mentioned before, the category $\stCC_w$ is a triangulated
category with shift functor $\Omega_w^{-1}$.
Recall that $V_\ii = V_1 \oplus \cdots \oplus V_r$ is
a basic $\CC_w$-maximal rigid module.
Set $T_\ii := I_w \oplus \Omega_w^{-1}(V_\ii)$.
In Section~\ref{algorithm} we construct a sequence of mutations which
starts in $V_\ii$ and ends in $T_\ii$.
This mutation sequence is crucial for the proof of some of our results.
(For example, it helps to show that the coordinate rings $\C[N(w)]$ and
$\C[N^w]$ are generated by the set of cluster variables.) 

Now let $R = R_1 \oplus \cdots \oplus R_r$ be
any $\CC_w$-maximal rigid $\LL$-module, which is mutation equivalent
to $V_\ii$.
Suppose that we know a sequence of mutations starting in $V_\ii$
and ending in $R$.
Then we can use the mutation sequence from $V_\ii$ to $T_\ii$
to obtain a mutation sequence between $R$ and
$I_w \oplus \Omega_w^{-1}(R)$, and between $R$ and
$I_w \oplus \Omega_w(R)$, see Section~\ref{shiftcomb}.

\subsection{Unipotent subgroups and cells}\label{introbla}
Let $w$ be a Weyl group element and put
$\Delta^+_w:=\{\alpha\in\Delta^+\mid w(\alpha)<0\}$.
Let
\[
\n(w) = \bigoplus_{\alpha\in\Delta^+_w} \n_\alpha
\]
be the corresponding sum of root subspaces of $\n$, see Section~\ref{nw}.
This is a finite-dimensional nilpotent Lie algebra.
Let $N(w)$ be the corresponding finite-dimensional unipotent group,
see Section~\ref{ssec:nilp}.

The maximal Kac-Moody group attached to $\g$ as in~\cite[Chapter 6]{Ku}
comes with a pair of subgroups $N$ and $N_-$ 
(denoted by $\U$ and $\U_-$ in~\cite{Ku}). Note that later on for the
definition of generalized minors in Section~\ref{minorssection} we also
have to work with the minimal Kac-Moody group (denoted by 
$\mathcal{G}^{\operatorname{min}}$ in in~\cite[7.4]{Ku}).
We have
\[
N(w) = N \cap (w^{-1}N_-w).
\]
We also define the {\it unipotent cell}
\[
N^w :=  N\cap(B_-wB_-)
\]
where $B_-$ is the standard negative Borel subgroup of the Kac-Moody group. 

Every $\LL$-module $X$ in $\CC_w$ gives rise to a linear form
$\delta_X \in \MM^* \equiv U(\n)_{\rm gr}^*$ and by means of the 
identification $U(\n)_{\rm gr}^* \equiv \C[N]$ to a regular function $\vph_X$
on $N$.

The following theorem, proved in Section~\ref{coordinaterings}, 
is our third main result.

\begin{Thm}\label{wident}
The algebras $\C[N(w)]$ and $\C[N^w]$ 
of regular functions on $N(w)$ and $N^w$, respectively,
have a cluster algebra structure.
For each reduced expression $\ii = (i_r,\ldots,i_1)$
of $w$, the tuple 
$((\varphi_{V_{\ii,1}},\ldots,\varphi_{V_{\ii,r}}), B(\Gamma_{V_\bi})^\circ)$
provides an initial seed of these cluster algebra structures.
The functions $\varphi_{V_{\ii,k}} \in \C[N]$ can be interpreted
as generalized minors.
We obtain natural cluster algebra isomorphisms
$$
\C[N(w)] \cong \RR(\CC_w) 
\text{ and }
\C[N^w] \cong \widetilde{\RR}(\CC_w).
$$
\end{Thm}

As a result, we have obtained a categorification in the sense 
of~\ref{ssec:categorif} of the cluster algebra structure on 
$\C[N(w)]$ and $\C[N^w]$.

\subsection{Example}\label{exampleintro}
We are going to illustrate some of the previous results on an example.
Let $Q$ be a quiver with underlying graph
$\xymatrix@-0.5pc{1 \ar@{-}[r] & 2 \ar@{-}[r] & 3 \ar@{-}[r] &4}$
and let
$\ii := (3,4,2,1,3,4,2,1)$.
This is a reduced expression of the Weyl group element
$w := s_3s_4s_2s_1s_3s_4s_2s_1$.
The category $\CC_w$ contains 18 indecomposable modules,
and 4 of these are $\CC_w$-projective-injective.
The stable category ${\stCC}_w$ is triangle equivalent to the
cluster category $\CC_Q$.

The maximal rigid module $V_\ii$ has 8 indecomposable direct summands,
namely
\begin{align*}
V_1 &= \bsm1\esm &
V_2 &= \bsm1\\&2\esm &
V_3 &= \bsm4\esm &
V_4 &= \bsm1\\&2&&4\\&&3\esm\\
V_5 &= I_{\ii,1} = \bsm&2\\1\esm &
V_6 &= I_{\ii,2} = \bsm&2&&4\\1&&3\\&2\esm &
V_7 &= I_{\ii,4} = \bsm1\\&2\\&&3\\&&&4\esm &
V_8 &= I_{\ii,3} = \bsm&2\\1&&3\\&2&&4\\&&3\esm
\end{align*}
Similarly, $T_\ii$ has 4 non-$\CC_w$-projective-injective 
indecomposable direct summands,
namely 
\begin{align*}
T_1 &= \bsm2\esm &
T_2 &= \bsm2&&4\\&3\esm &
T_3 &= \bsm1\\&2\\&&3\esm &
T_4 &= \bsm2\\&3\esm.
\end{align*}
Here we set $T_k := \Omega_w^{-1}(V_k)$ for $1 \le k \le 4$.

The group $N$ can be taken to be the group of upper
unitriangular $5\times 5$ matrices with complex coefficients.
Given two subsets $I$ and $J$ of $\{1,2,\ldots,5\}$
with $|I|=|J|$, we denote 
by $D_{IJ}\in\C[N]$ the regular function mapping an element
$x\in N$ to its minor $D_{IJ}(x)$ with row subset $I$ and 
column subset $J$.  
We get
\begin{align*}
\vph_{V_1} &= D_{\{1\},\{2\}} &
\vph_{V_2} &= D_{\{12\},\{23\}} &
\vph_{V_3} &= D_{\{1234\},\{1235\}} &
\vph_{V_4} &= D_{\{123\},\{235\}} \\
\vph_{V_5} &= D_{\{1\},\{3\}} &
\vph_{V_6} &= D_{\{12\},\{35\}} &
\vph_{V_7} &= D_{\{1234\},\{2345\}} &
\vph_{V_8} &= D_{\{123\},\{345\}} \\
\vph_{T_1} &= D_{\{12\},\{13\}} &
\vph_{T_2} &= D_{\{123\},\{135\}} &
\vph_{T_3} &= D_{\{123\},\{234\}} &
\vph_{T_4} &= D_{\{123\},\{134\}}.
\end{align*}
The unipotent subgroup $N(w)$ consists of all
$5\times 5$ matrices of the form
$$
\left[
\begin{matrix}
1&u_1&u_2&u_7&u_4\\
0&1&u_5&u_8&u_6\\
0&0&1&0&0\\
0&0&0&1&u_3\\
0&0&0&0&1
\end{matrix}
\right],
\qquad
(u_1,\ldots,u_8\in\C).
$$
The unipotent cell $N^w$ is a locally closed subset of $N$
defined by the following equations and inequalities:
\begin{multline*}
N^w = \{ x \in N \mid
D_{\{1\},\{4\}}(x) = D_{\{1\},\{5\}}(x) = D_{\{12\},\{45\}}(x) = 0,\
D_{\{1\},\{3\}}(x) \not= 0,\\
D_{\{12\},\{35\}}(x) \not= 0,\
D_{\{123\},\{345\}}(x) \not= 0,\
D_{\{1234\},\{2345\}}(x) \not= 0
\}.
\end{multline*}
Note that the 4 inequalities are given by the non-vanishing
of the 4 regular functions 
$\vph_{I_{\ii,j}}$, $1 \le j \le 4$
attached to the indecomposable $\CC_w$-projective-injective modules.
We have
\begin{align*}
D_{\{1\},\{4\}} &= \vph_{\bsm&&3\\&2\\1\esm} &
D_{\{1\},\{5\}} &= \vph_{\bsm&&&4\\&&3\\&2\\1\esm} &
D_{\{12\},\{45\}} &= \vph_{\bsm&&3\\&2&&4\\1&&3\\&2\esm}.
\end{align*}
Our results show that the polynomial algebra $\C[N(w)]$
has a cluster algebra structure, of which
$(\vph_{V_1},\vph_{V_2},\vph_{V_3},\vph_{V_4},
\vph_{I_{\ii,1}},\vph_{I_{\ii,2}},\vph_{I_{\ii,3}},\vph_{I_{\ii,4}})$
is a distinguished cluster.
Its coefficient ring is the polynomial ring in the four variables
$(\vph_{I_{\ii,1}},\vph_{I_{\ii,2}},\vph_{I_{\ii,3}},\vph_{I_{\ii,4}})$.
The cluster mutations of this algebra come from mutations of the
basic $\CC_w$-maximal rigid $\LL$-modules.
Moreover, if we replace the coefficient ring by the ring of
Laurent polynomials in the four variables
$(\vph_{I_{\ii,1}},\vph_{I_{\ii,2}},\vph_{I_{\ii,3}},\vph_{I_{\ii,4}})$,
we obtain the coordinate ring $\C[N^w]$.


\section{Kac-Moody Lie algebras}
\label{semican}


From now on, let $K = \C$ be the field of complex numbers.
In this section we recall known results on Kac-Moody
Lie algebras.

\subsection{Kac-Moody Lie algebras}\label{KMalg}
Let $\GG=(\GG_0,\GG_1,\gamma)$ be a finite graph (without loops).  
It has as set of vertices $\GG_0$, edges $\GG_1$ and 
$\gamma\df\GG_1\to \cP_2(\GG_0)$ determining
the adjacency of the edges; 
here $\cP_2(\GG_0)$ denotes the set
of two-element subsets of $\GG_0$. 
If $\GG_0 = \{ 1,2,\ldots,n \}$ we can assign
to $\GG$ a {\em symmetric generalized Cartan matrix} 
$C_\GG = (c_{ij})_{1 \le i,j \le n}$, which is an $n \times n$-matrix
with integer entries
$$
c_{ij} :=
\begin{cases}
2 & \text{if $i = j$},\\
-\abs{\gamma^{-1}(\{i,j\})} & \text{if $i \neq j$}.
\end{cases}
$$ 
Obviously, the assignment
$\GG \mapsto C_\GG$ induces a bijection between isomorphism classes of
graphs with vertex set $\{ 1,2,\ldots,n \}$ and symmetric generalized 
Cartan matrices
in $\Z^{n\times n}$ up to simultaneous permutation of rows and columns.

For a quiver $Q=(Q_0,Q_1,s,t)$ as defined in Section~\ref{prenil},
its underlying graph $\abs{Q}:=(Q_0,Q_1,q)$ is given
by $q(a)=\{s(a),t(a)\}$ for all $a \in Q_1$ \textit{i.e.}~it is obtained by 
``forgetting'' the orientation of the edges.
We write $C_Q := C_{\abs{Q}} := (c_{ij})_{i,j}$.

Let $\g := \g_Q := \g(C_Q)$ be the 
(symmetric) {\it Kac-Moody Lie algebra} (see~\cite{Ka})
associated to $Q$, which is defined as follows:
Let $\h$ be a $\C$-vector space of dimension $2n-\operatorname{rank}(C_Q)$,
and let $\Pi := \{ \alpha_1,\ldots,\alpha_n \} \subset \h^*$ and 
$\Pi^\vee := \{ \alpha_1^\vee,\ldots,\alpha_n^\vee \} \subset \h$
be linearly independent subsets of the vector spaces $\h^*$ and $\h$,
respectively, such that
$$
\alpha_i(\alpha_j^\vee) = c_{ij}
$$
for all $i,j$.

Let $\h^* = \h_1^* \oplus \h_2^*$ be a vector space
decomposition, where $\h_1^*$ is just the subspace with 
basis $\Pi$, and $\h_2^*$ is any
direct complement of $\h_1^*$ in $\h^*$.
Let $(-,-)\df \h^* \times \h^* \to \C$
be the standard bilinear form, defined by
$(\alpha_j,\alpha_j) := \alpha_i(\alpha_j^\vee)$,
$(\alpha_i,x) := (x,\alpha_i) := x(\alpha_i^\vee)$,
and $(x,y) := 0$ for all $x,y \in \h_2^*$ and $1 \le i,j \le n$.
Note that $\alpha_i(\alpha_j^\vee) = (\dimv(S_i),\dimv(S_j))_Q$,
where $(-,-)_Q$ is the bilinear form defined in Section~\ref{prenil}.

Now $\g = (\g,[-,-])$ is the Lie algebra over $\C$
generated by $\h$ and the symbols $e_i$
and $f_i\ (1 \le i \le n)$ satisfying the following 
defining relations:
\begin{itemize}

\item[(L1)]
$[h,h'] = 0$ for all $h,h' \in \h$,

\item[(L2)]
$[h,e_i] = \alpha_i(h)e_i$, and
$[h,f_i] = -\alpha_i(h)f_i$,

\item[(L3)]
$[e_i,f_i] = \alpha_i^\vee$ and
$[e_i,f_j] = 0$ for all $i \not= j$,

\item[(L4)]
$(\ad(e_i)^{1-c_{ij}})(e_j) = 0$ for all $i \not= j$,

\item[(L5)]
$(\ad (f_i)^{1-c_{ij}})(f_j) = 0$ for all $i \not= j$.
\end{itemize}
(For $x,y \in \g$ and $m \ge 1$
we set $\ad(x)(y) := \ad(x)^1(y) := [x,y]$
and
$\ad(x)^{m+1}(y) := \ad(x)^m([x,y])$.)

The Lie algebra $\g$ is finite-dimensional if and only if
$Q$ is a Dynkin quiver. In this case, this is the usual Serre presentation 
of the simple Lie algebra associated to $Q$.

Conversely, if $\g = \g(C)$ is a Kac-Moody Lie algebra defined by a symmetric
generalized Cartan matrix $C$, we say that $\g$ is of type $\GG$ if 
$C = C_\GG$.
This is well defined for symmetric Kac-Moody Lie algebras.
We call $\GG$ the {\it Dynkin graph} of $\g$.

For $\alpha \in \h^*$ let
$$
\g_\alpha := 
\{ x \in \g \mid [h,x] = \alpha(h)x \text{ for all $h \in \h$} \}.
$$
One can show that $\dm \g_\alpha < \infty$ for all $\alpha$.
By 
$$
R := \sum_{i=1}^n \Z\alpha_i
$$
we denote the {\it root lattice} 
of $\g$.
Define 
$R^+ :=  \N\alpha_1 \oplus \cdots \oplus \N\alpha_n$.
The {\it roots} of $\g$ 
are defined as the elements in
$$
\Delta := \{  \alpha \in R \setminus \{0\} \mid \g_\alpha \not= 0 \}.
$$
Set $\Delta^+ := \Delta \cap R^+$ and $\Delta^- := \Delta \cap (-R^+)$.
One can show that $\Delta = \Delta^+ \cup \Delta^-$.
The elements in $\Delta^+$ and $\Delta^-$ are the {\it positive roots}
and the {\it negative roots}, 
respectively.
The elements in $\{ \alpha_1,\ldots,\alpha_n \}$ are positive roots
of $\g$ and are called {\it simple roots}.

One has the triangular decomposition
$
\g = \n_- \oplus \h \oplus \n
$
with 
$$
\n_- = \bigoplus_{\alpha \in \Delta^+} \g_{-\alpha}
\text{\;\;\; and \;\;\;}
\n = \bigoplus_{\alpha \in \Delta^+} \g_\alpha.
$$
The Lie algebra $\n$ is generated by $e_1,\ldots,e_n$
with defining relations (L4).
Set $\n_\alpha := \g_\alpha$ if $\alpha \in R^+ \setminus \{0\}$.

For $1 \le i \le n$ define an element $s_i$
in the automorphism group ${\rm Aut}(\h^*)$ of $\h^*$ by
$$
s_i(\alpha) := \alpha - \alpha(\alpha_i^\vee)\alpha_i
$$
for all $\alpha \in \h^*$.
The subgroup $W \subset {\rm Aut}(\h^*)$ generated by $s_1,\ldots,s_n$
is the {\it Weyl group} 
of $\g$.
The elements $s_i$ are called
{\it Coxeter generators} 
of $W$.
The identity element of $W$ is denoted by 1.
The {\it length} $l(w)$ 
of some $w \not= 1$ in $W$ is 
the smallest number $t \ge 1$ such that $w = s_{i_t}\cdots s_{i_2} s_{i_1}$
for some $1 \le i_j \le n$.
In this case $(i_t,\ldots,i_2, i_1)$ is a
{\it reduced expression} for $w$.
Let $R(w)$ be the set of all reduced expressions for $w$.
We set $l(1) = 0$.

A root $\alpha \in \Delta$ is a {\it real root} 
if
$\alpha = w(\alpha_i)$
for some $w \in W$ and some $i$.
It is well known that $\dm \g_\alpha = 1$ if $\alpha$ is a real root. 
By $\Delta_{\rm re}$ 
we denote the set of real roots of $\g$.
Define $\Delta_{\rm re}^+ := \Delta_{\rm re} \cap \Delta^+$.

Finally, let us fix a basis 
$\{\vpi_j\mid 1\le j \le 2n - \operatorname{rank}(C_Q)\}$ of
$\mathfrak{h}^*$ such that
\[
\vpi_j(\alpha^\vee_i) = \delta_{ij},\qquad
(1\le i\le n,\ 1\le j\le 2n-\operatorname{rank}(C_Q)).  
\]
The $\vpi_j$ are the {\it fundamental weights}. 
We denote by 
$$
P := \{ \nu \in \h^* \mid \nu(\alpha_i^\vee) \in \Z \text{ for all }
1 \le i \le n\}
$$
the {\it integral weight lattice}, and we set
$$
P^+ := \{ \nu\in P\mid \nu(\alpha_i^\vee)\ge 0\text{ for all } 
1 \le i \leq n \}.
$$
The elements in $P^+$ are called {\it integral dominant weights}.
We have 
$$
P = \bigoplus_{j=1}^n \Z \vpi_j \oplus
\bigoplus_{j=n+1}^{2n-\operatorname{rank}(C_Q)} \C\vpi_j
\text{\;\;\; and \;\;\;}
P^+ = \bigoplus_{j=1}^n \N \vpi_j \oplus
\bigoplus_{j=n+1}^{2n-\operatorname{rank}(C_Q)} \C\vpi_j
$$
Define
$$
\overline{P} := \bigoplus_{j=1}^n \Z\vpi_j
\text{\;\;\; and \;\;\;}
\overline{P}^+ := \bigoplus_{j=1}^n \N\vpi_j.
$$ 
The lattice $\overline{P}$ can be naturally identified with the weight
lattice of the derived subalgebra 
$\g' := [\g,\g]$ of $\g$.

\subsection{The universal enveloping algebra $U(\n)$}\label{universal}
The universal enveloping algebra $U(\n)$ 
of the Lie algebra
$\n$ is the associative $\C$-algebra defined by generators
$E_1,\ldots,E_n$ and relations
$$
\sum_{k=0}^{1-c_{ij}} (-1)^k E_i^{(k)}E_jE_i^{(1-c_{ij}-k)} = 0
$$
for all $i \not= j$,
where the $c_{ij}$ are the entries of the generalized Cartan matrix
$C_Q$, and let
$$
E_i^{(k)} := E_i^k/k!.
$$
We have a canonical embedding 
$\iota\df \n \to U(\n)$ which maps $e_i$ to $E_i$ for all $1 \le i \le n$.
We consider $\n$ as a subspace of $U(\n)$, and we also identify
$e_i$ and $E_i$.

Let 
$$
J = 
\begin{cases}
\N_1& \text{if $\dim(\n) = \infty$},\\
\{ 1,2,\ldots,d \} & \text{if $\dim(\n) = d$}.
\end{cases}
$$
Let ${\rm P} := \{ p_i \mid i \in J \}$ be a $\C$-basis of $\n$
such that ${\rm P} \cap \n_\alpha$ is a basis of $\n_\alpha$
for all positive roots $\alpha$.
We assume that $\{ e_1,\ldots,e_n \} \subset {\rm P}$.
Thus $e_i$ is a basis vector of the
(1-dimensional) space $\n_{\alpha_i}$.
For $k \ge 0$ define 
$$
p_i^{(k)} := p_i^k/k!.
$$
Let $\N^{(J)}_{\phantom{0}}$ be the set of
tuples $(m_i)_{i \in J}$ of natural numbers $m_i$ such that
$m_i = 0$ for all but finitely many $m_i$.
For ${\bf m} = (m_i)_{i \ge 1} \in \N^{(J)}_{\phantom{0}}$ define
$$
p_{\bf m} := p_1^{(m_1)} p_2^{(m_2)}\cdots p_s^{(m_s)}
$$
where $s$ is chosen such that $m_j = 0$ for all $j > s$.

\begin{Thm}[Poincar{\'e}-Birkhoff-Witt]\label{PBW1}
The set
$$
\cP := \left\{ p_{\bf m} \mid {\bf m} \in \N^{(J)}_{\phantom{0}} \right\}
$$
is a $\C$-basis of $U(\n)$.
\end{Thm}

The basis $\cP$ is called a {\it PBW-basis} of $U(\n)$.
For $d = (d_1,\ldots,d_n) \in \N^n$ let
$U_d$ be the subspace of $U(\n)$ spanned by the 
elements of the form $e_{i_1}e_{i_2} \cdots e_{i_m}$,
where for each $1 \le i \le n$ the set
$\{ k \mid i_k = i, 1\le k \le m \}$ contains exactly
$d_i$ elements.
It follows that
$$
U(\n) = \bigoplus_{d \in \N^n} U_d.
$$
This turns $U(\n)$ into an $\N^n$-graded algebra.

Furthermore, $U(\n)$ is a cocommutative Hopf algebra
with comultiplication
$$
\Delta\df U(\n) \to U(\n) \otimes U(\n)
$$
defined by $\Delta(x) := 1 \otimes x + x \otimes 1$
for all $x \in \n$.
It is easy to check that 
\begin{equation}\label{eqPBWdelta}
\Delta(p_{\bf m}) = \sum_{\bf k} p_{\bf k} \otimes p_{{\bf m}-{\bf k}},
\end{equation}
where the sum is over all tuples ${\bf k} = (k_i)_{i\ge 1}$ with
$0\le k_i \le m_i$ for every $i$.

By $U_d^*$ we denote the vector space dual of $U_d$.
Define the {\it graded dual} of $U(\n)$ 
by
$$
U(\n)^*_{\rm gr} := \bigoplus_{d \in \N^n} U_d^*.
$$
It follows that
$U(\n)^*_{\rm gr}$ is a commutative associative $\C$-algebra
with multiplication defined via the comultiplication
$\Delta$ of $U(\n)$:
For $f',f'' \in U(\n)^*_{\rm gr}$ and
$x \in U(\n)$, we have
$$
(f' \cdot f'')(x) = \sum_{(x)} f'(x_{(1)})f''(x_{(2)}),
$$
where (using the Sweedler notation) we write
$$
\Delta(x) = \sum_{(x)} x_{(1)} \otimes x_{(2)}.
$$

Let
$
\cP^* := 
\left\{ p_{\bf m}^* \mid {\bf m} \in \N^{(J)}_{\phantom{0}} \right\}
$
be the dual PBW-basis of $U(\n)^*_{\rm gr}$, where 
$$
p_{\bf m}^*(p_{\bf n}) :=
\begin{cases}
1 & \text{if ${\bf m} = {\bf n}$},\\
0 & \text{otherwise}.
\end{cases}
$$
The element in $\cP^*$ corresponding to $p_i \in {\rm P}$ is
denoted by $p_i^*$.
It follows from (\ref{eqPBWdelta}) that 
$$
p_{\bf m}^* \cdot p_{\bf n}^* = p_{{\bf m}+{\bf n}}^*,
$$
that is, each element $p_{\bf m}^*$ in $\cP^*$
is equal to a monomial in the $p_i^*$'s.
Hence, the graded dual $U(\n)^*_{\rm gr}$ can be identified with the
polynomial algebra $\C[p_1^*,p_2^*,\ldots]$ (with countably 
many variables $p_i^*$).

\subsection{The Lie algebra $\n(w)$}\label{nw}
Let 
$$
\widehat{\n} := \prod_{\alpha \in \Delta^+} \n_\alpha
$$
be the completion of $\n$.
A subset $\Theta \subseteq \Delta^+$ is {\it bracket closed}
if for all $\alpha,\beta \in \Theta$ with $\alpha+\beta \in \Delta^+$
we have $\alpha+\beta \in \Theta$.
In this case, we define
$$
\widehat{\n}(\Theta) := \prod_{\alpha \in \Theta} \n_\alpha.
$$
Since $\Theta$ is bracket closed, $\widehat{\n}(\Theta)$
is a Lie subalgebra of $\widehat{\n}$.
One calls $\Theta$ {\it bracket coclosed}
if $\Delta^+ \setminus \Theta$ is bracket closed.

For $w \in W$ set 
$\Delta_w^+ := \left\{ \alpha \in \Delta^+ \mid w(\alpha) < 0 \right\}$.
It is well known that for each reduced expression
$(i_r,\ldots,i_2, i_1) \in R(w)$ we have
$$
\Delta_w^+  = \{\alpha_{i_1}, s_{i_1}(\alpha_{i_2}), \ldots,
s_{i_1}s_{i_2}\cdots s_{i_{r-1}}(\alpha_{i_r}) \}.
$$
For $1 \le k \le r$ set
$$
\beta_\ii(k) := 
\begin{cases}
\alpha_{i_1} & \text{if $k = 1$},\\
s_{i_1}s_{i_2}\cdots s_{i_{k-1}}(\alpha_{i_k}) & \text{otherwise}.
\end{cases}
$$
The set $\Delta_w^+$ contains $l(w)$ positive roots, all
of these are real roots, see for example~\cite[1.3.14]{Ku}.
The next lemma is also well known.

\begin{Lem}\label{bracketclosed1}
For every $w \in W$, the set
$\Delta_w^+$ is bracket closed
and bracket coclosed.
\end{Lem}

Let $\n(w) := \widehat{\n}(\Delta_w^+)$
be the {\it nilpotent Lie algebra associated to} $w$.
We have
$$
\n(w) = \bigoplus_{\alpha \in \Delta_w^+} \n_\alpha 
$$
and $\dm \n(w) = l(w)$.

Again, let
$\ii = (i_r,\ldots,i_1)$ be a reduced expression.
As in Section~\ref{universal} we choose a $\C$-basis 
$
{\rm P} = \{ p_j \mid j \in J \}
$
such that ${\rm P} \cap \n_\alpha$ is a basis of $\n_\alpha$
for all positive roots $\alpha$.
The resulting PBW-basis
$\cP = \{ p_{\bf m} \mid {\bf m} \in \N^{(J)}_{\phantom{0}} \}$ 
of $U(\n)$
is called $\ii$-{\it compatible} provided
the vector $p_k$ belongs to $\n_{\beta_\ii(k)}$ for
all $1 \le k \le r$.
In this case
$$
\cP_\ii := \left\{ p_1^{(m_1)}p_2^{(m_2)} \cdots p_r^{(m_r)} \mid m_k \ge 0 
\text{ for all } 1 \le k \le r \right\}
$$
is a PBW-basis of the universal enveloping algebra $U(\n(w))$
of $\n(w)$, and
$$
\cP_\ii^* := \left\{ (p_1^*)^{m_1}(p_2^*)^{m_2} \cdots 
(p_r^*)^{m_r} \mid m_k \ge 0 
\text{ for all } 1 \le k \le r \right\}
$$
is the corresponding dual PBW-basis of the graded dual $U(\n(w))_{\rm gr}^*$.

\subsection{Highest weight modules}
A $U(\g)$-module $M$ is a {\it weight module} or
$\h$-{\it diagonalizable} if
$$
M = \bigoplus_{\mu\in \h^*} M_\mu
$$
where 
$$
M_\mu := \left\{ m \in M \mid h \cdot m 
= \mu(h)m \text{ for all }
h \in \h \right\}.
$$
For each vector $v \in M_\mu$ let ${\rm wt}(v) := \mu$ be its 
{\it weight}.
Analogously, one defines when a right $U(\g)$-module is a 
{\it weight module}.

A $U(\g)$-module $M$ is a
{\it highest weight module} if the following hold:
\begin{itemize}

\item
$M$ is a weight module;

\item
There is a vector $v \in M$ with $U(\g) \cdot v = M$;

\item
$e_i \cdot v = 0$ for all $i$.

\end{itemize}

A right $U(\g)$-module $M$ is a {\it lowest weight module} if 
the following hold:
\begin{itemize}

\item
$M$ is a weight module;

\item
There is a vector $u \in M$ with $u \cdot U(\g) = M$;

\item
$u \cdot f_i = 0$ for all $i$.

\end{itemize}
When we work with right $U(\g)$-modules, we invert the usual ordering on
weights.
So if $M$ is a lowest weight right $U(\g)$-module, then
the vector $u$ (which is uniquely determined up to a non-zero scalar)
has actually the lowest weight of $M$.
Indeed, if $m \in M_\mu$ and $h \in \h$, then we have
$$
(m \cdot e_i) \cdot h = (\mu(h)m) \cdot e_i -
(\alpha_i(h)m) \cdot e_i = 
(\mu - \alpha_i)(h)(m \cdot e_i).
$$
Here we used that $[h,e_i] = he_i - e_ih = \alpha_i(h)e_i$.
So $m \cdot e_i$ has weight $\mu-\alpha_i$.

\subsection{Construction of highest weight modules}\label{dualVerma}
In this section we present some of our results from 
\cite{GLSVerma} in a form
convenient for our present purpose.
For $\nu \in P^+$ we write
\[
\hI_\nu := \bigoplus_{i=1}^n \hI_i^{\,\nu(\alpha_i^\vee)}.
\]
For $1 \le i \le n$ and a nilpotent $\LL$-module $X$ we denote by
$\G(i,X)$ the variety of submodules $Y$ of $X$ such that 
$X/Y \cong S_i$.
Similarly, if  
$$
\soc(X) = \bigoplus_{i=1}^n S_i^{m_i}
$$ 
and 
$\nu \in P^+$ is such that $\nu(\alpha_i^{\vee})\geq m_i$ for $1 \le i \le n$, 
then we have an embedding $X \hookrightarrow \hI_\nu$.
In this case, we denote by $\G(i,\nu,X)$ the variety of submodules $Y$ of 
$\hI_\nu$ such that $X\subset Y$ and $Y/X \cong S_i$.  
Hence, if
$\dimv(X) = \beta$ and $f \in \MM_{\beta-\alpha_i}$, 
we can form the following sum
\[
\Sigma := \sum_{m \in \C} m \, \chi_{\rm c}(\{ Y \in \G(i,X)\mid f(Y)=m \}).
\]
For convenience we shall denote such an expression by an
integral, for example, 
\[
\Sigma = \int_{Y\in\G(i,X)} f(Y).
\]
Similarly, there exists a partition 
$$
\G(i,X) = \bigsqcup_{j=1}^m A_j
$$ 
into constructible subsets such that
$\delta_Y=\delta_{Y'}$ for all $Y,Y'\in A_j$.  
Then, choosing  arbitrary $Y_j\in A_j$ for $j=1, \ldots, m$,
we can also denote by an integral the following element of
$\MM_{\beta-\alpha_i}^*$
\[
\int_{Y\in\G(i,X)}\delta_Y =\sum_{j=1}^m \chi_{\rm c}(A_j)\delta_{Y_j}.
\]

\begin{Thm} 
Let $\lambda\in P$ be an integral weight,
and let $M_{\rm low}(\lambda)$ be the lowest weight Verma 
right $U(\g)$-module (with underlying vector space $U(\n)$)
with lowest weight $\lambda$.
Under the  identifications 
$$
M_{\rm low}(\lambda) \equiv U(\n) \equiv \MM
$$ 
the corresponding right $U(\g)$-module structure
on $\MM$ is described as follows:
The generators $e_i \in \n,\ f_i \in \n_-,\ h \in \h$ act
on
$g \in \MM_\beta$ by
\begin{align*}
(g \cdot e_i)(X')  &= \int_{Y\in\G(i,X')} g(Y),\\
(g \cdot f_i)(X) &= \int_{Y\in\G(i,\nu,X)}g(Y) - 
(\nu-\lambda)(\alpha_i^\vee)g(X\oplus S_i),\\
g\cdot h &= (\lambda -\beta)(h) g,
\end{align*}
where $X'\in\LL_{\beta+\alpha_i}$, $X\in\LL_{\beta-\alpha_i}$ and 
$\nu\in P^+$ are as above.
\end{Thm}

Note that $g\cdot e_i=g*\mathbf{1}_i$ by our convention for the multiplication
in $\MM$. Moreover, 
the formula for $g\cdot f_i\in\MM_{\beta-\alpha_i}$ is 
in fact independent of the choice of $\nu$.

For each $\h$-diagonalizable right $U(\g)$-module 
$$
M = \bigoplus_{\mu\in \h^*} M_\mu
$$ 
one can consider the {\em dual} representation
$$
M^* = \bigoplus_{\mu\in\h^*} M^*_\mu
$$ 
defined by
$M^*_\mu := \Hom_\C (M_{\mu},\C)$. It acquires the structure of 
a \emph{left} $U(\g)$-module via
\[
(x\cdot\phi)(m):=\phi( m\cdot x),
\qquad
(x \in U(\g),\ m\in M).
\]
Consider the canonical epimorphism from the Verma module
$M_{\rm low}(\lambda)$ to the irreducible lowest weight right 
$U(\g)$-module $L_{\rm low}(\lambda)$. For the
corresponding dual representations we obtain an inclusion
$$
L^*_{\rm low}(\lambda) \hookrightarrow M_{\rm low}^*(\lambda).
$$ 
It is well known that
$L_{\rm low}^*(\lambda)$ is isomorphic to the irreducible highest weight 
left $U(\g)$-module $L(\lambda)$ with highest weight $\lambda$.
This yields the following realization of the integrable module
$L(\lambda)$
in terms of $\delta$-functions.

\begin{Thm} \label{thm:irred-l}
Let $\lambda\in P^+$ be an integral dominant weight.
The subspace 
$$
U(\lambda) := \Span_\C\langle \delta_X \mid 
X \text{ submodule of } \hI_\lambda\rangle
$$
of $U(\n)^*_{\rm gr}$ 
carries the above-mentioned structure of an irreducible highest weight 
left $U(\g)$-module $L(\lambda)$. 
For such $X$ with $\dimv(X) = \beta$ 
the action of the Chevalley generators of $U(\g)$ is 
given by
\begin{align*}
e_i\cdot \delta_X &=\int_{Y\in\G(i,X)} \delta_Y,\\
f_i\cdot \delta_X &=\int_{Y'\in\G(i,\lambda,X)}\delta_{Y'},\\
h\cdot \delta_X &=(\lambda-\beta)(h)\delta_X.
\end{align*}
\end{Thm}

Note that $U(\n)_{\rm gr}^*$ carries also a {\em right} $U(\n)$-module
structure coming from the left regular representation of $U(\n)$. In
order to describe it, we introduce the following definition.
For $X \in \LL_\beta$ we denote by $\G'(i,X)$ the variety of submodules
$Y$ of $X$ such that $\dimv(Y) = \alpha_i$. Each element of this space
is isomorphic to $S_i$ and clearly $\G'(i,X)$ is a projective space.
It is easy to see that
\[
\delta_X \cdot e_i = \int_{S\in\G'(i,X)} \delta_{X/S}.
\]
Under the above identification 
$M^*_{\rm low}(\lambda) \equiv U(\n)_{\rm gr}^*$, the
subspace of $U(\n)^*_{\rm gr}$ carrying the
$U(\g)$-module $L(\lambda)$ can be described as follows.

\begin{Cor}\label{embedding1}
For $\lambda \in P^+$ we have
\[
U(\lambda) = \left\{\phi\in U(\n)^*_{\rm gr}\mid
 \phi\cdot e_i^{\lambda(\alpha_i^\vee)+1}=0 \text{ for all } 
1 \le i \le n \right\}.
\]
\end{Cor}

\begin{proof}
The nilpotent $\LL$-module $X$ is isomorphic to a submodule of $\hI_\lambda$
if and only if 
$$
\delta_X \cdot e_i^{\lambda(\alpha_i^\vee)+1}=0
$$ 
for every $i$.
The claim then follows from Theorem~\ref{thm:irred-l}.   
\end{proof}

Note that for $\lambda,\mu \in P^+$ we have 
$U(\lambda) = U(\mu)$ if and only if 
$$
\lambda - \mu \in \bigoplus_{j=n+1}^{2n - \operatorname{rank}(C_Q)}
\C\vpi_j.
$$


\section{Unipotent groups}


\subsection{The group $N$ and its coordinate ring $\C[N]$}\label{DefN}
The completion
$\widehat{\n}$ of $\n$ defined in~\ref{nw},
is a pro-nilpotent pro-Lie algebra, see~\cite[Section 6.1.1]{Ku}.
Let $N$ be the pro-unipotent pro-group
with Lie algebra $\widehat{\n}$.
We refer to Kumar's book \cite[Section 4.4]{Ku} for 
all missing definitions.

We can assume that $N = \widehat{\n}$ as a set and that the
multiplication of $N$ is defined via the Baker-Campbell-Hausdorff formula.
Hence the exponential map ${\rm Exp}\df \widehat{\n} \to N$
is just the identity map.

Put ${\mathcal H} := U(\n)^*_{\rm gr}$.
This is a commutative Hopf algebra.
We can regard ${\mathcal H}$ as the coordinate ring $\C[N]$ of $N$, that is,
we can identify $N$ with the set 
$$
{\rm maxSpec}({\mathcal H}) \equiv \Hom_{\rm alg}({\mathcal H},\C)
$$ 
of $\C$-algebra
homomorphisms ${\mathcal H} \to \C$.
An element $f \in \Hom_{\rm alg}({\mathcal H},\C)$ is determined
by the images $c_i := f(p_i^*)$ for all $i \ge 1$.

It is well known (see e.g. \cite[\S 3.4]{Ab}) that 
$\Hom_{\rm alg}({\mathcal H},\C)$ can also be identified with
the group $G({\mathcal H}^\circ)$ of all group-like elements of
the dual Hopf algebra ${\mathcal H}^\circ$ of ${\mathcal H}$, by mapping 
$f\in \Hom_{\rm alg}({\mathcal H},\C)$ to
\[
y_f = \sum_{\bf m}\left( \prod_i c_i^{m_i}\right) p_{\bf m}\in
G({\mathcal H}^\circ).
\]
Note that the map $f\mapsto y_f$ does not depend on the
choice of the PBW-basis $\cP = \{ p_{\bf m} \mid {\bf m} \in \N^{(J)} \}$.
Note also that $G({\mathcal H}^\circ)$ is contained in the vector space dual 
${\mathcal H}^*$
of ${\mathcal H}$, which
is the completion $\widehat{U(\n)}$ of $U(\n)$ with respect
to its natural grading.
When we use this second identification, an element 
$x\in N = \widehat{\n}$ is simply represented by
the group-like element 
$$
\exp(x) := \sum_{k\ge 0} x^k/k!
$$
in $\widehat{U(\n)}$.
To summarize, we have ${\mathcal H}=U(\n)^*_{\rm gr}\equiv \C[N]$ and
\[
N\equiv {\rm maxSpec}({\mathcal H}) \equiv \Hom_{\rm alg}({\mathcal H},\C)
\equiv 
G({\mathcal H}^\circ)\subset {\mathcal H}^\circ\subset {\mathcal H}^*
\equiv \widehat{U(\n)}.
\]

\subsection{The unipotent groups $N(w)$ and $N'(w)$}\label{ssec:nilp}
Let $\Theta$ be a bracket closed subset of $\Delta^+$, and
let 
$$
N(\Theta) := {\rm Exp}(\widehat{\n}(\Theta))
$$ 
be the corresponding
pro-unipotent pro-group.
For example, if $\alpha \in \Delta_{\rm re}^+$, then
$\Theta_\alpha := \{ \alpha \}$ is bracket closed.
In this case, $N(\alpha)$ is called the 
{\it one-parameter subgroup}
of $N$ associated to $\alpha$.
We have an isomorphism of groups $N(\alpha) \cong (\C,+)$.

If $\Theta$ is bracket closed and bracket coclosed, then
set $N'(\Theta) := N(\Delta^+ \setminus \Theta)$.
In this case, the multiplication in $N$ 
yields a bijection~\cite[Lemma 6.1.2]{Ku}
$$
m\df N(\Theta) \times N'(\Theta) \to N.
$$

For $w \in W$ let
$N(w) := N(\Delta_w^+)$.
This is a 
unipotent algebraic group
of dimension $l(w)$, and its Lie algebra is $\n(w)$.
Again we can identify $U(\n(w))_{\rm gr}^* \equiv \C[N(w)]$.
Similarly, define $N'(w) := N'(\Delta_w^+)$.


\section{Evaluation functions and generating functions of Euler 
characteristics}


Recall the identifications $\MM^* \equiv U(\n)^*_{\rm gr} = \C[N]$.
To every $X \in \nil(\LL)$, we have associated a linear form
$\delta_X \in U(\n)^*_{\rm gr}$. 
We shall also denote the evaluation function $\delta_X$ 
by $\vph_X$ when we 
regard it as a function on $N$.
For $1 \le i \le n$ define $x_i\df \C \to N$ by
$$
x_i(t) := \exp(te_i) = \sum_{k \ge 0} \frac{(te_i)^k}{k!}
$$ 
The following formula shows how to evaluate $\vph_X$
on a product of $x_i(t)$'s.

\begin{Prop}\label{phi_form}
Let $X \in \nil(\LL)$, and let $\ii = (i_1,\ldots,i_k)$ be any
sequence with $1 \le i_j \le n$ for all $1 \le j \le k$.
We have  
\[
\vph_X(x_{i_1}(t_1)\cdots x_{i_k}(t_k)) =
\sum_{\aaa=(a_1,\ldots,a_k)\in\N^k}\chi_{\rm c}(\F_{\ii^\aaa,X})
\frac{t_1^{a_1}\cdots t_k^{a_k}}{a_1!\cdots a_k!}.
\] 
Here $\ii^\aaa$ is short for the sequence 
$(i_1,\ldots,i_1,\ldots,i_k,\ldots,i_k)$ consisting of $a_1$ letters
$i_1$ followed by $a_2$ letters $i_2$, etc.
\end{Prop}

\begin{proof}
By Section~\ref{DefN} we can regard 
$x_{i_1}(t_1)\cdots x_{i_k}(t_k)$ as an element of $\widehat{U(\nn)}$,
namely,
\[
x_{i_1}(t_1)\cdots x_{i_k}(t_k)=\sum_{\aaa=(a_1,\ldots,a_k)\in\N^k}
\frac{t_1^{a_1}\cdots t_k^{a_k}}{a_1!\cdots a_k!}
e_{i_1}^{a_1}\cdots e_{i_k}^{a_k}.
\]
It follows from the identification of $\vph_X$ with $\delta_X$ that
\[
\vph_X(x_{i_1}(t_1)\cdots x_{i_k}(t_k))=
\sum_{\aaa=(a_1,\ldots,a_k)\in\N^k}
\frac{t_1^{a_1}\cdots t_k^{a_k}}{a_1!\cdots a_k!}
\delta_X(e_{i_1}^{a_1}\cdots e_{i_k}^{a_k}).
\] 
Now, in the geometric realization $\MM$ of the enveloping 
algebra $U(\n)$ in terms of constructible functions, 
$e_{i_1}^{a_1}\cdots e_{i_k}^{a_k}$ becomes the convolution 
product ${\bf 1}_{i_1}^{a_1} \star \cdots \star {\bf 1}_{i_k}^{a_k}$
and it is easy to see that 
$$
\delta_X({\bf 1}_{i_1}^{a_1} \star \cdots \star {\bf 1}_{i_k}^{a_k})
= \chi_{\rm c}(\F_{\ii^\aaa,X}).
$$
This finishes the proof.
\end{proof}

\begin{Rem}
The formula for $\vph_X$ given in \cite[\S 9]{GLSRigid} 
involves descending flags instead of ascending flags of
submodules of $X$. This is because in the present paper
we have taken a convolution product $\star$ opposite to
that of our previous papers, see Remark~\ref{newconv}.
\end{Rem}

Proposition~\ref{phi_form} says that we can think of
the $\varphi$-functions $\varphi_X$ as generating
functions of Euler characteristics.

For $\ii = (i_1,\ldots,i_k)$ and $\ba = (a_1,\ldots,a_k)$ as above
and $X \in \nil(\LL)$
let $\F_{\ii,\ba,X}$ be the projective variety of partial
composition series of type $(\ii,\ba)$ of $X$.
Thus an element of $\F_{\ii,\ba,X}$ is a 
chain
$$
0 = X_0 \subseteq X_1 \subseteq \cdots \subseteq X_k = X
$$
of submodules of $X$ such that $X_j/X_{j-1} \cong S_{i_j}^{a_j}$
for all $1 \le j \le k$.
There is an obvious surjective morphism 
$\pi_{\ii,\ba}\df \F_{\ii^\ba,X} \to \F_{\ii,\ba,X}$
whose fibers are all isomorphic to
$$
\F\left(\C^{a_1}\right) \times \cdots \times
\F\left(\C^{a_k}\right),
$$ 
where $\F\left(\C^m\right)$ is the variety of complete
flags of subspaces in $\C^m$.
In particular, we have
$$
\chi_{\rm c}\left(\F_{\ii^\ba,X}\right) = 
\chi_{\rm c}\left(\F_{\ii,\ba,X}\right) a_1! \cdots a_k!.
$$
Summarizing, we get
$$
\vph_X(x_{i_1}(t_1)\cdots x_{i_k}(t_k)) =
\sum_{\aaa=(a_1,\ldots,a_k)\in\N^k}\chi_{\rm c}(\F_{\ii,\aaa,X})
t_1^{a_1}\cdots t_k^{a_k}.
$$


\section{Generalized minors}\label{minorssection}


\subsection{Generalized minors}\label{genminors}
We start with some generalities on Kac-Moody groups. 
Let $\Gmin$ be the Kac-Moody group with $\Lie(\Gmin)=\g$ defined in
\cite[7.4]{Ku}.
It has a refined Tits system 
\[
(\Gmin, \Norm_{\Gmin}(H),N\cap\Gmin,N_-,H).
\] 
Write $\Nmin :=\Gmin\cap N$. 
Moreover, $\Gmin$ is
an affine ind-variety in a unique way
\cite[7.4.8]{Ku}.

For any real root $\alpha$ of $\g$, the one-parameter
subgroup $N(\alpha)$ is contained in $\Gmin$, 
and the $N(\alpha)$ together with $H$ generate $\Gmin$ as a group.
We have an anti-automorphism $g\mapsto g^T$ of $\Gmin$ which maps
$N(\alpha)$ to $N(-\alpha)$ for each real root $\alpha$, and fixes $H$.
We have another anti-automorphism $g\mapsto g^\iota$ which fixes 
$N(\alpha)$ for every real root $\alpha$, and $h^\iota = h^{-1}$
for every $h\in H$.

For each $\gamma \in \h^*$ there is a character $H \to \C^*$,
$a \mapsto a^\gamma$ defined by
$\exp(h)^\gamma := e^{\gamma(h)}$ for all $h \in \h$.

For $1 \le i \le n$ we have a unique homomorphism 
$\vph_i \df SL_2(\C) \to \Gmin$ satisfying
\[
\vph_i
\begin{pmatrix}
1& t \cr 0& 1
\end{pmatrix} 
=
\exp(te_i),
\qquad
\vph_i
\begin{pmatrix}
1& 0 \cr t& 1
\end{pmatrix} 
=
\exp(tf_i),
\qquad
(t\in\C).
\]
We define 
\[
\overline{s}_i :=
\vph_i
\begin{pmatrix}
0 & -1 \cr 1& 0
\end{pmatrix}. 
\]
For $w \in W$, we define 
$\overline{w} := \overline{s}_{i_r}\cdots\overline{s}_{i_1}$,
where $(i_r,\ldots,i_1)$ is a reduced expression for $w$.
Thus, we choose for every $w\in W$ a particular 
representative $\overline{w}$ of $w$ in the normalizer
$\Norm_\Gmin(H)$.

Let $L(\lambda)$ denote the irreducible highest weight $\g$-module
with highest weight $\lambda \in P^+$. 
Let $u_\lambda$ be a highest weight vector of $L(\lambda)$. 
This is an integrable module, so it is also a representation
of $\Gmin$.
For a reduced expression $\ii = (i_r,\ldots,i_1)$ of
a Weyl group element $w$, the vector
$$
\overline{s}_{i_1} \cdots \overline{s}_{i_r}(u_\lambda) \in L(\lambda)
$$
is an {\it extremal weight vector} of $L(\lambda)$, \textit{i.e.}~it belongs
to the extremal weight space $L(\lambda)_{w(\lambda)}$.
For a $U(\g)$-module $V$ and a
weight vector $v \in V_\mu$ define 
$$
f_i^{\rm max}v := f_i^{(m)}v
$$ 
where $m \ge 0$ is maximal such that $f_i^{(m)}v \not= 0$.
Similarly, define $e_i^{\rm max}$.
The following results can be found in \cite[Section 4.4.3]{J}:
We have
$$
\overline{s}_{i_1} \overline{s}_{i_2} \cdots \overline{s}_{i_r}(u_\lambda) = 
f_{i_1}^{\rm max}f_{i_2}^{\rm max} \cdots f_{i_r}^{\rm max}(u_\lambda)
$$
and
$$
e_{i_1}f_{i_2}^{\rm max} \cdots f_{i_r}^{\rm max}(u_\lambda) = 0.
$$
Furthermore, 
$$
{\rm wt}\left(\overline{s}_{i_1} \cdots 
\overline{s}_{i_r}(u_{\lambda})\right) 
= 
{\rm wt}\left(\overline{s}_{i_2} \cdots \overline{s}_{i_r}(u_{\lambda})\right) 
-  b_1\alpha_{i_1}
$$
where 
$b_1 := -(s_{i_1} \cdots s_{i_r}(\lambda),\alpha_{i_1})
= (s_{i_2} \cdots s_{i_r}(\lambda),\alpha_{i_1})$.

We have the following analogue of the Gaussian decomposition.

\begin{Prop}\label{Gauss}
Let $G_0$ be the subset $N_-\cdot H\cdot \Nmin$ of $\Gmin$.
\begin{itemize}

\item[(i)]
The subset $G_0$ is dense open in $\Gmin$ and each element $g\in G_0$ 
admits a unique
factorization $g=[g]_- [g]_0 [g]_+$ with $[g]_-\in N_-$, $[g]_0\in H$ and
$[g]_+\in \Nmin$.

\item[(ii)]
The map $g\mapsto [g]_+$ (resp.\@ $g\mapsto [g]_0$)
is a morphism of ind-varieties from $G_0$ to 
$\Nmin$ (resp.\@ to $H$). 

\end{itemize}
\end{Prop}

Part (i) follows from the fundamental properties of a refined Tits 
system~\cite[Theorem 5.2.3]{Ku}.
For part (ii), see \cite[Proposition 7.4.11]{Ku}.

Following Fomin and Zelevinsky~\cite{FZ5} we can now define for each 
$\vpi_j$ a generalized minor 
$\Delta_{\vpi_j,\vpi_j}$ as the regular function on $\Gmin$ such that 
\[
\Delta_{\vpi_j,\vpi_j}(g) = [g]_0^{\vpi_j},\qquad (g \in G_0). 
\]
For $w\in W$, we also define $\Delta_{\vpi_j,w(\vpi_j)}$ by
\[
\Delta_{\vpi_j,w(\vpi_j)}(g) := \Delta_{\vpi_j,\vpi_j}(g\overline{w}).
\]
The generalized minors $\Delta_{\vpi_j,\vpi_j}(g)$ 
have the following alternative description.

\begin{Prop}
Let $g \in \Gmin$.
The coefficient of $u_{\vpi_j}$ in the projection of $gu_{\vpi_j}$ 
on the weight
space $L(\vpi_j)_{\vpi_j}$ is equal to $\Delta_{\vpi_j,\vpi_j}(g)$.
\end{Prop}

\begin{proof}
Set $u_j := u_{\vpi_j}$.
Let $g = [g]_-[g]_0[g]_+ \in G_0$.
We have $[g]_+u_j = u_j$, and $[g]_0u_j=[g]_0^{\vpi_j} u_j$.
The result then follows from the fact that $[g]_-u_j$ is
equal to $u_j$ plus elements in lower weights.
\end{proof}

\begin{Prop}
We have 
$$
G_0 = 
\left\{ g \in \Gmin \mid \Delta_{\vpi_j,\vpi_j}(g) \neq 0 
\text{ for all } 1 \le j \le n \right\}.
$$
\end{Prop}

\begin{proof} 
Set $u_j := u_{\vpi_j}$.
We use the Birkhoff decomposition \cite[Theorem 5.2.3]{Ku}
\[
\Gmin = \bigsqcup_{w\in W} N_-\overline{w}H\Nmin,
\]
where $G_0$ is the subset of the right-hand side corresponding
to $w=e$.
If $g=[g]_-[g]_0[g]_+\in G_0$, then 
$\Delta_{\vpi_j,\vpi_j}(g)=[g]_0^{\vpi_j}\neq 0$.
Conversely,
if $g \not\in G_0$ we have $g= n_-whn$ for some $n_-\in N_-$,
$n\in \Nmin$, $h\in H$ and $w\not = e$.
Then for some $j$ we have $w(\vpi_j)\not = \vpi_j$ and
$\overline{w}hnu_j$ is a multiple of the extremal weight vector 
$\overline{w}u_j$.
Since the projection of $n_-\overline{w}u_j$ 
on the highest weight space  
of $L(\vpi_j)$ is zero, it follows that $\Delta_{\vpi_j,\vpi_j}(g)=0$.  
Finally, note that for any $j>n$ the minor $\Delta_{\vpi_j,\vpi_j}$ does 
not vanish on $\Gmin$. Indeed, the corresponding highest weight irreducible
module $L(\vpi_j)$ is one-dimensional since $\vpi_j(\alpha_i^\vee) = 0$
for any $i$.
Hence in the above description of $G_0$, we may omit the 
minors $\Delta_{\vpi_j,\vpi_j}$ with $j>n$.
\end{proof}

\subsection{The module $L(\lambda)$ as a subspace of $\C[N]$}
\label{adaptedordering}
For $w \in W$ and $1 \le j \le n$, we denote by 
$$
D_{\vpi_j,w(\vpi_j)}
$$ 
the restriction of the
generalized minor $\Delta_{\vpi_j,w(\vpi_j)}$ to $\Nmin$.
For example, $D_{\vpi_j,\vpi_j}$ is equal to the constant
function $1$.
In Section~\ref{catminors} we are going to show that each
(restricted) generalized minor 
$D_{\vpi_j,w(\vpi_j)}$ can be identified with a generating function
$\vph_X$ for a certain $\LL$-module $X$.
In order to do this, we need to recall some results on 
Kac-Moody groups.

Let $G' := [\Gmin,\Gmin]$ be the group constructed by 
Kac and Peterson \cite{KP1}, see \cite[Section 7.4.E (1)]{Ku}.
The associated Lie algebra is $\g' = [\g,\g]$.

Let $\C[G']_{\rm s.r.}$ denote the algebra of strongly regular
functions on $G'$ \cite[\S 2C]{KP1}.
Define the invariant ring
\[
\C[N_-\backslash G']_{\rm s.r.}
:=
\left\{f\in \C[G']_{\rm s.r.}\mid f(ng)=f(g)
\text{ for all } n\in N_-,\ g \in G' \right\}.
\]
This ring is endowed with the usual left action of $G'$
given by
\[
(g\cdot f)(g') := f(g'g),\qquad (f\in \C[N_-\backslash G']_{\rm s.r.},
\,g,g'\in G').
\] 
It was proved by Kac and Peterson \cite[Corollary 2.2]{KP1} that
as a left $G'$-module, it decomposes as follows
\[
\C[N_-\backslash G']_{\rm s.r.} = 
\bigoplus_{\lambda \in \overline{P}^+} L(\lambda).
\]
This is a multiplicity-free decomposition, in which the 
irreducible highest weight module $L(\lambda)$ is carried by the subspace
\[
S(\lambda) = \left\{ f \in \C[N_-\backslash G']_{\rm s.r.}
\mid f(hg) = \Delta_\lambda(h)f(g) 
\text{ for all } h\in H,\ g \in G' \right\},
\]
where we denote 
$$
\Delta_\lambda : = \prod_{j=1}^n 
\Delta_{\vpi_j,\vpi_j}^{\lambda(\alpha_j^\vee)}.
$$
Clearly, $\Delta_\lambda$ is contained in $S(\lambda)$, 
and it is a highest weight vector.
Moreover, for any $w \in W$, the 1-dimensional extremal
weight space of $S(\lambda)$ with weight $w(\lambda)$ is spanned
by 
$$
\Delta_{w(\lambda)} := 
\prod_{j=1}^n \Delta_{\vpi_j,w(\vpi_j)}^{\lambda(\alpha_j^\vee)}.
$$ 

Now consider the restriction map 
\[
\rho\df \C[N_-\backslash G']_{\rm s.r.} \to \C[\Nmin]_{\rm s.r.}
\]
given by restriction of functions from $G'$ to $\Nmin$.

\begin{Lem}
For every $\lambda\in \overline{P}^+$, the restriction 
$$
\rho_\lambda\df S(\lambda) \to \C[N^{\min}]_{\rm s.r.}
$$ 
of $\rho$ to $S(\lambda)$ is injective.
\end{Lem}

\begin{proof}
Let $B_-'$ be the Borel subgroup of $G'$ with unipotent radical
$N_-$.
We have 
$$
\Nmin \subset G_0 \cap G' = B_-'\Nmin.
$$
It follows that the natural projection from $G'$ onto
$B_-'\backslash G'$ restricts to an embedding of $\Nmin$,
with image the open subset of the flag variety 
$\X = B_-'\backslash G'$ defined by the non-vanishing of the minors
$\Delta_{\vpi_j,\vpi_j}$.
Now $\C[N_-\backslash G']_{\rm s.r.}$ can be regarded as the 
multi-homogeneous coordinate ring of $\X$ with homogeneous
components $S(\lambda)$, where $\lambda$ runs through $\overline{P}^+$.
It follows that $\C[\Nmin]$ can be identified with the
subring of degree 0 homogeneous elements of the localized
ring obtained from $\C[N_-\backslash G']_{\rm s.r.}$
by formally inverting the element 
$$
\Delta := \prod_{j=1}^n \Delta_{\vpi_j,\vpi_j}.
$$
Therefore, the restriction $\rho_\lambda$ of $\rho$ to every 
homogeneous piece $S(\lambda)$ is an embedding.
\end{proof}

It follows that we can transport the $G'$-module structure from
$S(\lambda)$ to $\rho(S(\lambda))$ by setting
\[
g\cdot \vph := \rho(g\cdot \rho_\lambda^{-1}(\vph)),
\qquad
(g \in G',\ \vph\in \rho(S(\lambda))).
\]
In this way, we can identify the highest weight module $L(\lambda)$
with the subspace $\rho(S(\lambda))$ of $\C[\Nmin]_{\rm s.r.}$. 
The highest weight vector is now $\rho(\Delta_\lambda)=1$, and the 
extremal weight vectors are the (restricted) generalized minors 
$$
D_{w(\lambda)} := \prod_{j=1}^n D_{\vpi_j,w(\vpi_j)}^{\lambda(\alpha_j^\vee)},
$$
for $w\in W$. 

At this point, we note that a strongly regular function on $\Nmin$
is just the same as an element of $U(\n)^*_{\rm gr}$.
Indeed, the elements of $\C[N^{\min}]_{\rm s.r.}$ are the restrictions
to $\Nmin$ of the linear combinations of matrix coefficients 
of the irreducible integrable representations 
$L(\lambda)$ with $\lambda \in \overline{P}^+$ of $G'$,
see \cite[Lemma 4.2]{KP1}.
Now, by Theorem~\ref{thm:irred-l}, we can realize every  
$L(\lambda)$ as a subspace of $U(\n)^*_{\rm gr}$, and every 
$f\in U(\n)^*_{\rm gr}$ belongs to such a subspace for 
$\lambda=\sum_{i=1}^nl_i\vpi_i$ with the $l_i\gg 0$.
It follows that each element of $U(\n)^*_{\rm gr}$
can be seen as a matrix coefficient for some   
$L(\lambda)$, and vice versa.
We can therefore identify
\[
\C[\Nmin]_{\rm s.r.} \equiv U(\n)^*_{\rm gr} \equiv \C[N].
\]
Moreover, these two ways of embedding $L(\lambda)$ in $\C[N]$ coincide.

\begin{Lem}\label{embedding2}
Let $\lambda \in \overline{P}^+$.
Under the identification $U(\n)^*_{\rm gr} \equiv \C[\Nmin]_{\rm s.r.}$,
the subspace $U(\lambda)$ defined in
Theorem~\ref{thm:irred-l}
coincides
with $\rho(S(\lambda))$.
\end{Lem}

\begin{proof}
The natural right action of $U(\n)$ on $U(\n)^*_{\rm gr}$ defined
before Corollary~\ref{embedding1} coincides with the
right action of $U(\n)$ on $\C[\Nmin]_{\rm s.r.}$ obtained by differentiating
the right regular representation of $\Nmin$:
\[
(f\cdot n)(x) = f(nx),\qquad (x,n\in\Nmin,\ f\in\C[\Nmin]_{\rm s.r.}).
\]
Consider first the case of a fundamental weight  $\lambda=\vpi_j$.
It is easy to check that 
\[
\Delta_{\vpi_j,\vpi_j}(x_i(t)g)=
\left\{
\begin{array}{ll}
\Delta_{\vpi_j,\vpi_j}(g)&
\mbox{if $i\not = j$}, \\
\Delta_{\vpi_j,\vpi_j}(g)
+t \Delta_{\vpi_j,\vpi_j}(\overline{s}_j g)&
\mbox{if $i = j$}.
\end{array}
\right.
\]
Now, the subspace $\rho(S(\lambda))$ is spanned by the functions
$n \mapsto \Delta_{\vpi_j,\vpi_j}(ng),\ (n \in N_-, \ g \in G')$.
By differentiating the previous equation with respect to $t$
and setting $t=0$, we obtain that 
\[
\rho(S(\lambda)) \subseteq 
\left\{ f\in \C[\Nmin]_{\rm s.r.}\mid f\cdot e_i = 0 
\text{ for } i\not =j, \ 
f\cdot e_j^2 = 0 \right\}.
\] 
Hence, using Corollary~\ref{embedding1}, we see that 
$\rho(S(\lambda))$ is contained in the embedding of $L(\vpi_j)$ into the
dual Verma module $M_{\rm low}^*(\vpi_j)$. 
Since these spaces have
the same graded dimensions, they must coincide.
The case of a general $\lambda \in \overline{P}^+$ 
follows using the fact that
$$
\Delta_\lambda = \prod_{j=1}^n \Delta_{\vpi_j,\vpi_j}^{\lambda(\alpha_j^\vee)}
$$
and that the $e_i$'s act as derivations on $\C[\Nmin]_{\rm s.r.}$.
\end{proof}


\section{The coordinate rings $\C[N(w)]$ and $\C[N^w]$}
\label{coordinaterings}


\subsection{The coordinate ring $\C[N(w)]$ as a ring of invariants}\label{invariant_ring}

Again, we fix a reduced expression $\ii = (i_r,\ldots,i_1)$ of
a Weyl group element $w$.
Assume that
$$
\cP = \{ p_{\bf m} \mid {\bf m} \in \N^{(J)}_{\phantom{0}} \}
$$ 
is an $\ii$-compatible PBW-basis of $U(\n)$.
Note that this PBW-basis of $U(\n)$ and also the corresponding dual PBW-basis 
of $U(\n)_{\rm gr}^*$
are homogeneous with respect to the (root lattice) $\N^n$-grading  of $U(\n)$. 
We write $\abs{\mm} = d \in \N^n$ in case $p_\mm$ is a homogeneous element
of degree $d \in \N^n$.
Let us denote by
$(\ee_i)_{i\in J}$ the usual coordinate vectors of $\Z^{(J)}_{\phantom{0}}$.
For example, $\abs{\ee_k}=\beta_\ii(k)$ for $1 \le k \le r$. 

The multiplication $\mu\df U(\n)\otimes U(\n) \to U(\n)$ is given
by its effect on the PBW-basis, say
$$
p_\mm \cdot p_{\nn} = 
\sum_{\abs{\kk} = \abs{\mm+\nn}} \!c_{\mm,\nn}^{\kk}\, p_{\kk}.
$$

Next, the comultiplication $\mu^*\df \C[N] \to \C[N] \otimes \C[N]$ is a
ring homomorphism, so it is determined by the value on the generators
$p^*_i = p^*_{\ee_i}$. 
By construction, we have
\[
\mu^*(p^*_i) = \sum_{\abs{\mm+\nn}=|\ee_i|} 
c_{\mm,\nn}^{\ee_i}\left(p^*_\mm\otimes p^*_{\nn}\right)
\]

\begin{Lem} \label{lem:coord}
Let $1 \le i \le r$ and $0 \neq \nn \in \N^{(J)}_{\phantom{0}}$ such
that $n_j=0$ for $1\leq j\leq r$. Then $c_{\mm,\nn}^{\ee_i} = 0$.
\end{Lem}

\begin{proof}
Let $\mm=\mm^<+\mm^>$ such that $m^<_j=0$ for $j>r$ and $m^>_j=0$ 
for $1\le j\leq r$, so $p_{\mm}=p_{\mm^<}\cdot p_{\mm^>}$.  
Since $\Delta_w^+$ is bracket closed and coclosed
we have
\[
p_{\mm^>}\cdot p_\nn=\sum_{\abs{\kk'}=\abs{\mm^>+\nn}} c_{\mm^>,\nn}^{\kk'} p_{\kk'}
\]
with $k'_j=0$ for $1\leq j\leq r$.
Thus
\[
p_\mm\cdot p_{\nn}=\sum_{\abs{\kk'}=\abs{\mm^>+\nn}}
c_{\mm^>,\nn}^{\kk'} p_{\kk'+\mm^<}.
\]
Putting $\kk=\kk'+\mm^<$ we get 
$c_{\mm,\nn}^\kk=c_{\mm^>,\nn}^{\kk'}$.
Thus, if in our situation $c_{\mm,\nn}^{\kk}\neq 0$ then
$k_j\neq 0$ for some $k>r$.
\end{proof}

Now, let us turn to the subgroups $N(w)$ and $N'(w)$. 
Consider the
ideals 
\[
I(w) := (p^*_{r+1},p^*_{r+2},\ldots ),\quad I'(w) := (p^*_1,\ldots,p^*_r)
\]
in $\C[N]$. 
Then we have
\begin{align*}
N(w) & =\{\nu\in\Hom_{\text{alg}}(\C[N],\C)\mid \nu(I(w))=0\},\text{ and}\\
N'(w)& =\{\nu'\in\Hom_{\text{alg}}(\C[N],\C)\mid \nu'(I'(w))=0\}.
\end{align*}
In other words we have canonically $\C[N(w)]=\C[N]/I(w)$ and 
$\C[N'(w)]=\C[N]/I'(w)$.

We consider the action of $N'(w)$ on $N$ via {\em right} multiplication. 
By definition, this comes from the 
left action of $N'(w)$ on $\C[N]$ given by
\[
\nu'\cdot f=(\text{id}\otimes\nu')\mu^*(f)
\]
for $f\in\C[N]$ and $\nu'\in N'(w)$. 
(Here we identify $\C[N]\otimes \C \equiv \C[N]$ in the canonical way.)

We denote by $\C[N]^{N'(w)}$ the 
invariant subring for this group action.

\begin{Prop}\label{prop17}
Consider the injective ring homomorphism 
$$
\tilde{\pi}^*_w\df \C[N(w)]\to \C[N]
$$
defined by $p^*_i+I(w)\mapsto p^*_i$ for $1\leq i\leq r$. 
The corresponding
morphism (of schemes) $\tilde{\pi}_w\df N\to N(w)$ 
is  $N'(w)$-invariant and is a
retraction for the inclusion of $N(w)$ into $N$. As a consequence,
$\tilde{\pi}^*_w$ identifies $\C[N(w)]$ with 
$\C[N]^{N'(w)}=\C[p_1^*,\ldots,p_r^*]$.
\end{Prop}

\begin{proof}
We have
\[
\mu^*(p^*_i)=1\otimes p^*_i + p^*_i\otimes 1+
\sum_{\abs{\mm+\nn}=\abs{\ee_i}} c_{\mm,\nn}^{\ee_i} 
\left(p^*_\mm\otimes p^*_\nn\right)
\]
where in the last sum $\abs{\mm}\neq 0\neq\abs{\nn}$.
Thus for $1\leq i\leq r$ and $\nu'\in N'(w)$ we get
\[\nu'\cdot p^*_i = 1\cdot 0 + p^*_i \cdot 1 +
\sum_{\abs{\mm+\nn}=\abs{\ee_i}} c_{\mm,\nn}^{\ee_i} p^*_\mm\cdot\nu'(p^*_\nn)
\]
with the last sum vanishing by Lemma~\ref{lem:coord} and
the definition of $N'(w)$. 
In other words, $p^*_i \in \C[N]^{N'(w)}$ for $1 \le i \le r$.
Thus, $\tilde{\pi}_w\df N \to N(w)$ is $N'(w)$-invariant,
that is, $\tilde{\pi}_w(nn') = \tilde{\pi}_w(n)$ for any $n'\in N'(w)$.

Now, since the multiplication map $N(w)\times N'(w) \to N$ is bijective,
each $N'(w)$-orbit on $N$ is of the form $n\cdot N'(w)$ for 
a unique $n\in N(w)$.
We conclude that the inclusion $N(w)\hookrightarrow N$ is a section for
$\tilde{\pi}_w$. 
Our claim follows.
\end{proof}

\subsection{The coordinate ring $\C[N^w]$ as a localization 
of $\C[N]^{N'(w)}$}\label{coordinateunipotent}
Let us now consider the groups $N(w)$ and $N'(w)$ introduced in 
Section~\ref{ssec:nilp}.

\begin{Lem}
We have 
\begin{align*}
N(w) &= N\cap (w^{-1} N_- w),\\
N'(w) &= N \cap (w^{-1}N w),\\
N'(w) \cap \Nmin &= \Nmin \cap (w^{-1}\Nmin w).
\end{align*}
\end{Lem}

\begin{proof}
This follows from \cite[5.2.3]{Ku} and \cite[6.2.8]{Ku}.
\end{proof}

It follows that $\Delta_{\vpi_j,w^{-1}(\vpi_j)}$ is invariant under the 
action of
$N'(w) \cap \Nmin$ on $\Gmin$ via right multiplication. 
Indeed, for $g\in \Gmin$ and $n'\in N'(w) \cap \Nmin$, we have
$n'\overline{w}^{-1}=\overline{w}^{-1}n''$ for some $n''\in N'(w) \cap \Nmin$,
hence
\begin{align*}
\Delta_{\vpi_j,w^{-1}(\vpi_j)}(gn') &= 
\Delta_{\vpi_j,\vpi_j}(gn'\overline{w}^{-1})=
\Delta_{\vpi_j,\vpi_j}(g\overline{w}^{-1}n'')\\
&=\Delta_{\vpi_j,\vpi_j}(g\overline{w}^{-1})
=\Delta_{\vpi_j,w^{-1}(\vpi_j)}(g).
\end{align*}
Define
\[
O_w := \left\{ n \in \Nmin \mid \Delta_{\vpi_j,w^{-1}(\vpi_j)}(n) \neq 0 
\text{ for all } 1 \le j \le n \right\}.
\]
This is the open subset of $\Nmin$ consisting of elements $n$ such
that $\overline{w}n^T\in G_0$.
Indeed,
\[
\Delta_{\vpi_j,w^{-1}(\vpi_j)}(n)=
\Delta_{\vpi_j,\vpi_j}(n\overline{w}^{-1})=
\Delta_{\vpi_j,\vpi_j}((n\overline{w}^{-1})^T)=
\Delta_{\vpi_j,\vpi_j}(\overline{w}n^T),
\]
since $\overline{w}^{-1}=\overline{w}^T$.
Following \cite[Section 5]{BZ}, we can now define the map  
$\tilde{\eta}_w\df O_w\to \Nmin$
given by 
\[
\tilde{\eta}_w(z) := [\overline{w}z^T]_+.
\]
Recall that $N^w = N \cap (B_-wB_-)$, see Section~\ref{introbla}.

\begin{Prop}\label{BZprop}
The following properties hold:
\begin{itemize}

\item[(i)] The map $\tilde{\eta}_w$ is a morphism of ind-varieties.

\item[(ii)] The image of $\tilde{\eta}_w$ is $N^w$.

\item[(iii)] $\tilde{\eta}_w(x) = \tilde{\eta}_w(y)$ if and only if 
$x=yn'$ for some $n'\in N'(w)\cap \Nmin$.

\item[(iv)] $\tilde{\eta}_w$ restricts to a bijective 
morphism $N(w)\cap O_w \to N^w$.

\item[(v)] We have $N^w\subset O_w$, and 
$\tilde{\eta}_w$ restricts to a bijection 
$\eta_w\df N^w \to N^w$.

\item[(vi)] The inverse of $\eta_w$ is given by
$\eta_w^{-1}(x) = \eta_{w^{-1}}(x^\iota)^\iota$ for $x \in N^w$.
It follows that $\eta_w$ is an automorphism of $N^w$.

\end{itemize}
\end{Prop}

\begin{proof}
Property (i) follows from Proposition~\ref{Gauss}~(ii).
Next, we have
\[
[\overline{w}z^T]_+=([\overline{w}z^T]_0^{-1}
[\overline{w}z^T]_-^{-1})\overline{w}z^T \in B_-\overline{w}B_{-}.
\]
This shows that the image of $\tilde{\eta}_w$ is contained in $N^w$.
The rest of (ii) and (iii) is proved as in \cite[Proposition 5.1]{BZ}.
Property~(iv) follows from (ii), (iii), and the decomposition
$\Nmin = N(w) \times (N'(w)\cap \Nmin)$. 
Finally, (v) and (vi) are proved exactly in the same way as
in \cite[Propositions 5.1, 5.2]{BZ}.  
\end{proof}

\begin{Prop}\label{coordNw}
The map $\tilde{\pi}_w$ restricts to a morphism   
$\pi_w\df N^w\to O_w\cap N(w)$.
This is an isomorphism with inverse
$$
\eta^{-1}_w\tilde{\eta}_w\df O_w\cap N(w)\to N^w.
$$
In particular, $N^w$ is an affine
variety with coordinate ring identified to the localized
ring $\C[N]^{N'(w)}_{\Delta_w}$, 
where 
$$
\Delta_w := \prod_{j=1}^n \Delta_{\vpi_j,w^{-1}(\vpi_j)}.
$$
\end{Prop}

\begin{proof}
By Proposition~\ref{BZprop}~(iv) and (v), we know that 
$\eta^{-1}_w\tilde{\eta}_w$ is a bijection.
On the other hand $\tilde{\pi}_w(N^w) \subseteq O_w \cap N(w)$
because $N^w\subset O_w$.
Now, by Proposition~\ref{BZprop}~(iii), we have
$$
\tilde{\eta}_w(\pi_w(x)) = \tilde{\eta}_w(x) = \eta_w(x)
$$
for every $x\in N^w$. Hence
$\eta^{-1}_w\tilde{\eta}_w\pi_w(x)=x$ for every $x$ in
$N^w$. 
So we have $\eta^{-1}_w\tilde{\eta}_w\pi_w = \id_{{N^w}}$,
and this proves that $\pi_w$ is the inverse of 
$\eta^{-1}_w\tilde{\eta}_w$.

These maps are morphisms of varieties so they induce 
isomorphisms 
\[
\C[N^w] \xrightarrow{\sim} \C[N(w)\cap O_w]=\C[N(w)]_{\Delta_w}
\xrightarrow{\sim} 
\C[N]_{\Delta_w}^{N'(w)}.
\]
\end{proof}

The following commutative diagram displays the different morphisms
appearing in Propositions~\ref{BZprop} and \ref{coordNw}:
$$
\xymatrix{
N \ar[r]^>>>>>{\tilde{\pi}_w} & N(w) && 
O_w \ar[r]^>>>>>{\tilde{\eta}_w} & \Nmin\\
& N^w \ar[r]_>>>>>{\pi_w}^>>>>>>\sim \ar[ul]_>>>>>>>>>>\iota 
\ar@/_2pc/[rr]_{\eta_w} & 
N(w) \cap O_w \ar[r]^>>>>>\sim
\ar[ul]_>>>>>>>>>\iota \ar[ur]^>>>>>>>>>>\iota & N^w \ar[ur]^>>>>>>>>>\iota
}
$$
(The arrows labelled with $\iota$ are inclusion maps.)


\section{The modules $V_k$ and $M_k$}


For the entire section, we fix a reduced expression
$\ii = (i_r,\ldots,i_1)$ of a Weyl group element $w$, and
as before let
$V_\ii = V_1 \oplus \cdots \oplus V_r$ and 
$M_\ii = M_1 \oplus \cdots \oplus M_r$.
Recall that for each $1 \le k \le r$ there is a short exact sequence
$$
0 \to V_{k^-} \to V_k \to M_k \to 0
$$
of $\LL$-modules.

\subsection{Generalized minors as $\varphi$-functions}
\label{catminors}
For $1 \le k \le r$ set
$$
w^{-1}_{\le k} := s_{i_1}\cdots s_{i_k}.
$$

\begin{Prop}\label{minors}
For $1 \le k \le r$ we have
\[
\vph_{V_k}=
D_{\vpi_{i_k},w^{-1}_{\le k}(\vpi_{i_k})}.
\]
In particular, we have $\vph_{I_{\ii,j}} = D_{\vpi_j,w^{-1}(\vpi_j)}$
for every $1 \le j \le n$.
\end{Prop}

\begin{proof} 
Using Lemma~\ref{embedding2}, we can realize the fundamental module
$L(\vpi_{i_k})$ as the subspace $\rho(S(\vpi_{i_k}))$ of $\C[N]$.
Then using Theorem~\ref{thm:irred-l}, the definition of $V_k$ 
(see Section~\ref{terminal}) 
and the discussion in 
Section~\ref{genminors},
we can check that the function $\vph_{V_k}$ is an extremal weight vector
of weight $w^{-1}_{\le k}(\vpi_{i_k})$ in $L(\vpi_{i_k})$, hence
it coincides with $D_{\vpi_{i_k},w^{-1}_{\le k}(\vpi_{i_k})}$ up to a scalar.
Moreover, its image under $e_{i_k}^{\rm max}\cdots e_{i_1}^{\rm max}$
is equal to $1$, so the normalizations agree and we have
$\vph_{V_k} = D_{\vpi_{i_k},w^{-1}_{\le k}(\vpi_{i_k})}$. 
\end{proof}

\begin{Cor}\label{dimVk2}
For $1 \le k \le r$ we have
$\dimv(V_k) = \vpi_{i_k} - s_{i_1}s_{i_2} \cdots s_{i_k}(\vpi_{i_k})$.
\end{Cor}

\begin{proof}
The statement follows from the following general fact:
Assume that $\delta_X \in U(\lambda)$ for some weight 
$\lambda \in P^+$ and
some $\LL$-module $X$.
When we consider $\delta_X$ 
as an element of $L(\lambda) \equiv U(\lambda)$,
Theorem~\ref{thm:irred-l} implies that
${\rm wt}(\delta_X) = \lambda - \dimv(X)$.
\end{proof}

Recall that for 
$1 \le k \le r$ we defined
$$
\beta_\ii(k) = 
\begin{cases}
\alpha_{i_1} & \text{if $k=1$},\\
s_{i_1} \cdots s_{i_{k-1}}(\alpha_{i_k}) & \text{otherwise}.
\end{cases}
$$

\begin{Cor}
For $1 \le k \le r$ we have
$\dimv(M_k) = \beta_\ii(k)$.
\end{Cor}

\begin{proof}
By Corollary~\ref{dimVk2}
we know that
$\dimv(V_k) = \vpi_{i_k} - s_{i_1}s_{i_2} \cdots s_{i_k}(\vpi_{i_k})$
for each $1 \le k \le r$.
By the definition of $M_k$ we have
\begin{align*}
\dimv(M_k) &= \dimv(V_k) - \dimv(V_{k^-}) \\
&= s_{i_1}s_{i_2} \cdots s_{i_{k^-}}(\vpi_{i_k}) - 
s_{i_1}s_{i_2} \cdots s_{i_k}(\vpi_{i_k}) \\
&= s_{i_1}s_{i_2} \cdots s_{i_{k^-}}\left(\vpi_{i_k} - 
s_{i_{k^-+1}} \cdots s_{i_k}(\vpi_{i_k})\right).
\end{align*}
But 
$$
s_j(\vpi_{i_k}) = 
\begin{cases}
\vpi_{i_k} & \text{if $j \not= i_k$},\\ 
\vpi_{i_k}-\alpha_{i_k} & \text{if $j = i_k$}.
\end{cases}
$$
It follows that
\begin{align*}
\dimv(M_k) &= s_{i_1}s_{i_2} \cdots s_{i_{k^-}}\left(\vpi_{i_k} - 
\vpi_{i_k} + s_{i_{k^-+1}} \cdots s_{i_{k-1}}(\alpha_{i_k})\right) \\
&= s_{i_1}s_{i_2}\cdots s_{i_{k-1}}(\alpha_{i_k}).
\end{align*}
This finishes the proof.
\end{proof}

\begin{Cor}\label{bracketclosed2}
We have
$\Delta_w^+ = \{ \dimv(M_1),\ldots,\dimv(M_r) \}$.
\end{Cor}

\subsection{Example}
Let $Q$ be a quiver with underlying graph
$$
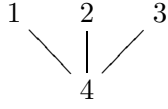
\xymatrix@-0.7pc{
1 \ar@{-}[dr] & 2\ar@{-}[d] & 3 \ar@{-}[dl]\\
&4
}
$$
Let $w$ be the Weyl group element
$s_3s_4s_2s_1s_4$.
The set of reduced expressions for $w$ is
$R(w) = \{ (3,4,2,1,4), (3,4,1,2,4) \}$.
We have
$$
\Delta_w^+ = \left\{ 
\bsm 0&0&0\\&1\esm, \bsm 1&0&0\\&1\esm,\bsm 0&1&0\\&1\esm,
\bsm 1&1&0\\&1\esm,\bsm 1&1&1\\&2\esm
 \right\}.
$$
Let $\ii=(3,4,2,1,4)$.
We get
$$
V_\ii = V_1 \oplus \cdots \oplus V_5
=
\bsm 4 \esm \oplus 
\bsm &4\\1 \esm \oplus
\bsm 4\\2 \esm \oplus
\bsm &4\\1&2\\&4 \esm \oplus 
\bsm &4\\1&2\\&4\\&&3\esm
$$
and
$$
M_\ii = M_1 \oplus \cdots \oplus M_5
=
\bsm 4 \esm \oplus 
\bsm &4\\1 \esm \oplus
\bsm 4\\2 \esm \oplus
\bsm &4\\1&2 \esm \oplus 
\bsm &4\\1&2\\&4\\&&3\esm.
$$
Note that $\add(M_\ii)$ is neither closed under factor modules
nor under submodules.
We have
$$
\CC_w = \add\left(V_\ii \oplus \bsm &4\\1&2 \esm\right).
$$
We can think of $\CC_w$ as a categorification of a cluster
algebra of type $\A_1$ with four coefficients.

\subsection{Example}\label{sect22_5}
Let $Q$ be a quiver with underlying graph 
$\xymatrix@-0.5pc{1 \ar@<0.5ex>@{-}[r]\ar@<-0.5ex>@{-}[r] & 2 \ar@{-}[r]& 3}$
Then $\ii := (i_7,\ldots,i_1) := (3,1,2,3,1,2,1)$ is a reduced expression
of a Weyl group element $w \in W_Q$.
The indecomposable direct summands of $V_\ii$ are
\begin{align*}
V_1 &= {\bsm1\esm} &
V_2 &= {\bsm1&&1\\&2\esm} &
V_3 &= {\bsm1&&1&&1\\&2&&2\\&&1\esm} \\
V_4 &= {\bsm1&&1\\&2\\&3\esm} &
V_5 &= 
{\bsm&\\&\\1&&1&&1&&1&&1&&1\\&2&&2&&&&2&&2\\&&1&&&3&&&1\\&&&&&2\esm} \\
V_6 &= 
{\bsm&\\&\\1&&1&&1&&1&&1&&1&&1&&1&&1\\&2&&2&&&&2&&2&&&&2&&2\\
&&1&&&3&&&1&&&3&&&1
\\&&&&&2&&&&&&2\\&&&&&&&&1\esm} &
V_7 &=
 {\bsm&\\&\\1&&1&&1&&1\\&2&&2&&2\\&&1&&1\\&&&2\\&&&3\esm}.
\end{align*}
Here, the $\LL$-modules are represented by their socle filtration.
The indecomposable $\CC_w$-projective-injective modules 
are $V_5, V_6$ and $V_7$.
The corresponding functions $\vph_{V_k}$ are given by
\begin{align*}
\vph_{V_1} &= D_{\vpi_1,s_1(\vpi_1)} &
\vph_{V_2} &= D_{\vpi_2,s_1s_2(\vpi_2)} &
\vph_{V_3} &= D_{\vpi_1,s_1s_2s_1(\vpi_1)} \\
\vph_{V_4} &= D_{\vpi_3,s_1s_2s_1s_3(\vpi_3)} &
\vph_{V_5} &= D_{\vpi_2,s_1s_2s_1s_3s_2(\vpi_2)} \\
\vph_{V_6} &= D_{\vpi_1,s_1s_2s_1s_3s_2s_1(\vpi_1)} &
\vph_{V_7} &= D_{\vpi_3,s_1s_2s_1s_3s_2s_1s_3(\vpi_3)}.
\end{align*}
%

\subsection{Example}
We continue to discuss the example from Section~\ref{exampleintro}.
Thus $Q$ is a quiver with underlying graph
$\xymatrix@-0.5pc{1 \ar@{-}[r] & 2 \ar@{-}[r] & 3 \ar@{-}[r] &4}$.
Note that the Weyl group $W_Q$ is the symmetric group $S_5$,
and the generators $s_i$ are the transpositions $(i,i+1)$.
The generalized minors become ordinary minors. 
More precisely, for $w \in S_5$ and $i \in \{1,2,3,4,5\}$ we have 
\[
\Delta_{\varpi_i,w(\varpi_i)} = 
\Delta_{\{1,2,\ldots,i\},w(\{1,2,\ldots,i\})},
\]
since we may identify $S_5$ with the group of permutation matrices in $\GL_5$.
Here $\Delta_{I,J}$ denotes the minor in $\C[{\rm SL}_5]$ with
row set $I$ and column set $J$.
As in Section~\ref{exampleintro} let
$w := s_3s_4s_2s_1s_3s_4s_2s_1$ and
$\ii := (i_8,\ldots,i_1) := (3,4,2,1,3,4,2,1)$.
We get
\begin{multline*}
x_\ii(t) :=
x_3(t_8) x_4(t_7) x_2(t_6) x_1(t_5) x_3(t_4) x_4(t_3) x_2(t_2) x_1(t_1) = \\
= \left(\begin{array}{ccccc}
1 &t_5+t_1&t_5t_2&0      &0\\
0 &1      &t_6+t_2&t_6 t_4&t_6 t_4 t_3\\
0 &0      &1      &t_8+t_4&t_8(t_7 +t_3)+t_4 t_3\\
0 &0      &0      &1      &t_7+t_3\\
0 &0      &0      &0      &1
\end{array}\right).
\end{multline*}
A straightforward calculation shows:
\begin{align*}
D_{\varpi_1,w_{\leq 1}^{-1}(\varpi_1)}&=
D_{\{1\},\{2\}} = t_5+t_1,\\
D_{\varpi_2,w_{\leq 2}^{-1}(\varpi_2)}&=
D_{\{1,2\},\{2,3\}} = t_6(t_5+t_1)+t_2t_1,\\
D_{\varpi_4,w_{\leq 3}^{-1}(\varpi_4)}&=
D_{\{1,2,3,4\},\{1,2,3,5\}} = t_7+t_3,\\
D_{\varpi_3,w_{\leq 4}^{-1}(\varpi_3)}&=
D_{\{1,2,3\},\{2,3,5\}} =
t_8(t_7(t_6(t_5+t_1)+t_2t_1)+t_6t_3(t_5+t_1)+t_3t_2t_1)\\
&\quad\quad\quad\quad\quad\quad\quad\quad+t_4t_3t_2t_1,\\
D_{\varpi_1,w_{\leq 5}^{-1}(\varpi_1)}&=
D_{\{1\},\{3\}} = t_5t_2,\\
D_{\varpi_2,w_{\leq 6}^{-1}(\varpi_2)}&=
D_{\{1,2\},\{3,5\}} = t_6t_5t_4t_3t_2,\\
D_{\varpi_4,w_{\leq 7}^{-1}(\varpi_4)}&=
D_{\{1,2,3,4\},\{2,3,4,5\}} = t_7t_4t_2t_1,\\
D_{\varpi_3,w_{\leq 8}^{-1}(\varpi_3)}&=
D_{\{1,2,3\},\{3,4,5\}} = t_8t_7t_6t_5t_4t_2.
\end{align*}
Here
the evaluation of the minors is always on 
$x_\ii(t)$.
Due to the structure of the modules $V_k$ described in 
Section~\ref{exampleintro}, 
we could also use
Proposition~\ref{phi_form} and calculate directly that
$$
\varphi_{V_k}(x_\ii(t))=
D_{\varpi_{i_k},w_{\le k}^{-1}(\varpi_{i_k})}(x_\ii(t))
$$
for all $1 \le k \le 8$.

\subsection{Refined socle and top series}\label{refined}
For any $\LL$-module $X \in \CC_w$ there exists
a unique chain
$$
0 = X_r \subseteq \cdots \subseteq X_1 \subseteq X_0 = X
$$
of submodules of $X$ such that
$X_{k-1}/X_k = \soc_{S_{i_k}}(X/X_k)$.
This is called the {\it refined socle series} of type $\ii$ of $X$.
Define 
$$
{\bf s}_\ii(X) := (p_r,\ldots,p_1)
$$
where $p_k := \dim(X_{k-1}/X_k)$ for $1 \le k \le r$.
Similarly, there exists a unique chain
$$
0 = Y_r \subseteq \cdots \subseteq Y_1 \subseteq Y_0 = X
$$
of submodules of $X$ 
such that
$Y_{k-1}/Y_k = \tp_{S_{i_k}}(Y_{k-1})$
for all $1 \le k \le r$.
This is called the {\it refined top series} of type $\ii$ of $X$.
Define 
$$
{\bf t}_\ii(X) := (q_r,\ldots,q_1)
$$
where $q_k := \dim(Y_{k-1}/Y_k)$ for $1 \le k \le r$.
(For a simple module $S$ and a module $M$ let $\tp_S(M)$ 
be the intersection of all submodules $U$ of $M$ with $M/U \cong S$.)

The existence of refined socle and top series of type $\ii$
of $X \in \CC_w$ comes
from the fact that $V_\ii$ generates the category $\CC_w$.
It follows directly from the definitions that
each module $V_k$ has a refined socle and top series of type $\ii$.
Now one easily checks that this property also holds for
factor modules of modules in $\add(V_\ii)$.

The uniqueness of refined socle and top series of type $\ii$
implies the following result:

\begin{Lem}
Let $\ii = (i_r,\ldots,i_1)$ be a reduced expression of $w$,
and let $X \in \CC_w$.
Set ${\bf s} := {\bf s}_\ii(X) = (p_r,\ldots,p_1)$ and
${\bf t} := {\bf t}_\ii(X) = (q_r,\ldots,q_1)$.
Then the following hold:
\begin{itemize}

\item[(i)]
We have
$$
\F_{\ii^{\bf s},X} \cong \prod_{k=1}^r \F\left(\C^{p_k}\right)
\text{\;\;\; and \;\;\;}
\F_{\ii^{\bf t},X} \cong \prod_{k=1}^r \F\left(\C^{q_k}\right).
$$
In particular,
$$
\chi_{\rm c}(\F_{\ii^{\bf s},X}) = \prod_{k=1}^r p_k!
\text{\;\;\; and \;\;\;}
\chi_{\rm c}(\F_{\ii^{\bf t},X}) = \prod_{k=1}^r q_k!.
$$

\item[(ii)]
$\F_{\ii,{\bf s},X}$ and $\F_{\ii,{\bf t},X}$ both consist of a single point.
In particular,
$\chi_{\rm c}(\F_{\ii,{\bf s},X}) = 1$
and
$\chi_{\rm c}(\F_{\ii,{\bf t},X}) = 1$.

\end{itemize}
\end{Lem}

Observe that
$(i_k,\ldots,i_t)$ is a reduced expression for the Weyl group element
$w_{k,t} := s_{i_k} s_{i_{k-1}} \cdots s_{i_t}$ for all $1 \le t \le k \le r$.
Set $\jj := (i_r,\ldots,i_2)$.
For $1 \le k \le r$ define 
$$
b_k := b_{\ii,k} := -(s_{i_k} \cdots s_{i_r}(\vpi_{i_r}),\alpha_{i_k})
= (s_{i_{k+1}} \cdots s_{i_r}(\vpi_{i_r}),\alpha_{i_k}),
$$
and set $\bb_\ii := (b_r,\ldots,b_1)$.

\begin{Prop}\label{Vk1}
For $\ii$ and $\jj$ as above, the following hold:
\begin{itemize}

\item[(i)]
$\tp_{S_{i_1}}(V_{\jj,r-1}) = 0$;

\item[(ii)]
$\tp_{S_{i_1}}(V_{\ii,r}) = S_{i_1}^{b_1}$;

\item[(iii)]
$\sss_\ii(V_{\ii,r}) = \ttt_\ii(V_{\ii,r}) = b_\ii$.

\end{itemize}

\end{Prop}

\begin{proof}
For $r=1$ the statements are obvious.
Thus assume $r \ge 2$.
Let $u_{\vpi_{i_r}}$ be a highest weight vector
in $L(\vpi_{i_r})$.
Since $\ii = (\jj,i_1)$ is a reduced expression, 
we know from Section~\ref{genminors} that
\begin{equation}\label{Eq5}
e_{i_1}(\overline{s}_{i_2} \cdots \overline{s}_{i_r}(u_{\vpi_{i_r}})) = 0.
\end{equation}
By Proposition~\ref{minors}
we can identify 
$\overline{s}_{i_2} \cdots \overline{s}_{i_r}(u_{\vpi_{i_r}})$ with
$\varphi_{V_{\jj,r-1}}$.
We have $\tp_{S_{i_1}}(V_{\jj,r-1}) \cong S_{i_1}^c$ for some $c \ge 0$.
Let
$U$ be the unique submodule
such that $V_{\jj,r-1}/U = \tp_{S_{i_1}}(V_{\jj,r-1})$.
We get 
$$
e_{i_1}^{(c)}\varphi_{V_{\jj,r-1}} = \varphi_U \not= 0.
$$
But if $c \ge 1$, then equation~(\ref{Eq5}) yields
$e_{i_1}^{(c)}\varphi_{V_{\jj,r-1}} = 0$, a contradiction.
This implies $c = 0$.
So we proved (i).
To show (ii) we use that $\varphi_{V_{\ii,r}}$ can be identified with
$$
\overline{s}_{i_1}\overline{s}_{i_2} \cdots 
\overline{s}_{i_r}(u_{\vpi_{i_r}}) = 
f_{i_1}^{\max}
\left(\overline{s}_{i_2} \cdots \overline{s}_{i_r}(u_{\vpi_{i_r}})\right) =
f_{i_1}^{\max}\left(\varphi_{V_{\jj,r-1}}\right).
$$
We have
${\rm wt}\left(\overline{s}_{i_1} \cdots 
\overline{s}_{i_r}(u_{\vpi_{i_r}})\right) 
= 
{\rm wt}\left(\overline{s}_{i_2} \cdots 
\overline{s}_{i_r}(u_{\vpi_{i_r}})\right) 
- b_1\alpha_{i_1}$,
see Section~\ref{genminors}.
This implies (ii).
Finally, it follows by induction on $r$ that
$\sss_\ii(V_{\ii,r}) = \ttt_\ii(V_{\ii,r}) = b_\ii$.
This finishes the proof.
\end{proof}

\subsection{Computation of the 
Euler characteristics $\chi_{\rm c}(\F_{\kk,V_k})$}
\label{Euler1}
By Proposition~\ref{phi_form}, to evaluate $\vph_{V_k}$
on $x_{j_1}(t_1)\cdots x_{j_p}(t_p)$, we need to know
the Euler characteristic $\chi_{\rm c}(\F_{\kk,V_k})$ for arbitrary
types $\kk$ of composition series.
These Euler characteristics can in turn be calculated via
a simple algorithm that we shall now describe.

To this end, it will be convenient to embed $U(\n)^*_{\rm gr} \equiv \C[N]$
in the shuffle algebra $F^*$, as explained in \cite[\S 2.8]{LeMZ}.
As a $\C$-vector space, $F^*$ has a basis consisting of all words
\[
w[\kk] := w[k_1,k_2,\ldots,k_s] := w_{k_1}w_{k_2}\cdots w_{k_s},
\qquad (1\le k_1,\ldots,k_s\le n,\ s \ge 0),
\]
in the letters $w_1,\ldots, w_n$.
The multiplication in $F^*$ is the classical 
commutative shuffle product $\shuffle$ of words 
with unit the empty word $w[]$, see \textit{e.g.}~\cite{Reu} and 
\cite[\S 2.5]{LeMZ}.
By \cite[Propositions 9 and 10]{LeMZ}, for any $X \in \nil(\LL)$ 
the image of $\vph_X$ in 
this embedding is just the {\it generating function} 
\[
g_X := \sum_{\kk} \chi_{\rm c}(\F_{\kk,X}) w[\kk]
\]
of the Euler characteristics $\chi_{\rm c}(\F_{\kk,X})$ for all types
$\kk$ of composition series. 
(The Euler characteristic
$\chi_{\rm c}(\F_{\kk,X})$ is equal to the coefficient of $t_1 \cdots t_s$ in
$\vph_X(x_{k_1}(t_1)\cdots x_{k_s}(t_s))$.)

Let $\la \in P^+$ and $1 \le i \le n$. 
Define endomorphisms $\rho_\la(e_i), \rho_\la(f_i)$ of the vector 
space $F^*$ by 
\begin{eqnarray*}
\rho_\la(e_i) (w[j_1,\ldots ,j_k]) &:=& \delta_{i,j_k} w[j_1, \ldots
,j_{k-1}],\label{eq2}\\[2mm]
\rho_\la(f_i) (w[j_1,\ldots ,j_k]) &:=& \sum_{l=0}^k 
(\la-\alpha_{j_1}-\cdots
-\alpha_{j_l})(\alpha_i^\vee)\,w[j_1,\ldots,j_l,i,j_{l+1},\ldots , j_k].
\end{eqnarray*}

\begin{Prop}\label{prop1}
The formulas above extend to a representation  
$\rho_\la\df U(\g) \to \End_\C(F^*)$ of $U(\g)$.
This turns $F^*$ into a $U(\g)$-module.
The image of $\C[N]$ in its embedding in $F^*$
is a $U(\g)$-submodule 
isomorphic to the dual Verma module $M_{\rm low}^*(\la)$,
see Section~\ref{dualVerma}.
In particular the set
$$
\left\{\rho_\la(f_{i_1}\cdots f_{i_s})(w[]) \mid 1\le i_1,\ldots ,i_s \le n,\;
s \ge 0 \right\}
$$
spans the irreducible module $L(\la)$, considered as
a submodule of $M_{\rm low}^*(\la)$.
\end{Prop}

The above formulas for $\rho_\la(e_i)$ and $\rho_\la(f_i)$
can be obtained by specializing $q$ to 1 in the formulas of the
proof of \cite[Proposition 50]{LeMZ}. 
We omit the details.

By Proposition~\ref{minors}, for $1 \le k \le r $ we have
\[
\vph_{V_k}=
D_{\vpi_{i_k},w^{-1}_{\le k}(\vpi_{i_k})}.
\]
By Section~\ref{genminors} we know that
$\vph_{V_k}$ is obtained by acting on the highest 
weight vector $u_{\vpi_{i_k}}$ 
of $L(\vpi_{i_k})$ with the product 
$f_{i_1}^{(b_1)}\cdots f_{i_k}^{(b_k)}$ of divided powers of the Chevalley 
generators, where $b_k = b_{\ii,k}$ is defined as in 
Section~\ref{refined}.
Therefore we have
\begin{equation}\label{calculEuler}
g_{V_k} = 
\rho_{\vpi_{i_k}}\left(f_{i_1}^{(b_1)}\cdots f_{i_k}^{(b_k)}\right)(w[]).
\end{equation}
Hence to calculate the generating function $g_{V_k}$ one only needs
to apply 
$b_1 + \cdots + b_k = \dim(V_k)$ times the above combinatorial formula 
for $\rho_{\varpi_{i_k}}(f_i)$.
Thus we have obtained an algorithm for calculating all 
Euler characteristics $\chi_{\rm c}(\F_{\kk,V_k})$.

\subsection{Example}
We continue the example of Section~\ref{sect22_5}.
Clearly, we have 
\[
g_{V_1}=\rho_{\vpi_1}(f_1)(w[])= \vpi_1(\alpha_1^\vee)w[1]=w[1].
\]
Similarly
\[
g_{V_2}=\rho_{\vpi_2}(f_1^{(2)}f_2)(w[]).
\]
Now we calculate successively
\[
\begin{array}{rcl}
\rho_{\vpi_2}(f_2)(w[])&=&\vpi_2(\alpha_2^\vee)\,w[2]\ =\ w[2],\\
\rho_{\vpi_2}(f_1)(w[2])&=&\vpi_2(\alpha_1^\vee)\,w[1,2]
+(\vpi_2-\alpha_2)(\alpha_1^\vee)\,w[2,1]
\ =\ 2\,w[2,1],\\
\rho_{\vpi_2}(f_1)(2\,w[2,1])&=&
2(\vpi_2(\alpha_1^\vee)\,w[1,2,1]
+(\vpi_2-\alpha_2)(\alpha_1^\vee)\,w[2,1,1]\\
&&+(\vpi_2-\alpha_2-\alpha_1)(\alpha_1^\vee)\,w[2,1,1])
\\
&=&4\,w[2,1,1].
\end{array}
\]
Hence, taking into account that $f_1^{(2)}=f_1^{2}{/}2$, we get
\[
g_{V_2}=2\,w[2,1,1].
\]
Similar applications of formula~(\ref{calculEuler}) yield the following results
\[
\begin{array}{rcl}
g_{V_3}&=& \rho_{\vpi_1}\left(f_1^{(3)}f_2^{(2)}f_1\right)(w[]) 
= 4\,w[1,2,1,2,1,1]+12\,w[1,2,2,1,1,1],\\
g_{V_4}&=& \rho_{\vpi_3}\left(f_1^{(2)}f_2f_3\right)(w[])
= 2\,w[3,2,1,1],\\
g_{V_7}&=&\rho_{\vpi_3}\left(f_1^{(4)}f_2^{(3)}f_1^{(2)}f_2f_3\right)(w[])\\
&=& 288\,w[3,2,1,1,2,2,2,1,1,1,1]+144\,w[3,2,1,1,2,2,1,2,1,1,1]\\
       && +96\,w[3,2,1,2,1,2,2,1,1,1,1]+48\,w[3,2,1,1,2,2,1,1,2,1,1]\\
       && +48\,w[3,2,1,2,1,1,2,2,1,1,1]+48\,w[3,2,1,2,1,2,1,2,1,1,1]\\
       && +48\,w[3,2,1,1,2,1,2,2,1,1,1]+16\,w[3,2,1,2,1,2,1,1,2,1,1]\\
       && +16\,w[3,2,1,2,1,1,2,1,2,1,1]+16\,w[3,2,1,1,2,1,2,1,2,1,1].
\end{array}
\]
The generating functions $g_{V_5}$ and $g_{V_6}$ are too large to be included
here. 
For example $g_{V_5}$ is a linear combination of 402 words.

\subsection{The modules $M[b,a]$}\label{intervalmod}
For $1 \le k \le r$ let
\begin{align*}
k^- &:= \max\{ 0,1 \le s \le k-1 \mid i_s = i_k \},\\
k^+ &:= \min\{ k+1 \le s \le r,r+1 \mid i_s = i_k \},\\
k_{\rm min} &:= \min\{ 1 \le s \le r \mid i_s = i_k \},\\
k_{\rm max} &:= \max\{ 1 \le s \le r \mid i_s = i_k \}.
\end{align*}
Set $k^{(0)} := k$, and for an integer $m$ define 
$k^{(m-1)} := (k^{(m)})^-$ and
$k^{(m+1)} := (k^{(m)})^+$.
For $1 \le j \le n$ and $1 \le k \le r+1$ let 
$$
k^-(j) := \max\{ 0,1 \le s \le k-1 \mid i_s = j \},
$$
and
$$
k[j] := |\{ 1 \le s \le k-1 \mid i_s = j \}|, 
$$
and set $t_j := (r+1)[j]$.

For $1 \le a \le b \le r$ with $i_a = i_b$
define $M[b,a] := V_b/V_{a^-}$.
(For convenience, we define $V_0 = V_{r+1} = 0$.)
We have a short exact sequence
$$
0 \to M[a^-,b_{\rm min}] \to M[b,b_{\rm min}] \to M[b,a] \to 0.
$$
Note that $a_{\rm min} = b_{\rm min}$, since we assume $i_a = i_b$.
For $1 \le k \le r$ we have $M[k,k_{\rm min}] = V_k$ and
$M[k,k] = M_k$.
One can visualize a module $M[b,a]$ by
$$
\begin{tabular}{c}
$M_b$\\
\hline
$M_{b^-}$\\
\hline
$\cdots$\\
\hline
$M_a$
\end{tabular}
$$
We have 
$$
V_\ii = \bigoplus_{k=1}^r M[k,k_{\rm min}].
$$
For each $k$ we have a short exact sequence
$$
0 \to M[k,k_{\rm min}] \to M[k_{\rm max},k_{\rm min}] \to 
M[k_{\rm max},k^+] \to 0.
$$
Note that $M[k_{\rm max},k_{\rm min}] = I_{\ii,i_k}$ 
is $\CC_w$-projective-injective.
Define 
$$
T_k := T_{\ii,k} :=
\begin{cases}
V_k & \text{if $k^+ = r+1$},\\
M[k_{\rm max},k^+] & \text{otherwise}.
\end{cases}
$$
Thus if $k^+ \not= r+1$, then $\Omega_w^{-1}(V_k) = T_k$.
Define $T_\ii := T_1 \oplus \cdots \oplus T_r$.
In other words, we have
$$
T_\ii = \bigoplus_{k=1}^r M[k_{\rm max},k] = I_w \oplus \Omega_w^{-1}(V_\ii).
$$

\subsection{Computation of $\dm \Hom_\LL(V_k,M_s)$}

\begin{Lem}\label{M1}
Let $1 \le k,s \le r$.
\begin{itemize}

\item[(i)]
If $k \le s$, then we have
$$
\dm \Hom_\LL(V_k,M_s) = \dm \Hom_\LL(M_k,M_s) 
\cong
\begin{cases}
0 & \text{if $k < s$},\\
1 & \text{if $k = s$}.
\end{cases}
$$

\item[(ii)]
If $k > s$, then
$$
\dm \Hom_\LL(V_k,M_s) = 
\begin{cases}
\sum_{m \ge 0,k^{(-m)} > s}\; 
(M_{k^{(-m)}},M_s)_Q & \text{if $i_k \not= i_s$},\\
1 + \sum_{m \ge 0,k^{(-m)} > s}\; 
(M_{k^{(-m)}},M_s)_Q & \text{if $i_k = i_s$}.
\end{cases}
$$

\item[(iii)]
We have
$$
\dm \Hom_\LL(V_k,V_s) = \dm \Hom_\LL(V_k,M_s \oplus M_{s^-} \oplus
\cdots \oplus M_{s_{\rm min}}). 
$$

\end{itemize}
\end{Lem}

\begin{proof}
We have short exact sequences
$$
\eta:\;\;\; 0 \to V_{k^-} \xrightarrow{\iota_k} V_k 
\xrightarrow{\pi_k} M_k \to 0
\text{\;\;\; and \;\;\;}
\psi:\;\;\; 0 \to V_{s^-} \xrightarrow{\iota_s} V_s 
\xrightarrow{\pi_s} M_s \to 0.
$$
First, assume that $k < s$.
Then the module $M_k$ is contained in $\CC_{(i_s,\ldots,i_1)}$
and also in $\CC_{(i_{s-1},\ldots,i_1)}$
Now $V_s$ is $\CC_{(i_s,\ldots,i_1)}$-projective-injective
and $V_{s^-}$ is $\CC_{(i_{s-1},\ldots,i_1)}$-projective-injective.
This implies 
$$
\dm \Hom_\LL(M_k,V_{s^-}) = \dm \Hom_\LL(M_k,V_s)
\text{\;\;\; and \;\;\;}
\dm \Ext_\LL^1(M_k,V_{s^-}) = 0.
$$
Now apply $\Hom_\LL(M_k,-)$ to the sequence $\psi$
and get $\Hom_\LL(M_k,M_s) = 0$.
Next, apply $\Hom_\LL(-,M_s)$ to $\eta$.
We have $\Hom_\LL(M_k,M_s) = 0$ and by induction we also
get $\Hom_\LL(V_{k^-},M_s) = 0$.
This implies $\Hom_\LL(V_k,M_s) = 0$.
 
Next, let $k = s$.
We apply $\Hom_\LL(-,M_k)$ to $\eta$.
Since $\Hom_\LL(V_{k^-},M_k) = 0$, we get
$\dm \Hom_\LL(V_k,M_k) = \dm \Hom_\LL(M_k,M_k)$.

Applying $\Hom_\LL(V_k,-)$ to $\eta$ gives an exact sequence 
$$
0 \to \Hom_\LL(V_k,V_{k^-}) \xrightarrow{\Hom_\LL(V_k,\iota_k)} 
\Hom_\LL(V_k,V_k) \xrightarrow{\Hom_\LL(V_k,\pi_k)} \Hom_\LL(V_k,M_k) \to 0. 
$$
Here we use that $V_{k^-}$ is contained in $\CC_{(i_k,\ldots,i_1)}$
and $V_k$ is $\CC_{(i_k,\ldots,i_1)}$-projective-injective.
Thus every homomorphism $h\df V_k \to M_k$ factors through $\pi_k$.
In other words, there exists some $g\df V_k \to V_k$ such that
$\pi_k \circ g = h$.
Now $V_k$ is indecomposable, so the endomorphism ring $\End_\LL(V_k)$
is local.
Therefore $g = \lambda{\rm id}_{V_k} + g'$ for some nilpotent
endomorphism $g'$ and some $\lambda \in K$.
Now we easily see that the image of $g'$ is contained in 
$\iota_k(V_{k^-})$.
Thus $h = \lambda \pi_k$.
This implies $\dm \Hom_\LL(V_k,M_k) = 1$.

Finally, assume that $k > s$.
Then Lemma~\ref{extandhom} yields
\begin{align*}
\dm \Ext_\LL^1(V_k,M_s) &= \dm \Hom_\LL(V_k,M_s) + \dm \Hom_\LL(M_s,V_k) 
- (V_k,M_s)_Q\\
&= \dm \Hom_\LL(V_k,M_s) + \dm \Hom_\LL(M_s,V_k) \\
&\;\;\;\;- (V_{k^-},M_s)_Q - (M_k,M_s)_Q\\
&= \dm \Hom_\LL(V_k,M_s) + \dm \Hom_\LL(M_s,V_k) 
+ \dm \Ext_\LL^1(V_{k^-},M_s)\\
&\;\;\;\;- \dm \Hom_\LL(V_{k^-},M_s) - 
\dm \Hom_\LL(M_s,V_{k^-}) - (M_k,M_s)_Q.
\end{align*}
Since $s < k$, we have
$\dm \Hom_\LL(M_s,V_{k^-}) = \dm \Hom_\LL(M_s,V_k)$
and $\Ext_\LL^1(V_k,M_s) = \Ext_\LL^1(V_{k^-},M_s) = 0$.
Thus we get
$$
\dm \Hom_\LL(V_k,M_s) = (M_k,M_s)_Q + \dm \Hom_\LL(V_{k^-},M_s).
$$
The result follows by induction.

To prove (iii) we just apply $\Hom_\LL(V_k,-)$ to the short exact sequence
$0 \to V_{s^-} \to V_s \to M_s \to 0$, and then use induction.
\end{proof}

Note that in general we have 
$\dm \Hom_\LL(V_k,M_s) \not= \dm \Hom_\LL(M_k,M_s)$.

\begin{Cor}\label{M2}
For $1 \le k \le r$ we have
$\Ext_\LL^1(M_k,M_k) = 0$.
\end{Cor}

\begin{proof}
Again we use the short exact sequence
$$
\eta:\;\;\; 0 \to V_{k^-} \to V_k \to M_k \to 0.
$$
The three modules in this sequence are contained in 
$\CC_{(i_k,\ldots,i_1)}$.
In particular, $V_k$ is
$\CC_{(i_k,\ldots,i_1)}$-projective-injective.
This implies $\Ext_\LL^1(V_k,M_k) = 0$.
We have $\Hom_\LL(V_{k^-},M_k) = 0$ by Lemma~\ref{M1}.
Thus, applying the functor $\Hom_\LL(-,M_k)$ to $\eta$
we get $\Ext_\LL^1(M_k,M_k) = 0$.
\end{proof}

\begin{Cor}\label{M2'}
For $1 \le k \le r$ with $k^- \not= 0$ we have
$\dm \Ext_\LL^1(M_k,V_{k^-}) = 1$.
\end{Cor}

\begin{proof}
Apply $\Hom_\LL(M_k,-)$ to the sequence $\eta$
appearing in the proof of Corollary~\ref{M2}.
\end{proof}


\section{The $\add(M_\ii)$-stratification of $\CC_w$}\label{Mstrata}


\subsection{The stratification}
Let $\aaa = (a_1,\ldots,a_r)$ be a tuple of nonnegative
integers, and let
$\CC_{M_\ii,\aaa}$ be the category of all 
$\LL$-modules $X$ such that there exists a
chain 
$$
0 = X_0 \subseteq X_1 \subseteq \cdots \subseteq X_r = X
$$
of submodules of $X$ with $X_k/X_{k-1} \cong M_k^{a_k}$ for
all $1 \le k \le r$.

\begin{Lem}\label{M3}
If $X$ is a module in $\CC_{M_\ii,\aaa}$ and $\CC_{M_\ii,\bb}$,
then $\aaa = \bb$.
\end{Lem}

\begin{proof}
Let $\aaa = (a_1,\ldots,a_r)$ and $\bb = (b_1,\ldots,b_r)$.
There is a short exact sequence
$$
0 \to X_{r-1} \to X \to M_r^{a_r} \to 0.
$$
Lemma~\ref{M1} and induction shows that $\Hom_\LL(X_{r-1},M_r) = 0$.
Thus $\dm \Hom_\LL(X,M_r) = a_r$.
Similarly, we get $\dm \Hom_\LL(X,M_r) = b_r$.
Thus $a_r = b_r$, and by induction we get $a_k = b_k$ for all $1 \le k \le r$.
\end{proof}

Define 
$$
\CC_{M_\ii} := \bigcup_{\aaa \in \N^r} \CC_{M_\ii,\aaa}.
$$

\begin{Lem}\label{M4}
We have $\CC_w = \CC_{M_\ii}$.
\end{Lem}

\begin{proof}
The category $\CC_w$ contains all $M_k$, and $\CC_w$ is closed
under extensions. 
This implies $\CC_{M_\ii} \subseteq \CC_w$.

Vice versa, assume $X \in \CC_w$.
By Proposition~\ref{cormutation1} 
there exists a short exact sequence
$$
\varepsilon:\;\;\; 
0 \to V'' \xrightarrow{f} V' \xrightarrow{g} X \to 0
$$
with $V',V'' \in \add(V_\ii)$ and
$g$ is a minimal right $\add(V_\ii)$-approximation.
We call $\varepsilon$ a {\it minimal $\add(V_\ii)$-resolution} of length
at most one.
Since $V_r$ is $\CC_w$-projective-injective, by the minimality of
$g$ we know that $V''$ does not contain a direct summand
isomorphic to $V_r$.
Let $U$ be the unique submodule of $V'$ such that
$V'/U \cong M_r^{a_r}$ with $a_r$ maximal. 
Clearly, we have
$$
U \cong V_{r^-}^{a_r} \oplus V'/V_r^{a_r}.
$$
By Lemma~\ref{M1} and induction, the image of $f$ is contained in $U$.
We have $V'/\Ima(f) \cong X$. 
Let $X_{r-1} := g(U)$.
We get $X/X_{r-1} \cong M_r^{a_r}$, and by passing to the restriction
maps, we obtain a short exact sequence
$$
0 \to V'' \xrightarrow{f'} V_{r^-}^{a_r} \oplus V'/V_r^{a_r} \to 
X_{r-1} \to 0.
$$
This is an $\add(V_\ii)$-resolution of $X_{r-1}$.
By possibly deleting a direct summand of $f'$ of the form
${\rm id}\df V_{r^-}^a \to V_{r^-}^a$, this yields again
a minimal $\add(V_\ii)$-resolution of length at most one of $X_{r-1}$.
The result follows by induction.
\end{proof}

For $X \in \CC_{M_\ii,\aaa}$ set
$$
M_\ii(X) := M_1^{a_1} \oplus \cdots \oplus M_r^{a_r}.
$$
Recall that $B_\ii := \End_\LL(V_\ii)^\op$.

For a $\LL$-module $X \in \CC_w$ 
we want to compute the dimension vector of the $B_\ii$-module
$\Hom_\LL(V_\ii,X)$.
The indecomposable projective $B_\ii$-modules are
the modules $\Hom_\LL(V_\ii,V_k)$, $1 \le k \le r$.
Thus the entries of the dimension vector
$\dimv_{B_\ii}(\Hom_\LL(V_\ii,X))$
are 
$$
\dm \Hom_{B_\ii}(\Hom_\LL(V_\ii,V_k),\Hom_\LL(V_\ii,X))
$$
where $1 \le k \le r$.
By Corollaries~\ref{fullyfaithful} and \ref{fullyfaithful3} we have
$$
\Hom_{B_\ii}(\Hom_\LL(V_\ii,V_k),\Hom_\LL(V_\ii,X)) \cong
\Hom_\LL(V_k,X).
$$
For $1 \le k \le r$ define
$$
\Delta_k := \Hom_\LL(V_\ii,M_k).
$$
(In Section~\ref{qhsection} we prove that $B_\ii$ is a quasi-hereditary
algebra and that the $\Delta_k$ are the corresponding
standard modules.)
The following result follows directly from Lemma~\ref{M1}.

\begin{Lem}\label{M4'}
The dimension vectors $\dimv_{B_\ii}(\Delta_k)$, $1 \le k \le r$
are linearly independent.
\end{Lem}

\begin{Lem}\label{M4''}
For all $1 \le k \le r$ we have 
$$
\dimv_{B_\ii}(\Hom_\LL(V_\ii,V_k)) = 
\dimv_{B_\ii}(\Delta_k) + 
\dimv_{B_\ii}(\Delta_{k^-}) + \cdots + 
\dimv_{B_\ii}(\Delta_{k_{\rm min}}).
$$
\end{Lem}

\begin{proof}
Use the short exact sequence 
$$
0 \to V_{k^-} \to V_k \to M_k \to 0
$$
and an induction on $k$.
\end{proof}

The next result shows that Lemma~\ref{M4''} is just a special case of
a general fact.

\begin{Prop}\label{M6}
For a $\LL$-module $X \in \CC_w$ and
$\aaa = (a_1,\ldots,a_r)$ the following are equivalent:
\begin{itemize}

\item[(i)]
$X \in \CC_{M_\ii,\aaa}$;

\item[(ii)]
There exists a short exact sequence
$$
0 \to \bigoplus_{k=1}^r V_{k^-}^{a_k} \to \bigoplus_{k=1}^r V_k^{a_k} \to X
\to 0;
$$

\item[(iii)]
$\dimv_{B_\ii}(\Hom_\LL(V_\ii,X)) = 
\dimv_{B_\ii}(\Hom_\LL(V_\ii,M_\ii(X))) = 
\sum_{k=1}^r a_k\; \dimv_{B_\ii}(\Delta_k)$.

\end{itemize}
\end{Prop}

\begin{proof}
${\rm (i)} \implies {\rm (ii)}$:
Assume $X \in \CC_{M_\ii,\aaa}$ with $\aaa = (a_1,\ldots,a_r)$.
By induction we get the following diagram of morphisms with exact
row and columns.
$$
\xymatrix{
& 0 \ar[d] && 0 \ar[d]\\
& \bigoplus_{k=1}^{r-1} V_{k^-}^{a_k} \ar[d] && 
V_{r^-}^{a_r} \ar[d]\\
& \bigoplus_{k=1}^{r-1} V_k^{a_k} \ar[d]^f && 
V_r^{a_r} \ar[d]^g\\
0 \ar[r] & X_{r-1} \ar[d]\ar[r]^>>>>>>{\iota} & 
X \ar[r]^>>>>>\pi & M_r^{a_r} \ar[d]\ar[r] & 0 \\
& 0 && 0
}
$$
Since $V_r$ is $\CC_w$-projective-injective, there exists a homomorphism
$g'$ such that $\pi \circ g' = g$.
Then $[f,g']\df \bigoplus_{k=1}^r V_k^{a_r} \to X$ is an epimorphism.
Let $Z := \Ker([f,g'])$.
The Snake Lemma yields an exact sequence
$$
\bigoplus_{k=1}^{r-1} V_{k^-}^{a_k} \xrightarrow{h'} Z 
\xrightarrow{h''} V_{r^-}^{a_r}.
$$
Clearly, $h''$ is an epimorphism, since $f$ is an epimorphism.
For dimension reasons $h'$ is a monomorphism.
Thus we get a short exact sequence
$$
0 \to \bigoplus_{k=1}^{r-1} V_{k^-}^{a_k} \xrightarrow{h'} Z 
\xrightarrow{h''} V_{r^-}^{a_r} \to 0.
$$
Applying $\Hom_\LL(V_\ii,-)$ 
to this sequence yields an exact sequence of $B_\ii$-modules
with a projective end term.
Thus this sequence splits, and we get
$Z = \bigoplus_{k=1}^r V_{k^-}^{a_k}$.
So we constructed a short exact sequence
$$
\eta_X\df \;\;\; 
0 \to \bigoplus_{k=1}^r V_{k^-}^{a_k} \to  
\bigoplus_{k=1}^r V_k^{a_k}
\to X \to 0.
$$

${\rm (ii)} \implies {\rm (iii)}$:
Apply $\Hom_\LL(V_\ii,-)$ to the short exact sequence $\eta_X$.
Since $V_\ii$ is rigid, this yields a short exact sequence
of $B_\ii$-modules, and we get
\begin{align*}
\dimv_{B_\ii}(\Hom_\LL(V_\ii,X)) &= 
\dimv_{B_\ii}(\Hom_\LL(V_\ii,\bigoplus_{k=1}^r V_k^{a_k})) - 
\dimv_{B_\ii}(\Hom_\LL(V_\ii,\bigoplus_{k=1}^r V_{k^-}^{a_k})) \\
&= \sum_{k=1}^r a_k\, 
(\dimv_{B_\ii}(\Hom_\LL(V_\ii,V_k)) - 
\dimv_{B_\ii}(\Hom_\LL(V_\ii,V_{k^-}))) \\
&= \sum_{k=1}^r a_k\, \dimv_{B_\ii}(\Delta_k).
\end{align*}
This implies (iii).

${\rm (iii)} \implies {\rm (i)}$:
Let $X \in \CC_w$, and assume $\dimv_{B_\ii}(\Hom_\LL(V_\ii,X)) = 
\sum_{k=1}^r a_k\, \dimv_{B_\ii}(\Delta_k)$.
Set $\aaa = (a_1,\ldots,a_r)$.
We know that $X \in \CC_{M_\ii,\bb}$ for some
$\bb = (b_1,\ldots,b_r)$.
By the implication ${\rm (i)} \implies {\rm (iii)}$ we get 
$\dimv_{B_\ii}(\Hom_\LL(V_\ii,X)) = 
\sum_{k=1}^r b_k\, \dimv_{B_\ii}(\Delta_k)$.
Since the vectors 
$\dimv_{B_\ii}(\Delta_1),\ldots,
\dimv_{B_\ii}(\Delta_r)$
are linearly independent, we get
$a_k = b_k$ for all $k$.
\end{proof}

\begin{Cor}\label{M7}
For $X,Y \in \CC_w$ we have
$\dimv_{B_\ii}(\Hom_\LL(V_\ii,X)) = \dimv_{B_\ii}(\Hom_\LL(V_\ii,Y))$ 
if and only
if $X,Y \in \CC_{M_\ii,\aaa}$ for some $\aaa$.
\end{Cor}

\begin{proof}
By Lemma~\ref{M4'}
the dimension vectors $\dimv_{B_\ii}(\Delta_k)$ 
are linearly independent.
Now use Proposition~\ref{M6}.
\end{proof}

A short exact sequence 
$
\eta\df 0 \to X \to Y \to Z \to 0
$
of $\LL$-modules is called $M_\ii$-{\it split} 
if
$M_\ii(X) \oplus M_\ii(Z) \cong M_\ii(Y)$.
Recall that $F_{V_\ii} := \Hom_\LL(V_\ii,-)$.

\begin{Cor}\label{QTexact}
For a short exact sequence 
$\eta\df 0 \to X \to Y \to Z \to 0$ of $\LL$-modules in $\CC_w$
the following are equivalent:
\begin{itemize}

\item[(i)]
$\eta$ is $F_{V_\ii}$-exact;

\item[(ii)]
$\eta$ is $M_\ii$-split.

\end{itemize}
\end{Cor}

\begin{proof}
Clearly, $\eta$ is $F_{V_\ii}$-exact if and only if 
$$
\dimv_{B_\ii}(\Hom_\LL(V_\ii,X)) + \dimv_{B_\ii}(\Hom_\LL(V_\ii,Z)) =
\dimv_{B_\ii}(\Hom_\LL(V_\ii,Y)).
$$
By Proposition~\ref{M6} this happens if and only if
$M_\ii(X) \oplus M_\ii(Z) \cong M_\ii(Y)$.
\end{proof}

\subsection{Example}
Let $Q$ be a quiver with underlying graph
$$
\xymatrix@-0.7pc{
1 \ar@{-}[r] & 2\ar@{-}[r] & 3
}
$$
and let $w_0$ be the longest Weyl group element in $W_Q$.
Thus we have $\CC_{w_0} = \md(\LL)$. 
The short exact sequences
$$
\eta'\df
0 \to \bsm2\esm \to \bsm1&\\&2\esm \oplus \bsm&3\\2&\esm \to 
\bsm1&&3\\&2&\esm \to 0
\text{\;\;\; and \;\;\;}
\eta''\df
0 \to \bsm1&&3\\&2&\esm \to \bsm&2&\\1&&3\\&2&\esm \to 
\bsm2\esm \to 0
$$
are exchange sequences in $\md(\LL)$.
Let $\ii = (1,2,1,3,2,1)$ and $\jj = (2,1,2,3,2,1)$ be reduced
expressions of $w_0$.
We get 
$$
M_\ii = \bsm1\esm \oplus 
\bsm1\\&2\esm \oplus 
\bsm1\\&2\\&&3\esm \oplus 
\bsm2\esm \oplus 
\bsm2\\&3\esm \oplus 
\bsm3\esm
\text{\;\;\; and \;\;\;}
M_\jj = 
\bsm1\esm \oplus 
\bsm1\\&2\esm \oplus 
\bsm1\\&2\\&&3\esm \oplus 
\bsm3\esm \oplus 
\bsm&3\\2\esm \oplus 
\bsm2\esm.
$$
Now one easily observes that $\eta'$ is $M_\ii$-split and not $M_\jj$-split,
and $\eta''$ is $M_\jj$-split but not $M_\ii$-split.


\section{Quasi-hereditary algebras associated to reduced expressions}
\label{qhsection}


\subsection{Quasi-hereditary algebras}\label{qh4}
Let $A$ be a finite-dimensional algebra.
By $P_1,\ldots,P_r$ and $Q_1,\ldots,Q_r$ and
$S_1,\ldots,S_r$ we denote the indecomposable 
projective, indecomposable injective and simple $A$-modules,
respectively,
where $S_i = {\rm top}(P_i)  = \soc(Q_i)$. 

For a class $\U$ of $A$-modules let $\F(\U)$ 
be the class of all $A$-modules 
$X$ which
have a filtration 
$$
0 = X_0 \subseteq X_1 \subseteq \cdots \subseteq X_t = X
$$
of submodules 
such that all factors $X_j/X_{j-1}$ belong to $\U$ for all $1 \le j \le t$.
Such a filtration is called a $\U$-{\it filtration} of $X$.
We call these modules the $\U$-{\it filtered modules}.

Fix a bijective map $\omega\df \{ S_1,\ldots,S_r \} \to \{ 1,\ldots,r \}$.
Let $\Delta_i$ be the largest factor module of $P_i$ such that
$[\Delta_i:S_j] = 0$ for all $j$ with $\omega(S_j) > \omega(S_i)$,
and set 
$$
\Delta = \{ \Delta_1,\ldots,\Delta_r \}.
$$
The modules $\Delta_i$ are called {\it standard modules}.
The algebra $A$ is called {\it quasi-hereditary} 
if 
$\End_A(\Delta_i) \cong K$ for all $i$, and if
${_A}A$ belongs to $\F(\Delta)$.
Quasi-hereditary algebras first occured in Cline, Parshall
and Scott's \cite{CPS} study of highest weight categories. 

Note that the definition of a quasi-hereditary algebra 
depends on the chosen ordering of the simple modules.
If we reorder them, it could happen that our algebra 
is no longer quasi-hereditary.

Now assume $A$ is a quasi-hereditary algebra, and let $\F(\Delta)$ be the
subcategory of $\Delta$-filtered $A$-modules.
For $X \in \F(\Delta)$ let $[X:\Delta_i]$ be the number of times
that $\Delta_i$ occurs as a factor in a $\Delta$-filtration
of $X$.
Then 
$$
\dimv_\Delta(X) = ([X:\Delta_1],\ldots,[X:\Delta_r]) \in \N^r
$$ 
is the $\Delta$-{\it dimension vector}
of $X$.
Let $\nabla_i$ be the largest submodule of $Q_i$ such that
$[\nabla_i:S_j] = 0$ for all $j$ with $\omega(S_j) > \omega(S_i)$,
and let 
$$
\nabla = \{ \nabla_1,\ldots,\nabla_r \}.
$$
The modules $\nabla_i$ are called {\it costandard modules}.
The following results (and the missing definitions) 
can be found in \cite{Ri5,Ri4}:
\begin{itemize}

\item[(i)]
There is a unique (up to isomorphism) basic tilting module
$T(\Delta\cap\nabla)$ 
over $A$ such that 
$$
\add(T(\Delta\cap\nabla)) = \F(\Delta) \cap \F(\nabla).
$$

\item[(ii)]
$\F(\Delta)$ is closed under extensions and
under direct summands.

\item[(iii)]
$[P_i:\Delta_j] = [\nabla_j:S_i]$ for all $1 \le i,j \le r$.

\item[(iv)]
If $X \in \F(\Delta)$, then
$[X:\Delta_i] = \dm \Hom_A(X,\nabla_i)$ for all $i$.

\item[(v)]
$\Hom_A(\Delta_i,\Delta_j) = 0$ for all $i,j$ with $\omega(S_i) > \omega(S_j)$.

\item[(vi)]
$\Ext_A^1(\Delta_i,\Delta_j) = 0$ for all $i,j$ with 
$\omega(S_i) \ge \omega(S_j)$.

\item[(vii)]
The $\F(\Delta)$-projective modules are the projective $A$-modules.
The $\F(\nabla)$-injective modules are the injective $A$-modules.

\item[(viii)]
The $\F(\Delta)$-injective modules are the modules
in $\add(T(\Delta\cap\nabla))$.
The $\F(\nabla)$-projective modules are the modules
in $\add(T(\Delta\cap\nabla))$.

\item[(ix)]
If $\Ext_A^1(X,\nabla_i) = 0$ for all $i$, then
$X \in \F(\Delta)$.
Similarly, if $\Ext_A^1(\Delta_i,Y) = 0$ for all $i$, then
$Y \in \F(\nabla)$.

\end{itemize}
The module $T(\Delta\cap\nabla)$ 
is called the {\it characteristic tilting module} of $A$.
In general, $T(\Delta\cap\nabla)$ is not a classical tilting module.
(Here a tilting module is called {\it classical} provided its
projective dimension is at most one.)
The endomorphism algebra $\End_A(T(\Delta\cap\nabla))$ is
called the {\it Ringel dual} of $A$.
It is again a quasi-hereditary algebra in a natural way, see \cite{Ri5}.

Following Ringel \cite{Ri7}, the finite-dimensional algebra $A$ is 
{\it strongly quasi-hereditary} if there is a bijective map 
$\omega\df \{ S_1,\ldots,S_r \} \to \{ 1,\ldots,r \}$
such that for each $1 \le k \le r$ there is a short exact
sequence
$$
0 \to R_k \to P_k \to D_k \to 0 
$$
satisfying the following two properties:
\begin{itemize}

\item[(1)]
$R_k$ is a direct sum of indecomposable projective $A$-modules
$P_j$ with $\omega(j) > \omega(k)$; 

\item[(2)]
$[D_k:S_j] = 
\begin{cases}
0 & \text{if $\omega(j) > \omega(k)$},\\
1 & \text{if $j = k$}.
\end{cases}
$
\end{itemize}
Each strongly quasi-hereditary algebra is quasi-hereditary with
$\Delta_k = D_k$ for all $k$.
Furthermore, we have $\pdim(\Delta_k) \le 1$ for all $k$.
If each of the modules $R_k$ is indecomposable, then one
easily checks that $A$ is $\Delta$-serial, \textit{i.e.}~each $P_k$
has a unique $\Delta$-filtration.

\subsection{The algebra $B_\ii$ is quasi-hereditary}\label{qh2}
As before, let $V_\ii = V_1 \oplus \cdots \oplus V_r$ and
$M_\ii = M_1 \oplus \cdots \oplus M_r$.
Set $B_\ii := \End_\LL(V_\ii)^\op$.
For $1 \le k \le r$ let $S(k)$ be the (simple) top of
the indecomposable $B_\ii$-modules $P_k := \Hom_\LL(V_\ii,V_k)$.
As before, define $\Delta_k := \Hom_\LL(V_\ii,M_k)$, and set
$$
\Delta := \left\{ \Delta_1,\ldots,\Delta_r \right\}.
$$
Define
$\omega\df \{ S(1),\ldots,S(r) \} \to \{ 1,\ldots,n \}$
by $\omega(S(k)) := r-k+1$.

The following theorem was first proved in \cite[Section 16]{GLSUni1}
for adaptable Weyl group
elements. 
Later the statement was generalized to arbitrary Weyl group 
elements by Iyama and Reiten~\cite{IR}.
Here we present a proof for the general case, which is very similar
to our original proof of the adaptable case.

\begin{Thm}\label{qh1}
Let $\ii$ be a reduced expression of a Weyl group element $w$.
The following hold:
\begin{itemize}

\item[(i)]
The algebra $B_\ii= \End_\LL(V_\ii)^\op$ is strongly quasi-hereditary 
and $\Delta$-serial with standard modules
$\Delta = \left\{ \Delta_1,\ldots,\Delta_r \right\}$;

\item[(ii)]
The functor $\Hom_\LL(V_\ii,-)$ yields an equivalence of categories
$F_\ii\df \CC_w \to \F(\Delta)$;

\item[(iii)]
$T(\Delta\cap\nabla) = \Hom_\LL(V_\ii,T_\ii)$.

\end{itemize}
\end{Thm}

\begin{proof}
(i):
We know that for each $1 \le k \le r$ there is a short exact
sequence
$$
\eta:\;\;\; 0 \to V_{k^-} \xrightarrow{\iota_k} V_k \to M_k \to 0.
$$
We apply the functor $\Hom_\LL(V_\ii,-)$ to this sequence and obtain
a short exact sequence
$$
0 \to \Hom_\LL(V_\ii,V_{k^-}) \to P_k \xrightarrow{H} \Delta_k \to 0
$$
of $B_\ii$-modules.
Let $\omega(S(j)) \ge \omega(S(k))$, 
and let $F\df \Hom_\LL(V_\ii,V_j) \to \Delta_k$
be a homomorphism of $B_\ii$-modules.
Since $\Hom_\LL(V_\ii,V_j)$ is a projective $B_\ii$-module, there
is a homomorphism $G\df \Hom_\LL(V_\ii,V_j) \to P_k$ such that
$H \circ G = F$.
There exists a $\LL$-module homomorphism $g\df V_j \to V_k$ such that
$G = \Hom_\LL(V_\ii,g)$.
Assume $\omega(S(j)) > \omega(S(k))$.
Since $j < k$, we know that $\Ima(g) \subseteq \iota_k(V_{k^-})$.
Thus $\Ima(G) \subseteq \Ima(\Hom_\LL(V_\ii,\iota_k)) = \Ker(H)$.
But this implies $F = 0$.
Therefore we have $[\Delta_k:S(j)] = 0$.
Next, we consider the case $\omega(S(j)) = \omega(S(k))$.
The endomorphism ring $\End_\LL(V_k)$ is local, and we work over
an algebraically closed field.
Thus $g = \lambda{\rm id}_{V_k} + g'$ with $g'$ nilpotent and $\lambda \in K$.
We have $\soc(V_k) \subseteq \Ker(g')$.
This implies $\Ima(g') \subseteq \iota_k(V_{k^-})$.
Thus $F = H \circ G = H \circ \Hom_\LL(V_\ii,\lambda{\rm id}_{V_k})$.
In other words, $\Hom_{B_\ii}(P_k,\Delta_k)$ is 1-dimensional.
This finishes the proof of (i).

(ii):
For $X,Z \in \CC_w$ we have a functorial isomorphism
$$
\Ext_{F_{V_\ii}}^1(X,Z) \to 
\Ext_{B_\ii}^1(\Hom_\LL(V_\ii,X),\Hom_\LL(V_\ii,Z)).
$$
Thus the image of the functor
$$
\Hom_\LL(V_\ii,-)\df \CC_w \to \md(B_\ii)
$$
is extension closed.
Clearly, for all $1 \le k \le r$ the standard module 
$\Delta_k$ is in $\Hom_\LL(V_\ii,\CC_w)$.
It follows that $\F(\Delta) \subseteq \Hom_\LL(V_\ii,\CC_w)$.

Now let $X \in \CC_w$.
By Lemma~\ref{M4} we know that 
$X \in \CC_{M_\ii,\aaa}$ for some $\aaa = (a_1,\ldots,a_r)$.
Thus there is a short exact sequence
$$
\eta:\;\;\; 0 \to X_{r-1} \to X \to M_r^{a_r} \to 0.
$$
We claim that $\eta$ is $F_{V_\ii}$-exact:
Clearly, $\eta$ is $F_{V_r}$-exact, since $V_r$ is
$\CC_w$-projective-injective and $X_{r-1} \in \CC_w$.
Since $\Hom_\LL(V_k,M_r) = 0$ for all $k < r$, it follows
that $\eta$ is also $F_{V_k}$-exact for all such $k$.
Clearly, $\Hom_\LL(V_\ii,M_r^{a_r})$ is contained in $\F(\Delta)$.
By induction also $\Hom_\LL(V_\ii,X_{r-1})$ is in $\F(\Delta)$.
Since $\F(\Delta)$ is closed under extensions, and since
$\eta$ is $F_{V_\ii}$-exact, we get that $\Hom_\LL(V_\ii,X)$ is
in $\F(\Delta)$.
So we proved that $\F(\Delta) = \Hom_\LL(V_\ii,\CC_w)$.
Now Corollary~\ref{fullyfaithful} and Lemma~\ref{fullyfaithful3} 
show that the restriction functor $F_\ii\df \CC_w \to \F(\Delta)$ 
is an equivalence of categories.

(iii):
It is enough to show that $\Ext_\LL^1(\Delta_k,T_\ii) = 0$
for all $1 \le k \le r$, see Section~\ref{qh4}.
Recall that all indecomposable direct summands of $T_\ii$ are of
the form $M[s_{\rm max},s]$ where $1 \le s \le r$.
We fix such an $s$.

For each $1 \le k \le r$ there is a short exact sequence
$$
\eta:\;\;\; 0 \to M[k^-,k_{\rm min}] \to M[k,k_{\rm min}] \to M_k \to 0.
$$
Applying $\Hom_\LL(V_\ii,-)$ yields a projective resolution
$$
0 \to \Hom_\LL(V_\ii,M[k^-,k_{\rm min}]) \to 
\Hom_\LL(V_\ii,M[k,k_{\rm min}]) \to \Hom_\LL(V_\ii,M_k) \to 0
$$
of $B_\ii$-modules.

If $k \le s$, then 
$\Hom_\LL(M[k^-,k_{\rm min}],M[s_{\rm max},s]) = 0$.
Since $F_\ii$ is an equivalence, we get
$$
\Hom_{B_\ii}(\Hom_\LL(V_\ii,M[k^-,k_{\rm min}]),
\Hom_\LL(V_\ii,M[s_{\rm max},s])) = 0.
$$
This implies $\Ext_{B_\ii}^1(\Delta_k,\Hom_\LL(V_\ii,M[s_{\rm max},s])) = 0$.

Next, assume that $k > s$.
We have a short exact sequence
$$
\psi:\;\;\; 0 \to M[s^-,s_{\rm min}] \to M[s_{\rm max},s_{\rm min}] 
\to M[s_{\rm max},s] \to 0.
$$
Applying $\Hom_\LL(-,M_k)$ yields 
$\Ext_\LL^1(M[s_{\rm max},s],M_k) = 0$.
Thus $\Ext_\LL^1(M_k,M[s_{\rm max},s]) = 0$.
This implies 
$$
\Ext_{F_{V_\ii}}^1(M_k,M[s_{\rm max},s]) = 
\Ext_{B_\ii}^1(\Delta_k,\Hom_\LL(V_\ii,M[s_{\rm max},s])) = 0.
$$
Here we used that 
$\Ext_\LL^1(M[s_{\rm max},s_{\rm min}],M_k) = 0$ (since
$M[s_{\rm max},s_{\rm min}]$ is $\CC_w$-projective-injective),
and $\Hom_\LL(M[s^-,s_{\rm min}],M_k) = 0$ by Lemma~\ref{M1}.
This finishes the proof of (iii).
\end{proof}

\begin{Cor}\label{M5}
The modules $\Hom_\LL(V_\ii,I_{\ii,j})$, $1 \le j \le n$ are
the indecomposable $\F(\Delta)$-projective-injectives modules.
\end{Cor}

\begin{proof}
This follows from Theorem~\ref{qh1}, (iii) and Section~\ref{qh4}.
\end{proof}

Each $\Delta$-filtration of the indecomposable projective
$B_\ii$-module $\Hom_\LL(V_\ii,V_k)$ looks as follows:
$$
\begin{tabular}{c}
$\Delta_k$\\
\hline
$\Delta_{k^-}$\\
\hline
$\cdots$\\
\hline
$\Delta_{k^{\min}}$
\end{tabular}
$$
(We just displayed the factors of the (unique) $\Delta$-filtration of
$\Hom_\LL(V_\ii,V_k)$.)

We can now reformulate parts of Proposition~\ref{M6} as follows:

\begin{Prop}\label{M6'}
For a $\LL$-module $X \in \CC_w$ and
$\aaa = (a_1,\ldots,a_r)$ the following are equivalent:
\begin{itemize}

\item[(1)]
$X \in \CC_{M_\ii,\aaa}$;

\item[(2)]
$\dimv_\Delta(F_\ii(X)) = (a_1,\ldots,a_r)$.

\end{itemize}
\end{Prop}

\begin{proof}
Since $\Delta_k = \Hom_\LL(V_\ii,M_k)$, it is clear that (iii) in
Proposition~\ref{M6} and (2) are equivalent.
\end{proof}

We know that $B_\ii$ is an algebra of finite global dimension.
Thus one can define the Ringel form
$$
\bil{X,Y}_{B_\ii} := \bil{\dimv(X),\dimv(Y)}_{B_\ii} := 
\sum_{j \ge 0} (-1)^j \dm \Ext_{B_\ii}^j(X,Y).
$$
The next lemma gives the values of $\bil{-,-}_{B_\ii}$
applied to standard modules.

\begin{Lem}\label{M8}
For $1 \le k,s \le r$ we have
$$
\bil{\Delta_k,\Delta_s}_{B_\ii} = \dm \Hom_{B_\ii}(\Delta_k,\Delta_s) -
\dm \Ext_{B_\ii}^1(\Delta_k,\Delta_s) =
\begin{cases}
0 & \text{if $k < s$},\\
1 & \text{if $k = s$},\\
(M_k,M_s)_Q & \text{if $k > s$}.
\end{cases}
$$
\end{Lem}

\begin{proof}
As before, for $1 \le t \le r$ we set $P_t := \Hom_\LL(V_\ii,V_t)$
and $\Delta_t := \Hom_\LL(V_\ii,M_t)$.
We know that $\pdim(\Delta_t) \le 1$ for all $t$.
Thus
$$
\bil{\Delta_k,\Delta_s}_{B_\ii} = \dm \Hom_{B_\ii}(\Delta_k,\Delta_s) -
\dm \Ext_{B_\ii}^1(\Delta_k,\Delta_s).
$$

The cases $k < s$ and $k = s$ are clear, see Section~\ref{qh4}.
Thus, assume $k > s$.
The short exact sequence
$$
0 \to V_{k^-} \to V_k \to M_k \to 0
$$
yields a projective resolution
$$
0 \to P_{k^-} \to P_k \to \Delta_k \to 0
$$
of $\Delta_k$.
We apply $\Hom_\LL(-,M_s)$ and obtain an exact sequence
$$
0 \to \Hom_{B_\ii}(\Delta_k,\Delta_s) \to \Hom_{B_\ii}(P_k,\Delta_s) \to
\Hom_{B_\ii}(P_{k^-},\Delta_s) \to \Ext_{B_\ii}^1(\Delta_k,\Delta_s) \to 0.
$$
This implies
\begin{align*}
\bil{\Delta_k,\Delta_s}_{B_\ii} &= 
\dm \Hom_{B_\ii}(P_k,\Delta_s) - \dm \Hom_{B_\ii}(P_{k^-},\Delta_s) \\
&= \dm \Hom_\LL(V_k,M_s) - \dm \Hom_\LL(V_{k^-},M_s)\\
&= (M_k,M_s)_Q.
\end{align*}
For the third equality we use
Lemma~\ref{M1}.
\end{proof}

\subsection{Example}
For an arbitrary $\CC_w$-maximal rigid $\LL$-module $T$, it seems
to be difficult to determine when $\End_\LL(T)^\op$ is
quasi-hereditary and when not.

Even if $Q$ is a quiver with underlying graph
$$
\xymatrix@-0.7pc{
1 \ar@{-}[r] & 2\ar@{-}[r] & 3
}
$$
there are maximal rigid modules whose endomorphism
algebra is not quasi-hereditary:
Let $w = w_0$ be the longest Weyl group element in $W_Q$.
Let $T$
be the $\CC_w$-maximal rigid $\LL$-module
$$
{\bsm2&\\&3\esm} \oplus {\bsm&2&\\1&&3\esm} \oplus 
{\bsm&2\\1&\esm} \oplus
{\bsm1&&\\&2&\\&&3\esm} \oplus {\bsm&2&\\1&&3\\&2&\esm} \oplus 
{\bsm&&3\\&2&\\1&&\esm}.
$$
The quiver of $\End_\LL(T)^\op$ looks as follows:
$$
\xymatrix@-0.7pc{
{\bsm2&\\&3\esm} \ar[r] & {\bsm&2&\\1&&3\esm} \ar[d] & 
{\bsm&2\\1&\esm} \ar[l]\\
{\bsm1&&\\&2&\\&&3\esm} \ar[u] & {\bsm&2&\\1&&3\\&2&\esm} \ar[l]\ar[r] & 
{\bsm&&3\\&2&\\1&&\esm} \ar[u]
}
$$
It is not difficult to show that $\End_\LL(T)^\op$ is not a quasi-hereditary
algebra.


\section{Mutations of clusters via dimension vectors}


\subsection{Dimension vectors of rigid modules}
Let $A$ be a finite-dimensional $K$-algebra.
For $m \ge 0$ let $A^m$ be the free $A$-module of rank $m$.
By $\md(A,m)$ we denote the affine variety of $m$-dimensional $A$-modules.
(One can define $\md(A,m)$ as the variety of $K$-algebra
homomorphisms $A \to M_m(K)$.)
If $U$ is a submodule of $A^m$ such that $A^m/U$ is $m$-dimensional,
then the {\it Richmond stratum} $\cS(U,A^m)$ is the 
subset of $\md(A,m)$ consisting of the modules $X$ such that
there exists a short exact sequence
$$
0 \to U \to A^m \to X \to 0,
$$
see \cite{R}.
A more general situation was studied by Bongartz \cite{Bo}.

\begin{Thm}[{\cite[Theorem 1]{R}}]
The Richmond stratum
$\cS(U,A^m)$ is a smooth, irreducible, locally closed subset 
of $\md(A,m)$, and
$$
\dm \cS(U,A^m) = \dm \Hom_A(U,A^m) - \dm \End_A(U).
$$
\end{Thm}

\begin{Prop}\label{rigiddim}
Assume that $\gldim(A) < \infty$.
Let $M$ and $N$ be rigid $A$-modules of projective dimension at most one.
If $\dimv(M) = \dimv(N)$, then $M \cong N$.
\end{Prop}

\begin{proof}
Let $m$ be the $K$-dimension of $M$ and $N$.
Thus, there are projective resolutions
$$
0 \to P \to A^m \to M \to 0
\text{\;\;\; and \;\;\;}
0 \to P' \to A^m \to N \to 0
$$
of $M$ and $N$, respectively.
Here we used that the projective dimensions of $M$ and $N$ are at most one.
Since $\dimv(M) = \dimv(N)$, we get $\dimv(P) = \dimv(P')$.
Since $A$ is a finite-dimensional algebra of finite global dimension,
its Cartan matrix is invertible.
In other words, the dimension vectors of the indecomposable projective
$A$-modules are linearly independent.
Thus we get $P \cong P'$.

Since $M$ and $N$ are rigid, their $\GL_m(K)$-orbits are open
in $\md(A,m)$.
In particular, these orbits are open in the Richmond stratum
$\cS(P,A^m)$.
But $\cS(P,A^m)$ is irreducible, and therefore it can contain at most 
one open orbit.
It follows that $M \cong N$.
\end{proof}

Now, let $\CC_w = \Gen(V_\ii)$ be defined as before, and let
$T = T_1 \oplus \cdots \oplus T_r$ 
be a fixed basic $\CC_w$-maximal rigid module and
set $B := \End_{\LL}(T)^\op$.

\begin{Cor}
Let $X$ and $Y$ be indecomposable rigid modules in $\CC_w$.
If 
$$
\dimv_B(\Hom_\LL(T,X)) = \dimv_B(\Hom_\LL(T,Y)), 
$$
then $X \cong Y$.
\end{Cor}

\begin{proof}
Use Corollary~\ref{pdimone} and Proposition~\ref{corquivershape},(vi),
and then apply Proposition~\ref{rigiddim}.
\end{proof}

\subsection{Mutations via dimension vectors}\label{dimvectormutation1}
We now explain how to calculate mutations of clusters via dimension vectors.
We start with some notation:
For 
${\mathbf d} = (d_1,\ldots,d_r)$ and ${\mathbf f} = (f_1,\ldots,f_r)$
in $\Z^r$ define
$$
\max\{{\mathbf d},{\mathbf f}\} := (h_1,\ldots,h_r)
$$
where $h_s = \max\{d_s,f_s\}$ for $1 \le s \le r$.
Set ${\rm Max}\{{\mathbf d},{\mathbf f}\} := {\mathbf d}$
if $d_s \ge f_s$ for all $s$.
In this case, we write
${\mathbf d} \ge {\mathbf f}$.
Of course, ${\rm Max}\{{\mathbf d},{\mathbf f}\} = {\mathbf d}$
implies
$\max\{{\mathbf d},{\mathbf f}\} = {\mathbf d}$.
By $|{\mathbf d}|$ we denote the sum of the entries of ${\mathbf d}$.

Let $\GG$ be a quiver without loops and without 2-cycles and with
vertices $1,\ldots,r$. 
Some of these vertices can be considered as {\it frozen vertices}, 
\textit{i.e.}~one cannot perform a mutation at these vertices.

Now replace each vertex $s$ of $\GG$ by some
${\mathbf d}_s \in \Z^r$.
Thus we obtain a new quiver $\GG'$ whose 
vertices are elements in $\Z^r$.

For $k$ not a frozen vertex, define
the mutation $\mu_{{\mathbf d}_k}(\GG')$ of $\GG'$ at
the vertex ${\mathbf d}_k$ in two steps:
\begin{itemize}

\item[(1)]
Replace the vertex ${\mathbf d}_k$ of $\GG'$ by
$$
{\mathbf d}_k^* := -{\mathbf d}_k + 
\max\left\{\sum_{{\mathbf d}_i \to {\mathbf d}_k} {\mathbf d}_i,
\sum_{{\mathbf d}_k \to {\mathbf d}_j}{\mathbf d}_j \right\}
$$
where the sums are taken over all arrows in $\GG'$ which start,
respectively end in the vertex ${\mathbf d}_k$;

\item[(2)]
Change the arrows of $\GG'$ following Fomin and Zelevinsky's
quiver mutation rule for the vertex ${\mathbf d}_k$.

\end{itemize}

Thus starting with $\GG'$ we can use iterated mutation
and obtain quivers whose vertices
are elements in $\Z^r$.

For example,
if for each $s$ we choose 
${\mathbf d}_s = -{\mathbf e}_s$, where ${\mathbf e}_s$ is
the $s$th canonical basis vector of $\Z^r$, then
the resulting vertices (\textit{i.e.}~elements in $\Z^r$) are the denominator
vectors of the cluster variables of the cluster algebra
$\cA(B(\GG)^\circ)$ associated to $\GG$,
compare with~\cite[Section 7, Equation (7.7)]{FZ3}.
(The variables attached to the frozen vertices serve as 
(non-invertible) coefficients.
To obtain the denominator vectors as defined in \cite{FZ3} one
has to ignore the entries corresponding to these $n$ coefficients.)
It is an open problem, if these denominator vectors actually parametrize
the cluster variables of $\cA(B(\GG)^\circ)$.

We will show that for an appropriate choice of $\GG$ and of 
the initial vectors ${\mathbf d}_s$,
the quivers obtained by iterated mutation of $\GG'$ are
in bijection with the seeds and clusters of $\cA(B(\GG)^\circ)$.
All resulting vertices (including the ${\mathbf d}_s$) will be elements
in $\N^r$, and we will show that for our particular choice of
initial vectors, we can
use ``{\rm Max}'' instead of ``{\rm max}'' in the formula above.
(This holds for all iterated mutations.)

For the rest of this section let $T = T_1 \oplus \cdots \oplus T_r$
be a basic $\CC_w$-maximal rigid $\LL$-module, and set 
$B := \End_\LL(T)^\op$.

\begin{Prop}\label{dimcount}
Let $R = R_1 \oplus \cdots \oplus T_r$ be a basic 
$\CC_w$-maximal rigid $\LL$-module.
Let 
$$
\eta'\df
0 \to R_k \xrightarrow{f'} R' \xrightarrow{g'} R_k^* \to 0 
\text{\;\;\;and \;\;\;}
\eta''\df
0 \to R_k^* \xrightarrow{f''} R'' \xrightarrow{g''} R_k \to 0
$$
be the two exchange sequences associated to an
indecomposable direct summand $R_k$ of $R$ which is not 
$\CC_w$-projective-injective.
Then 
$
\dm \Hom_\LL(T,R') \not= \dm \Hom_\LL(T,R''),
$
and we have
\begin{multline*}
\dimv_B(\Hom_\LL(T,R_k)) +  \dimv_B(\Hom_\LL(T,R_k^*)) = \\
= \max\{ \dimv_B(\Hom_\LL(T,R')), \dimv_B(\Hom_\LL(T,R'')) \}.
\end{multline*}
Furthermore, the following are equivalent:
\begin{itemize}

\item[(i)]
$\eta'$ is $F_T$-exact;

\item[(ii)] 
$\dm \Hom_\LL(T,R') > \dm \Hom_\LL(T,R'')$;

\item[(iii)]
$\dimv_B(\Hom_\LL(T,R')) \ge \dimv_B(\Hom_\LL(T,R''))$.

\end{itemize}
\end{Prop}

\begin{proof}
By Corollary~\ref{tilt} we know that $\Hom_\LL(T,R)$ is
a classical tilting module over $B$.
Thus we can apply \cite[Lemma 2.2]{H3} and assume without loss of 
generality that
$$
\Ext_B^1(\Hom_\LL(T,R_k),\Hom_\LL(T,R_k^*)) = 0.
$$
By Proposition~\ref{seslift},
\begin{align*}
1 = \dm \Ext_\LL^1(R_k^*,R_k) &\ge \dm \Ext_{F_T}^1(R_k^*,R_k)\\
&= \dm \Ext_B^1(\Hom_\LL(T,R_k^*),\Hom_\LL(T,R_k)) > 0.
\end{align*}
This implies $\Ext_\LL^1(R_k^*,R_k) = \Ext_{F_T}^1(R_k^*,R_k)$. 
Thus $\eta'$ is $F_T$-exact, and
$$
\eta\df 
0 \to \Hom_\LL(T,R_k) \xrightarrow{\Hom_\LL(T,f')} \Hom_\LL(T,R') 
\xrightarrow{\Hom_\LL(T,g')} \Hom_\LL(T,R_k^*) \to 0
$$
is a (non-split) short exact sequence.
If we apply $\Hom_\LL(T,-)$ to $\eta''$,
we obtain an exact sequence
$$
0 \to \Hom_\LL(T,R_k^*) \xrightarrow{\Hom_\LL(T,f'')} \Hom_\LL(T,R'') 
\xrightarrow{\Hom_\LL(T,g'')} \Hom_\LL(T,R_k).
$$
Now $\Hom_\LL(T,g'')$ cannot be an epimorphism, since that would 
yield a non-split extension and we know that 
$\Ext_B^1(\Hom_\LL(T,R_k),\Hom_\LL(T,R_k^*)) = 0$.
Thus for dimension reasons we get
$\dm \Hom_\LL(T,R') > \dm \Hom_\LL(T,R'')$.
Using the functors $\Hom_B(P,-)$ where $P$ runs through the
indecomposable projective $B$-modules, it also follows that
$\dimv_B(\Hom_\LL(T,R')) > \dimv_B(\Hom_\LL(T,R''))$.
Finally, the formula for dimension vectors follows from the
exactness of $\eta$.
\end{proof}

Proposition~\ref{dimcount} yields an easy combinatorial rule for the 
mutation of $\CC_w$-maximal rigid modules.
Let 
$R = R_1 \oplus \cdots \oplus R_r$ 
be a basic $\CC_w$-maximal rigid $\LL$-module.
Without loss of generality we assume that 
$R_{r-n+1},\ldots,R_r$
are $\CC_w$-projective-injective.
For $1 \le s \le r$ let
$
{\mathbf d}_s := \dimv_B(\Hom_\LL(T,R_s)).
$

As before, let $\GG_R$ be the quiver of $\End_\LL(R)^\op$.
The vertices of $\GG_R$ are labeled by the modules $R_s$.
For each $s$ we replace the vertex labeled by $R_s$ by the
dimension vector ${\mathbf d}_s$.
The resulting quiver is denoted by $\GG_R'$.

For $1 \le k \le r-n$ let
$$
0 \to R_k \to R' \to R_k^* \to 0
\text{\;\;\; and \;\;\;}
0 \to R_k^* \to R'' \to R_k \to 0
$$
be the two resulting exchange sequences.
We can now easily compute the dimension vector of the
$\End_\LL(T)^\op$-module $\Hom_\LL(T,R_k^*)$, namely 
Proposition~\ref{dimcount} yields that
$$
{\mathbf d}_k^* := \dimv_B(\Hom_\LL(T,R_k^*)) = 
\begin{cases}
- {\mathbf d}_k + 
\sum_{{\mathbf d}_i \to {\mathbf d}_k}{\mathbf d}_i & \text{if
$\sum_{{\mathbf d}_i \to {\mathbf d}_k} |{\mathbf d}_i| >
\sum_{{\mathbf d}_k \to {\mathbf d}_j} |{\mathbf d}_j|$},\\
- {\mathbf d}_k +  \sum_{{\mathbf d}_k \to {\mathbf d}_j}{\mathbf d}_j & 
\text{otherwise},
\end{cases}
$$
where the sums are taken over all arrows in $\GG_R'$ which start,
respectively end in the vertex ${\mathbf d}_k$.
More precisely, we have 
\begin{equation}\label{denform}
{\mathbf d}_k^* = -{\mathbf d}_k + 
\max\left\{\sum_{{\mathbf d}_i \to {\mathbf d}_k} {\mathbf d}_i,
\sum_{{\mathbf d}_k \to {\mathbf d}_j}{\mathbf d}_j \right\}
\end{equation}
and we know that
\begin{equation}\label{denform2}
\max\left\{\sum_{{\mathbf d}_i \to {\mathbf d}_k} {\mathbf d}_i,
\sum_{{\mathbf d}_k \to {\mathbf d}_j}{\mathbf d}_j \right\}
=
{\rm Max}\left\{\sum_{{\mathbf d}_i \to {\mathbf d}_k}{\mathbf d}_i,
\sum_{{\mathbf d}_k \to {\mathbf d}_j}{\mathbf d}_j \right\}.
\end{equation}

\begin{Rem}
Let $T = T_1 \oplus \cdots \oplus T_r$ be a basic $\CC_w$-maximal
rigid module, and
let $B^{(T)} := (\bil{S_i,S_j})_{1 \le i,j \le r}$
be the matrix of the Ringel form of the algebra $B := \End_\LL(T)^\op$.
Let $X$ be a $T$-reachable $\LL$-module, see
Section~\ref{reach}.
Set ${\mathbf d} := \dimv_B(\Hom_\LL(T,X)) \in \N^r$.
Define 
$$
\widetilde{g}_T(X) := {\mathbf d} \cdot B^{(T)},
$$
where ${\mathbf d}$ is considered as a row vector.
As explained in \cite[Section 4]{FK} the entries
of $\widetilde{g}_T(X)$, which correspond to the 
non-$\CC_w$-projective-injective
direct summands $T_k$ of $T$ 
form precisely the $g$-vector
of $\varphi_X$ with respect to the initial cluster
$(\delta_{T_1},\ldots,\delta_{T_r})$.
\end{Rem}

\subsection{Examples (Dimension vectors of $B_\ii$-modules)}
Let $Q$ be a quiver with underlying graph 
$\xymatrix@-0.5pc{1& 2 \ar@{-}[l]\ar@{-}[r] & 3}$
and let $\ii := (3,1,2,3,1,2)$.
Thus $\GG_\ii$ looks as follows:
$$
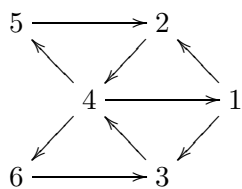
\xymatrix@-0.7pc{
5 \ar[rr] && 2 \ar[dl]\\
& 4 \ar[dl]\ar[ul]\ar[rr] && 1 \ar[ul]\ar[dl]\\
6 \ar[rr] && 3 \ar[ul]
}
$$
The following picture shows the quiver $\GG_{V_\ii}$ of 
$\End_\LL(V_\ii)^\op$ where the vertices corresponding to the modules
$V_k$.
$$
\xymatrix@-0.7pc{
{\bsm&&3\\&2\\1\esm} \ar@{<-}[dr]
&& {\bsm&2\\1\esm} \ar@{<-}[dr]\ar@{<-}[ll]\\
&{\bsm&2\\1&&3\\&2\esm} \ar@{<-}[ur]\ar@{<-}[dr]&& 
{\bsm2\esm} \ar@{<-}[ll]\\
{\bsm1\\&2\\&&3\esm} \ar@{<-}[ur] &&
{\bsm2\\&3\esm}\ar@{<-}[ll]\ar@{<-}[ur]
}
$$
Here is the quiver $\GG_{V_\ii}'$ whose vertices are
the dimension vectors $\dimv_{B_\ii}(\Hom_\LL(V_\ii,V_k))$:
$$
\xymatrix@-0.7pc{
{\bsm1&&1\\&1&&0\\1&&0\esm} \ar@{<-}[dr]
&& {\bsm0&&1\\&1&&0\\1&&0\esm} \ar@{<-}[dr]\ar@{<-}[ll]\\
&{\bsm1&&1\\&2&&1\\1&&1\esm} \ar@{<-}[ur]\ar@{<-}[dr]&& 
{\bsm0&&1\\&1&&1\\0&&1\esm} \ar@{<-}[ll]\\
{\bsm1&&0\\&1&&0\\1&&1\esm} \ar@{<-}[ur] &&
{\bsm1&&0\\&1&&0\\0&&1\esm}\ar@{<-}[ll]\ar@{<-}[ur]
}
$$

Next, let us look at an example of type $\widetilde{\A}_2$.
Thus, let $Q$ be a quiver with underlying graph
$$
\xymatrix@-0.5pc{
& 3 \ar@{-}[dl]\ar@{-}[dr]\\
1 \ar@{-}[rr] && 2
}
$$
and let $\ii := (3,2,1,3,2,1)$.
The quiver $\GG_{V_\ii}$ of $\End_\LL(V_\ii)^\op$ looks as follows:
$$
\xymatrix@-0.7pc{
{\bsm&&&&1\\&1&&2\\2&&3\\&1\esm} \ar@/_1.1cm/@{<-}[ddrr]\ar@{<-}[drrr]
&&
{\bsm1\esm} \ar@{<-}[ll] \\
&
{\bsm&&&&&&1\\&&&1&&2\\1&&2&&3\\&3&&1\\&&2\esm} 
\ar@{<-}[ul]\ar@{<-}[dr]
&&
{\bsm&1\\2\esm} 
\ar@{<-}[ll]\ar@{<-}[ul]
\\
{\bsm&&&&&&&1\\&&&&1&&2\\&1&&2&&3\\2&&3&&1\\&1&&2\\&&3\esm} 
\ar@{<-}[uu]\ar@{<-}[ur] 
&&
{\bsm&&&1\\1&&2\\&3\esm} \ar@{<-}[ll]\ar@{<-}[uu]\ar@{<-}[ur]
}
$$
Here is the quiver $\GG_{V_\ii}'$:
$$
\xymatrix@-0.7pc{
{\bsm 3&&1\\&4&&1\\5&&2\esm} \ar@/_1.1cm/@{<-}[ddrr]\ar@{<-}[drrr]
&&
{\bsm 2&&1\\&3&&1\\3&&2\esm}  \ar@{<-}[ll] \\
&
{\bsm 2&&0\\&3&&1\\4&&1\esm} 
\ar@{<-}[ul]\ar@{<-}[dr]
&&
{\bsm 2&&0\\&2&&1\\3&&1\esm} 
\ar@{<-}[ll]\ar@{<-}[ul]
\\
{\bsm 1&&0\\&2&&0\\3&&1\esm} 
\ar@{<-}[uu]\ar@{<-}[ur] 
&&
{\bsm 1&&0\\&2&&0\\2&&1\esm} 
\ar@{<-}[ll]\ar@{<-}[uu]\ar@{<-}[ur]
}
$$

\subsection{Example (Mutations via dimension vectors)}\label{ex14.3}
Let $Q$ be a quiver with underlying graph
$\xymatrix@-0.5pc{1 \ar@<0.5ex>@{-}[r]\ar@<-0.5ex>@{-}[r] & 
2 \ar@{-}[r] & 3}$ and let
$\ii := (i_7,\ldots,i_1) := (1,3,2,1,3,2,1)$ be a reduced expression.
As before, let
$V_\ii = V_1 \oplus \cdots \oplus V_7$.
The indecomposable $\CC_w$-projective-injectives are
$V_5$, $V_6$ and $V_7$.
Let us compute the dimension vectors $\dimv_{B_\ii}(\Hom_\LL(V_\ii,M_k))$.
\begin{align*}
\dimv(\Delta_1) &= {\bsm9&&3&&1\\&6&&2\\4&&2\esm} &
\dimv(\Delta_2) &= {\bsm6&&2&&0\\&4&&1\\3&&1\esm} &
\dimv(\Delta_3) &= {\bsm2&&0&&0\\&1&&0\\0&&1\esm} &
\dimv(\Delta_4) &= {\bsm3&&1&&0\\&2&&0\\2&&0\esm}\\
\dimv(\Delta_5) &= {\bsm2&&0&&0\\&1&&0\\1&&0\esm} &
\dimv(\Delta_6) &= {\bsm0&&0&&0\\&0&&0\\1&&0\esm} &
\dimv(\Delta_7) &= {\bsm1&&0&&0\\&0&&0\\0&&0\esm}.
\end{align*}
Here is the quiver $\GG_\ii$:
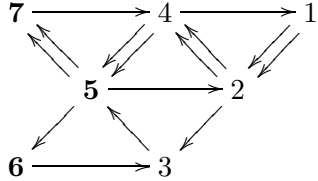
$$
\xymatrix@-0.7pc{
{\bf 7} \ar[rr]
&& 4 \ar@<0.5ex>[dl] \ar@<-0.5ex>[dl]\ar[rr]&&
1 \ar@<0.5ex>[dl] \ar@<-0.5ex>[dl] \\
& {\bf 5} \ar@<0.5ex>[ul]\ar@<-0.5ex>[ul]
\ar[rr]\ar[dl]&& 
2 \ar[dl]\ar@<0.5ex>[ul]\ar@<-0.5ex>[ul]\\
{\bf 6}\ar[rr] &&
3 \ar[ul]
}
$$
The following picture shows the quiver $\GG_{V_\ii}'$.
Its vertices are the dimension vectors of the
$\End_\LL(V_\ii)^\op$-modules $\Hom_\LL(V_\ii,V_k)$.
These dimension vectors can be constructed easily 
using Lemma~\ref{M1}.
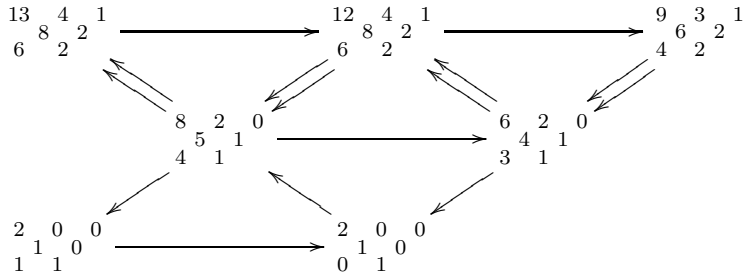
$$
\xymatrix@-0.7pc{
{\bsm13&&4&&1\\&8&&2\\6&&2\esm} \ar[rr]
&& {\bsm12&&4&&1\\&8&&2\\6&&2\esm} \ar@<0.5ex>[dl] \ar@<-0.5ex>[dl]\ar[rr]&&
{\bsm9&&3&&1\\&6&&2\\4&&2\esm}\ar@<0.5ex>[dl] \ar@<-0.5ex>[dl] \\
&{\bsm8&&2&&0\\&5&&1\\4&&1\esm} \ar@<0.5ex>[ul]\ar@<-0.5ex>[ul]
\ar[rr]\ar[dl]&& 
{\bsm6&&2&&0\\&4&&1\\3&&1\esm} \ar[dl]\ar@<0.5ex>[ul]\ar@<-0.5ex>[ul]\\
{\bsm2&&0&&0\\&1&&0\\1&&1\esm}\ar[rr] &&
{\bsm2&&0&&0\\&1&&0\\0&&1\esm}\ar[ul]
}
$$
Now let us mutate the $\LL$-module $V_4$.
We have
$$
\dimv_{B_\ii}(\Hom_\LL(V_\ii,V_4)) = {\bsm12&&4&&1\\&8&&2\\6&&2\esm}.
$$
We have to look at all arrows starting and ending in the corresponding
vertex of $\GG_{V_\ii}'$, and add up the entries of the 
attached dimension vectors, as explained in Section~\ref{dimvectormutation1}.
Since 
$$
\left|{\bsm13&&4&&1\\&8&&2\\6&&2\esm}\right| 
+ 2 \cdot \left|{\bsm6&&2&&0\\&4&&1\\3&&1\esm}\right| =  70 > 69 =
\left|{\bsm9&&3&&1\\&6&&2\\4&&2\esm}\right| + 
2 \cdot \left|{\bsm8&&2&&0\\&5&&1\\4&&1\esm}\right|,
$$ 
we get
$$
\dimv_{B_\ii}(\Hom_\LL(V_\ii,V_4^*)) = 
{\bsm13&&4&&1\\&8&&2\\6&&2\esm} 
+ 2 \cdot {\bsm6&&2&&0\\&4&&1\\3&&1\esm} 
-  {\bsm12&&4&&1\\&8&&2\\6&&2\esm}
= {\bsm13&&4&&0\\&8&&2\\6&&2\esm}
$$
and the quiver $\GG_{\mu_{V_4}(V_\ii)}'$ looks as follows:
$$
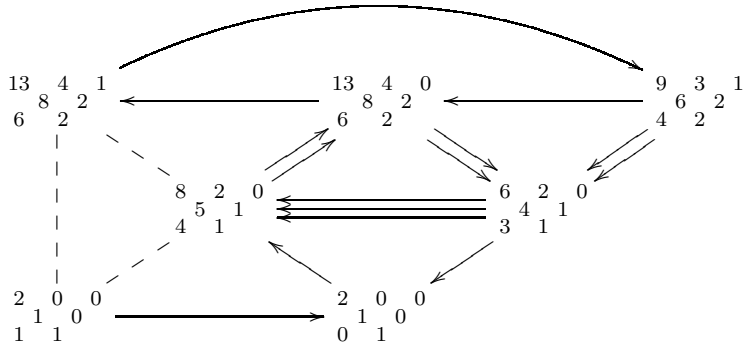
\xymatrix@-0.7pc{
&&&&\\
{\bsm13&&4&&1\\&8&&2\\6&&2\esm} \ar@{<-}[rr]\ar@{--}[dr]\ar@{--}[dd]
&& {\bsm 13&&4&&0\\&8&&2\\6&&2\esm} \ar@<0.5ex>@{<-}[dl] 
\ar@<-0.5ex>@{<-}[dl]\ar@{<-}[rr]&&
{\bsm9&&3&&1\\&6&&2\\4&&2\esm} \ar@/_3pc/@{<-}[llll]\ar@<0.5ex>[dl] 
\ar@<-0.5ex>[dl] \\
&{\bsm8&&2&&0\\&5&&1\\4&&1\esm} \ar@{<-}[dr]
\ar@<0.7ex>@{<-}[rr]\ar@<-0.7ex>@{<-}[rr]\ar@{<-}[rr] && 
{\bsm6&&2&&0\\&4&&1\\3&&1\esm}
\ar@<0.5ex>@{<-}[ul]\ar@<-0.5ex>@{<-}[ul]\ar[dl]\\
{\bsm2&&0&&0\\&1&&0\\1&&1\esm} \ar@{--}[ur] &&
{\bsm2&&0&&0\\&1&&0\\0&&1\esm}\ar@{<-}[ll]
}
$$
Note that we cannot control how the arrows between vertices corresponding
to the three indecomposable
$\CC_w$-projective-injectives
behave under mutation.
But this does not matter, because these arrows are not needed for the
mutation of seeds and clusters.
In the picture, we indicate the missing information by lines
of the form $\xymatrix{\ar@{--}[r]&}$.
This process can be iterated, and our theory says that each of the resulting
dimension vectors determines uniquely a cluster variable.

\subsection{Mutations via $\Delta$-dimension vectors}
\label{dimvectormutation2} 
Using Lemma~\ref{M1} we can explicitly compute the dimension vector of the
$B_\ii$-module $\Delta_s = \Hom_\LL(V_\ii,M_s)$ for all $1 \le s \le r$.
Recall that the $k$th entry of this dimension vector is just
$\dm \Hom_\LL(V_k,M_s)$.
Thus, the $K$-dimension of $\Delta_s$ is
$$
\dim(\Delta_s) = \dm \Hom_\LL(V_\ii,M_s) = \sum_{k=1}^r \dm \Hom_\LL(V_k,M_s).
$$
Define
$$
d_\Delta := (\dim(\Delta_1),\ldots,\dim(\Delta_r)).
$$

Now let $R = R_1 \oplus \cdots \oplus R_r$ 
be a basic $\CC_w$-maximal rigid $\LL$-module, and
suppose that
$R_k$ is not $\CC_w$-projective-injective. 
Then we can mutate $R$ in direction $R_k$.
We obtain two exchange sequences
$$
0 \to R_k \to R' \to R_k^* \to 0
\text{\;\;\; and \;\;\;}
0 \to R_k^* \to R'' \to R_k \to 0
$$
with $R',R'' \in \add(R/R_k)$.

For brevity, set 
$$
{\mathbf d}_s := \dimv_\Delta(\Hom_\LL(V_\ii,R_s))
$$
for all $1 \le s \le r$.
Similarly to the definition of $\GG_R'$ in Section~\ref{dimvectormutation1}
let $\GG_R''$ be the quiver which is obtained from the quiver
of $\End_\LL(R)^\op$ by replacing the vertex corresponding to $R_s$
by the $\Delta$-dimension vector ${\mathbf d}_s$.

For ${\mathbf d} = (d_1,\ldots,d_r)$ and ${\mathbf f} = (f_1,\ldots,f_r)$
in $\Z^r$ define
$$
{\mathbf d} \cdot {\mathbf f} := \sum_{i=1}^r d_if_i.
$$

\begin{Prop}
The $\Delta$-dimension vector of the $B_\ii$-module $\Hom_\LL(V_\ii,R_k^*)$
is
$$
{\mathbf d}_k^* :=
\begin{cases}
- {\mathbf d}_k + \sum_{{\mathbf d}_i \to {\mathbf d}_k} {\mathbf d}_i &
\text{if
$
\sum_{{\mathbf d}_i \to {\mathbf d}_k} {\mathbf d}_i \cdot d_\Delta >
\sum_{{\mathbf d}_k \to {\mathbf d}_j} {\mathbf d}_j \cdot d_\Delta
$},\\
- {\mathbf d}_k +  
\sum_{{\mathbf d}_k \to {\mathbf d}_j} {\mathbf d}_j &
\text{otherwise}.
\end{cases}
$$
Here
the sums are taken over all arrows of the quiver
of $\GG_R''$ which start,
respectively end in the vertex ${\mathbf d}_k$.
\end{Prop}

\begin{proof}
This follows immediately from our results in Section~\ref{dimvectormutation1}
\end{proof}

\subsection{Example (Mutations via $\Delta$-dimension vectors)}
\label{ex14.3cont}
We repeat Example~\ref{ex14.3}, but this time we work
with $\Delta$-dimension vectors.
Let $Q$ and $\ii$ be as before.
The following picture shows the quiver $\GG_{V_\ii}''$.
Its vertices are the $\Delta$-dimension vectors of the
$\End_\LL(V_\ii)^\op$-modules $\Hom_\LL(V_\ii,V_k)$.
$$
\xymatrix@-0.7pc{
{\bsm1&&1&&1\\&0&&0\\0&&0\esm} \ar[rr]
&& {\bsm0&&1&&1\\&0&&0\\0&&0\esm} \ar@<0.5ex>[dl] \ar@<-0.5ex>[dl]\ar[rr]&&
{\bsm0&&0&&1\\&0&&0\\0&&0\esm}\ar@<0.5ex>[dl] \ar@<-0.5ex>[dl] \\
&{\bsm0&&0&&0\\&1&&1\\0&&0\esm} \ar@<0.5ex>[ul]\ar@<-0.5ex>[ul]
\ar[rr]\ar[dl]&& 
{\bsm0&&0&&0\\&0&&1\\0&&0\esm} \ar[dl]\ar@<0.5ex>[ul]\ar@<-0.5ex>[ul]\\
{\bsm0&&0&&0\\&0&&0\\1&&1\esm}\ar[rr] &&
{\bsm0&&0&&0\\&0&&0\\0&&1\esm}\ar[ul]
}
$$
Again, let us mutate the $\LL$-module $V_4$.
We have
$$
\dimv_\Delta(\Hom_\LL(V_\ii,V_4)) = {\bsm0&&1&&1\\&0&&0\\0&&0\esm}.
$$
We have to look at all arrows starting and ending in the corresponding
vertex of $\GG_{V_\ii}''$, and to add up the entries of the 
attached $\Delta$-dimension vectors, as explained in the previous section.
In this example it is clear that the ingoing arrows yield 
the required larger dimension, since the calculation with 
outgoing arrows would 
produce a $\Delta$-dimension vector with negative entries, which is not
possible.
Thus the quiver $\GG_{\mu_{V_4}(V_\ii)}''$ looks as follows:
$$
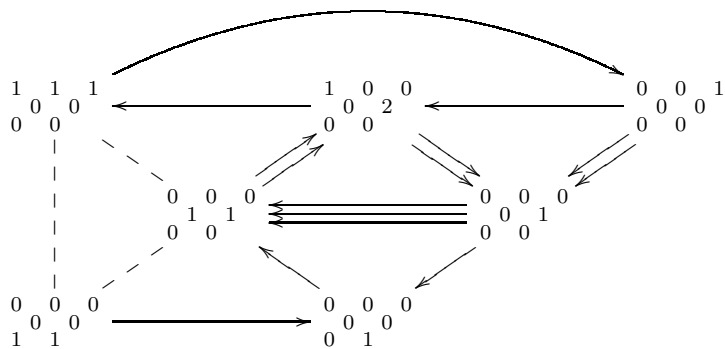
\xymatrix@-0.7pc{
&&&&\\
{\bsm1&&1&&1\\&0&&0\\0&&0\esm} \ar@{<-}[rr]\ar@{--}[dr]\ar@{--}[dd]
&& {\bsm 1&&0&&0\\&0&&2\\0&&0\esm} \ar@<0.5ex>@{<-}[dl] 
\ar@<-0.5ex>@{<-}[dl]\ar@{<-}[rr]&&
{\bsm0&&0&&1\\&0&&0\\0&&0\esm} \ar@/_3pc/@{<-}[llll]\ar@<0.5ex>[dl] 
\ar@<-0.5ex>[dl] \\
&{\bsm0&&0&&0\\&1&&1\\0&&0\esm} \ar@{<-}[dr]
\ar@<0.7ex>@{<-}[rr]\ar@<-0.7ex>@{<-}[rr]\ar@{<-}[rr] && 
{\bsm0&&0&&0\\&0&&1\\0&&0\esm}
\ar@<0.5ex>@{<-}[ul]\ar@<-0.5ex>@{<-}[ul]\ar[dl]\\
{\bsm0&&0&&0\\&0&&0\\1&&1\esm} \ar@{--}[ur] &&
{\bsm0&&0&&0\\&0&&0\\0&&1\esm}\ar@{<-}[ll]
}
$$


\section{A sequence of mutations from $V_\ii$ to $T_\ii$}
\label{section17}


\subsection{The algorithm}\label{algorithm}
Let $\ii := (i_r,\ldots,i_1)$ be a reduced expression of a Weyl group
element.
For $1 \le i,j \le n$ set
$$
q_{ij} := 
\begin{cases}
-c_{ij} & \text{if $i \not= j$},\\
0 & \text{otherwise}.
\end{cases}
$$
(The $c_{ij}$ are the entries of the Cartan matrix $C$ of our Kac-Moody
Lie algebra $\g$, see Section~\ref{KMalg}.
Note that this definition of $q_{ij}$ is equivalent to the one in 
Section~\ref{terminal}.)
As before, we define a quiver $\GG_\ii$ as follows:
The vertices of $\GG_\ii$ are $1,2,\ldots,r$.
For $1 \le s,t \le r$ there are $q_{i_s,i_t}$ arrows from $s$ to $t$ provided
$t^+ \ge s^+ > t > s$.
These are called the {\it ordinary arrows} of $\GG_\ii$.
Furthermore, for each $1 \le s \le r$ there is an arrow $s \to s^-$
provided $s^- > 0$.
These are the {\it horizontal arrows} of $\GG_\ii$.

As before, let
$V_\ii = V_1 \oplus \cdots \oplus V_r$ and
$M_\ii = M_1 \oplus \cdots \oplus M_r$.
We know that the quiver $\GG_\ii$ can be identified with the quiver
$\GG_{V_\ii}$ of the endomorphism algebra $B_\ii = \End_\LL(V_\ii)^\op$.
The vertices of $\GG_{V_\ii}$ are labeled by $V_1,\ldots,V_r$.
More precisely, the vertex $s$ of $\GG_\ii$ corresponds to the
vertex $V_s = M[s,s_{\rm min}]$ of $\GG_{V_\ii}$,
where $1 \le s \le r$.

Recall that for $1 \le j \le n$ and $1 \le k \le r+1$, we defined
\begin{align*}
k[j] &:= |\{ 1 \le s \le k-1 \mid i_s = j \}|,\\
t_j &:= (r+1)[j],\\
k_{\rm min} &:= \min\{ 1 \le s \le r \mid i_s = i_k \}.
\end{align*}
Now we describe an algorithm which yields a sequence
of mutations starting with $\GG_{V_\ii}$ and ending with
$\GG_{T_\ii}$ (see Section~\ref{intervalmod} 
for the definition of $T_\ii$).
The proof is done by induction on $r-n$.

Before going into details let us describe the general idea of this
algorithm. Assume that $Q$ is the linearly oriented quiver
$$
\xymatrix{
m \ar[r] & m-1 \ar[r] & \cdots \ar[r] & 2 \ar[r] & 1
}
$$
of type $\A_m$.
We would like to find a sequence of mutations
which transforms $Q$ into the quiver $Q^\op$
$$
\xymatrix{
m & m-1 \ar[l]& \cdots \ar[l]& 2 \ar[l] & \ar[l] 1
}
$$
with opposite linear orientation.
This can be done by  applying the following $m-1$
sequences of mutations:
$$
Q^{1} := \mu_{m-1} \cdots \mu_2\mu_1(Q),\;\;
Q^{2} := \mu_{m-2} \cdots \mu_2\mu_1(Q^{1}),\;\;
\cdots\;\;,\;\;\
Q^{m-1} := \mu_1(Q^{m-2}).
$$
Now one easily checks that $Q^{m-1} = Q^\op$.
If we delete all ordinary arrows of $\GG_\ii$ we obtain 
a disjoint union of linearly oriented quivers of type
$\A_{m_i}$ for various $m_i \ge 1$.
The main idea of the following algorithm is to
apply a sequence of mutations to $\GG_\ii$
which (in the same way as explained above) reverses the orientation of these 
subquivers of type $\A_{m_i}$ without causing too many changes
for the remaining ordinary arrows.

In the following, we just ignore the symbols of the form 
$M[a,b]$ in case $a<b$.

\noindent{\bf Step 1}:
We mutate the following 
$$
r_1 := t_{i_1}-1-1[i_1]
$$ 
vertices of $\GG_{V_\ii}^0 := \GG_{V_\ii}$ in the given order:
$$
M[1_{\rm min}^{(1[i_1])},1_{\rm min}^{(1[i_1])}],
M[1_{\rm min}^{(1[i_1]+1)},1_{\rm min}^{(1[i_1])}],
M[1_{\rm min}^{(1[i_1]+2)},1_{\rm min}^{(1[i_1])}],
\ldots,
M[1_{\rm min}^{(t_{i_1}-2)},1_{\rm min}^{(1[i_1])}].
$$
Under the identification $\GG_{V_\ii} \equiv \GG_\ii$, 
this sequence of mutations corresponds to the sequence of mutations 
$$
\overrightarrow{\mu_1} := \mu_{1_{\rm min}^{(r_1-1)}} \circ 
\cdots \circ \mu_{1_{\rm min}^{(1)}} \circ \mu_{1_{\rm min}}.
$$
We obtain a new quiver $\GG_{V_\ii}^1$ with $r_1$ new vertices
\begin{multline*}
M[1_{\rm min}^{(1[i_1]+1)},1_{\rm min}^{(1[i_1]+1)}],
M[1_{\rm min}^{(1[i_1]+2)},1_{\rm min}^{(1[i_1]+1)}],
M[1_{\rm min}^{(1[i_1]+3)},1_{\rm min}^{(1[i_1]+1)}],\\
\ldots,
M[1_{\rm min}^{(t_{i_1}-1)},1_{\rm min}^{(1[i_1]+1)}].
\end{multline*}
\noindent{\bf Step 2}:
We mutate the following 
$$
r_2 := t_{i_2}-1-2[i_2]
$$ 
vertices of $\GG_{V_\ii}^1$ in the following order:
$$
M[2_{\rm min}^{(2[i_2])},2_{\rm min}^{(2[i_2])}],
M[2_{\rm min}^{(2[i_2]+1)},2_{\rm min}^{(2[i_2])}],
M[2_{\rm min}^{(2[i_2]+2)},2_{\rm min}^{(2[i_2])}],
\ldots,
M[2_{\rm min}^{(t_{i_2}-2)},2_{\rm min}^{(2[i_2])}].
$$
This mutation sequence corresponds to
$$
\overrightarrow{\mu_2} := \mu_{2_{\rm min}^{(r_2-1)}} \circ 
\cdots \circ \mu_{2_{\rm min}^{(1)}} \circ \mu_{2_{\rm min}}.
$$
We obtain a new quiver $\GG_{V_\ii}^2$ with $r_2$ new vertices
\begin{multline*}
M[2_{\rm min}^{(2[i_2]+1)},2_{\rm min}^{(2[i_2]+1)}],
M[2_{\rm min}^{(2[i_2]+2)},2_{\rm min}^{(2[i_2]+1)}],
M[2_{\rm min}^{(2[i_2]+3)},2_{\rm min}^{(2[i_2]+1)}],\\
\ldots,
M[2_{\rm min}^{(t_{i_2}-1)},2_{\rm min}^{(2[i_2]+1)}].
\end{multline*}
\noindent{\bf Step k}:
We mutate the following 
$$
r_k := t_{i_k}-1-k[i_k]
$$ 
vertices of $\GG_{V_\ii}^{k-1}$ in the following order:
$$
M[k_{\rm min}^{(k[i_k])},k_{\rm min}^{(k[i_k])}],
M[k_{\rm min}^{(k[i_k]+1)},k_{\rm min}^{(k[i_k])}],
M[k_{\rm min}^{(k[i_k]+2)},k_{\rm min}^{(k[i_k])}],
\ldots,
M[k_{\rm min}^{(t_{i_k}-2)},k_{\rm min}^{(k[i_k])}].
$$
This mutation sequence corresponds to
$$
\overrightarrow{\mu_k} := \mu_{k_{\rm min}^{(r_k-1)}} \circ 
\cdots \circ \mu_{k_{\rm min}^{(1)}} \circ \mu_{k_{\rm min}}.
$$
We obtain a new quiver $\GG_{V_\ii}^k$ with $r_k$ 
new vertices
\begin{multline*}
M[k_{\rm min}^{(k[i_k]+1)},k_{\rm min}^{(k[i_k]+1)}],
M[k_{\rm min}^{(k[i_k]+2)},k_{\rm min}^{(k[i_k]+1)}],
M[k_{\rm min}^{(k[i_k]+3)},k_{\rm min}^{(k[i_k]+1)}],\\
\ldots,
M[k_{\rm min}^{(t_{i_k}-1)},k_{\rm min}^{(k[i_k]+1)}].
\end{multline*}
The algorithm stops when
all vertices are of the form
$M[k_{\rm max},k]$.
This will happen after 
$$
r(\ii) := \sum_{j=1}^n \frac{t_j(t_j-1)}{2}
$$ 
mutations.
Define 
$$
\mu_\ii := \overrightarrow{\mu_r} \circ \cdots \circ 
\overrightarrow{\mu_2} \circ \overrightarrow{\mu_1}.
$$
Thus we have
$$
\mu_\ii(V_\ii) = T_\ii.
$$

As an example, assume
$Q$ is a Dynkin quiver of type ${\mathbb E}_8$.
Thus the underlying graph of $Q$ looks as follows:
$$
\xymatrix@-0.7pc{
&& 7 \ar@{-}[d]\\
5 \ar@{-}[r] & 6 \ar@{-}[r] & 8 \ar@{-}[r] & 4 \ar@{-}[r] & 
3 \ar@{-}[r] & 2 \ar@{-}[r] & 1 
}
$$
Let $c:= s_8s_7s_6s_5s_4s_3s_2s_1$.
Then $w := c^{15}$ is the longest element in the Weyl group
$W$ of $Q$, and
$\ii := (8,\ldots,2,1,\ldots,8,\ldots,2,1)$ is a reduced expression
(with 120 entries) of $w$.
We get $t_j = 15$ for all 8 vertices $j$ of $Q$.
Then our algorithm says that starting with $V_\ii$
we reach $T_\ii$ after
$r(\ii) = 8 \cdot 105 = 840$ 
mutations.

We now want to describe what happens to the quiver
$\GG_{V_\ii}^{k-1}$ when we apply the mutation sequence
$\overrightarrow{\mu_k}$.
First, we need some notation:

For each $1 \le j \le n$ let 
\begin{align*}
p_j &:= \min\{ 1 \le s \le r \mid i_s = j \},\\
u_j &:= \min\{ 0,k \le s \le r \mid i_s = j \}.
\end{align*}
Note that $p_j^{(0)} = p_j$.
The sequence
$$
(p_j^{(0)},p_j^{(1)},\ldots,p_j^{(r_{u_j}-1)})
$$
of vertices of $\GG_{V_\ii}^{k-1}$ is called
the $j$-{\it chain} of $\GG_{V_\ii}^{k-1}$, provided
$u_j \not= 0$.
If $u_j = 0$, then we have an empty $j$-chain. 
The sequence
$$
(p_j^{(0)},p_j^{(1)},\ldots,p_j^{(t_j-1)})
$$ 
is the {\it extended} $j$-chain.

Each full subgraph of $\GG_{V_\ii}^{k-1}$ 
given by the vertices of a single extended $j$-chain looks as follows:
$$
\xymatrix@-0.7pc{
p_j^{(t_j-1)} & \cdots \ar[l] &
p_j^{(r_{u_j}+1)} \ar[l] & p_j^{(r_{u_j})} \ar[r]\ar[l] & 
p_j^{(r_{u_j}-1)} \ar[r] & \cdots \ar[r] & 
p_j^{(2)} \ar[r] & p_j^{(1)} \ar[r] & p_j^{(0)}
}
$$
The arrows of the extended $j$-chains ($1 \le j \le n$), are the 
{\it horizontal arrows} of
$\GG_{V_\ii}^{k-1}$.
In the mutation sequence
$$
\overrightarrow{\mu_r} \circ \cdots \circ 
\overrightarrow{\mu_{k+1}} \circ \overrightarrow{\mu_k}
$$
there are no mutations at the vertices 
$p_j^{(r_{u_j})},p_j^{(r_{u_j}+1)},\ldots, p_j^{(t_j-1)}$.
These are called the {\it frozen vertices} of $\GG_{V_\ii}^{k-1}$.

To describe the quiver $\GG_{V_\ii}^k$, it is enough to
study the effect of $\overrightarrow{\mu_k}$ on the $n-1$ full
subgraphs of $\GG_{V_\ii}^k$ which consist of the $i_k$-chain 
together with one extended $j$-chain, where
$1 \le j \le n$ and $j \not= i_k$.

For brevity, set $s = s^{(0)} = k_{\rm min}$, $t = t^{(0)} = p_j$.
Let $q = q_{i_k,j}$ be the number of edges between $i_k$ and $j$
in the underlying graph of $Q$.
The following picture shows how the arrows
between the $i_k$-chain and an extended $j$-chain in $\GG_{V_\ii}^{k-1}$
look like (we have $1 \le j \le n$ with $i_k \not= j$, and
we use the notation $\xymatrix{u \ar[r]|q & v}$
if there are $q$ arrows from $u$ to $v$):
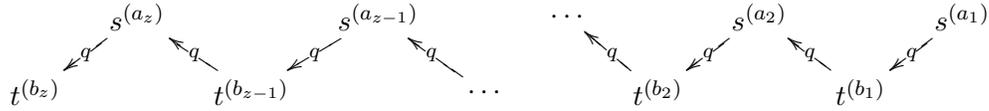
$$
\xymatrix@-1.0pc{
& s^{(a_z)} \ar[dl]|q &&  s^{(a_{z-1})} \ar[dl]|q && 
\cdots && s^{(a_2)} \ar[dl]|q &&
s^{(a_1)} \ar[dl]|q \\
t^{(b_z)} && t^{(b_{z-1})} \ar[ul]|q &&
\cdots \ar[ul]|q && t^{(b_2)} \ar[ul]|q && t^{(b_1)} \ar[ul]|q 
}
$$
Here $s^{(a_i)}$ belongs to the $i_k$-chain,
and $t^{(b_i)}$ belongs to the extended $j$-chain for all $1 \le i \le z$.
(The $q$ arrows from $s^{(a_z)}$ to $t^{(b_z)}$ do not exist necessarily.
But the first $q$ arrows between the $i_k$-chain and the $j$-chain (counted
from the right) always start at the $i_k$-chain.
We do not display any arrows between frozen vertices, they don't play
any role.)

The mutation sequence $\overrightarrow{\mu_k}$ consists of mutations
at the vertices $s^{(0)},s^{(1)},\ldots,s^{(r_k-1)}$.  
By definition, 
$$
\GG_{V_\ii}^k := \overrightarrow{\mu_k}\left(\GG_{V_\ii}^{k-1}\right).
$$
After applying $\overrightarrow{\mu_k}$, the horizontal arrows of the 
$i_k$-chain stay the same, except the arrow $s^{(r_k)} \ra s^{(r_k-1)}$
changes its orientation and becomes  $s^{(r_k)} \leftarrow s^{(r_k-1)}$.
The vertex $s^{(r_k-1)}$ becomes an additional frozen vertex of
$\GG_{V_\ii}^k$.

The arrows between the $i_k$-chain and the $j$-chain change as follows:
$$
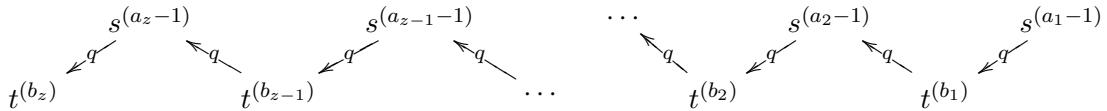
\xymatrix@-1.0pc{
& s^{(a_z-1)} \ar[dl]|q &&  s^{(a_{z-1}-1)} \ar[dl]|q && 
\cdots && s^{(a_2-1)} \ar[dl]|q &&
s^{(a_1-1)} \ar[dl]|q \\
t^{(b_z)} && t^{(b_{z-1})} \ar[ul]|q &&
\cdots \ar[ul]|q && t^{(b_2)} \ar[ul]|q && t^{(b_1)} \ar[ul]|q
}
$$
(In case $s^{(a_1)} = s$, the $q$ arrows from $s^{(a_1-1)}$ to $t^{(b_1)}$
do not exist.)

We illustrate this again in a more explicit example:
Here is a possible subgraph before we apply $\overrightarrow{\mu_k}$,
where $r_k = 8$ and $r_{u_j} = 6$:
$$
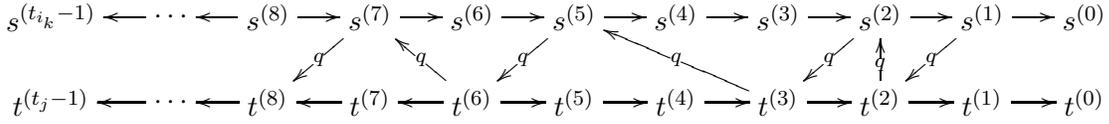
\xymatrix@-0.7pc{
s^{(t_{i_k}-1)} & \cdots \ar[l] & s^{(8)} \ar[l]\ar[r] &
s^{(7)} \ar[r]\ar[dl]|q & s^{(6)} \ar[r] & 
s^{(5)} \ar[r]\ar[dl]|q & s^{(4)} \ar[r] & 
s^{(3)} \ar[r] & s^{(2)} \ar[r]\ar[dl]|q & s^{(1)} \ar[r]\ar[dl]|q & s^{(0)}\\
t^{(t_j-1)} & \cdots \ar[l] & t^{(8)} \ar[l] & t^{(7)} \ar[l] & 
t^{(6)} \ar[l]\ar[r]\ar[ul]|q &
t^{(5)} \ar[r] & t^{(4)} \ar[r] & 
t^{(3)} \ar[r]\ar[ull]|q & t^{(2)} \ar[r]\ar[u]|>>>>q & 
t^{(1)} \ar[r] & t^{(0)}
}
$$
(The numbers $r_k$ and $r_{u_j}$ are determined by the orientation
of the horizontal arrows in the above picture.)

This is how it looks like after we applied  $\overrightarrow{\mu_k}$
to the $r_k$ vertices of the $i_k$-chain:
$$
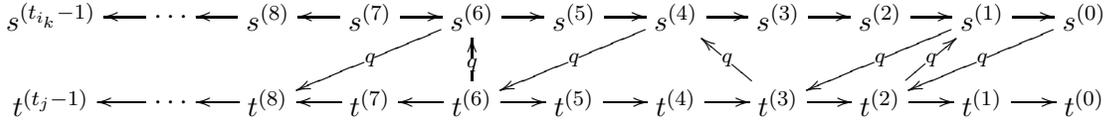
\xymatrix@-0.7pc{
s^{(t_{i_k}-1)} & \cdots \ar[l] & s^{(8)} \ar[l] &
s^{(7)} \ar[l]\ar[r] & s^{(6)} \ar[r]\ar[dll]|q & 
s^{(5)} \ar[r] & s^{(4)} \ar[r]\ar[dll]|q & 
s^{(3)} \ar[r] & s^{(2)} \ar[r] & s^{(1)} \ar[r]\ar[dll]|q & 
s^{(0)}\ar[dll]|q\\
t^{(t_j-1)} & 
\cdots \ar[l] & t^{(8)} \ar[l] & t^{(7)} \ar[l] & 
t^{(6)} \ar[l]\ar[r]\ar[u]|>>>>q &
t^{(5)} \ar[r] & t^{(4)} \ar[r] & 
t^{(3)} \ar[r]\ar[ul]|q & t^{(2)} \ar[r]\ar[ur]|q & 
t^{(1)} \ar[r] & t^{(0)}
}
$$
Again, possible arrows between frozen vertices are not shown.

Note that if we start with our initial 
$\CC_w$-maximal rigid module $V_\ii$, and if we only perform the $r(\ii)$ 
mutations described
in the algorithm, then we obtain the subset
$$
\left\{ M[b,a] \mid 1 \le a \le b \le r, i_a = i_b \right\}
$$
of the set of indecomposable rigid modules of $\CC_w$.
In particular, this subset contains all 
modules
$M_k = M[k,k]$ where $1 \le k \le r$.
The next theorem describes the precise exchange relation
obtained in each of the $r(\ii)$ steps of the algorithm above.

We use our description of mutations via $\Delta$-dimension
vectors from Section~\ref{dimvectormutation2} in order to show that the 
mutation $M[s,s_{\rm min}^{(k[i_s])}]^*$ of 
$M[s,s_{\rm min}^{(k[i_s])}]$ is indeed $M[s^+,s_{\rm min}^{(k[i_s]+1)}]$.

In  formula~(\ref{detform1new}) below we just write $M[b,a]$ instead of 
$\delta_{M[b,a]}$.
(Recall that for any $\LL$-module $X$ and any constructible function
$f \in \MM$ we have $\delta_X(f) := f(X)$.
This defines an element $\delta_X$ in $\MM^*$.)

\begin{Thm}[Generalized determinantal identities]\label{detthmnew}
Let $M_\ii = M_1 \oplus \cdots \oplus M_r$.
Then for $1 \le k,s \le r$ with $i_s = i_k$ we have
\begin{multline}\label{detform1new}
M[s,s_{\rm min}^{(k[i_s])}] 
\cdot
M[s^+,s_{\rm min}^{(k[i_s]+1)}] 
=
M[s^+,s_{\rm min}^{(k[i_s])}] 
\cdot
M[s,s_{\rm min}^{(k[i_s]+1)}] \;\;\;+\\
+ 
\prod_{t^+ \ge s^+ > t> s} M[t,t_{\rm min}^{(k[i_t])}]^{q_{i_si_t}}
\cdot 
\prod_{l^+ \ge s^+ > s > l > s_{\rm min}}
M[l,l_{\rm min}^{(k[i_l])}]^{q_{i_si_l}}.
\end{multline}
\end{Thm}

\begin{proof}
Formula~(\ref{detform1new}) is just an exchange relation
corresponding to the mutation of the module $M[s,s_{\rm min}^{(k[i_s])}]$
with $M[s,s_{\rm min}^{(k[i_s])}]^* = M[s^+,s_{\rm min}^{(k[i_s]+1)}]$.
More precisely, the mutation of $M[s,s_{\rm min}^{(k[i_s])}]$
happens during the mutation sequence $\overrightarrow{\mu_k}$, which is
part of the mutation sequence $\mu_\ii$.
\end{proof}

\begin{Rem}
It is not hard to see that the above theorem can be also stated as
follows:  For $1\leq t< s\leq r$ with $i_s=i_j=i$  we have
\[
M[s,t^+]\,M[s^-,t] = M[s,t]\, M[s^-,t^+] + \prod_{j\in I\setminus\{i\}}
M[s^-(j),t^+(j)]^{q_{ij}},
\]
where in addition to the notation in~\ref{intervalmod} we
set $t^+(j):=\min\{r+1,t+1\leq k\leq r\mid i_k=j\}$.
Fomin and Zelevinsky \cite[Theorem 1.17]{FZ5} prove generalized
determinantal identities associated to pairs of Weyl group elements
for all Dynkin cases (including the non-simply laced cases).
Using the material of Section~\ref{minorssection}, 
formula~(\ref{detform1new}) can be seen as a generalization of 
some of their identities to the symmetric Kac-Moody case.
\end{Rem}

\begin{Cor}
The functions $\delta_{M_1},\ldots,\delta_{M_r}$ are algebraically
independent.
In particular, $\C[\delta_{M_1},\ldots,\delta_{M_r}]$, the subalgebra
of $\MM^*$ generated by the $\delta_{M_k}$'s is just a polynomial
ring in $r$ variables.
\end{Cor}

\begin{proof}
Clearly, the functions 
$\delta_{M[1,1_{\rm min}]},\ldots,\delta_{M[r,r_{\rm min}]}$
are algebraically independent, since $V_k = M[k,k_{\rm min}]$ and
any product of the functions $\delta_{V_1},\ldots,\delta_{V_r}$ lies in the
dual semicanonical basis.
Here we use that $V_\ii$ is rigid and then we apply 
\cite[Theorem 1.1]{GLSSemi1}.
We claim that each function $\delta_{M[b,a]}$ with $1 \le a \le b \le r$ and
$i_a = i_b$ is a rational function in $\delta_{M_1},\ldots,\delta_{M_b}$.
In particular, each $\delta_{V_k}$ is a rational function in
$\delta_{M_1},\ldots,\delta_{M_r}$.
This implies that $\delta_{M_1},\ldots,\delta_{M_r}$ are
algebraically independent.

We prove our claim by induction on $r$ and on the {\it length} 
$l([b,a]) := |\{ a \le k \le b \mid i_k = i_b \}|$ 
of the {\it interval} $[b,a]$.
For $r=1$ the statement is clear.
Also, if $l([b,a]) = 1$, then $M[b,a] = M_b$ and we are done as well.
Thus assume by induction that our claim is true for all intervals
$[d,c]$ of length at most $m$ for some $m \ge 1$.
All intervals of length $m+1$ are of the form
$[b^+,a]$ for some $1 \le a \le b \le r$.
We have $a = b_{\rm min}^{(k[i_b])}$ for some $1 \le k \le r$.
We also  assume by induction that our claim holds for all intervals
$[d,c]$ with $b^+ > d$.
Our formula~(\ref{detform1new}) yields
\begin{multline}
M[b^+,a] = 
\frac{1}{M[b,a^+]} \cdot \left(M[b,a] \cdot M[b^+,a^+]\right) - \\
- \frac{1}{M[b,a^+]} \cdot \left(
\prod_{t^+ \ge b^+ > t> b} M[t,t_{\rm min}^{(k[i_t])}]^{q_{i_bi_t}}
\cdot 
\prod_{l^+ \ge b^+ > b > l > b_{\rm min}}
M[l,l_{\rm min}^{(k[i_l])}]^{q_{i_bi_l}} \right).
\end{multline}
The intervals on the right hand side of this equation all have
either length at most $m$, or they are of the form $[d,c]$ with
$b^+ > d$.
This finishes the proof.
\end{proof}

In fact, 
we will show that for any $\LL$-module
$X \in \CC_w$ we have
$\delta_X \in \C[\delta_{M_1},\ldots,\delta_{M_r}]$,
see Theorem~\ref{proofmain1}.
In particular,
for all
$1 \le k \le r$ the
rational function $\delta_{V_k}$ is
a polynomial
in $\delta_{M_1},\ldots,\delta_{M_r}$.

Another proof of the polynomiality of the functions $\delta_{M[b,a]}$ 
was found by Kedem and Di~Francesco~\cite[Lemma~B.7]{DFK},
using ideas of Fomin and Zelevinsky (in particular \cite[Lemma 4.2]{BFZ}).
We thank these four mathematicians
for communicating their insights to us at MSRI in March 2008.

\subsection{Example}
Let $Q$ be a quiver with underlying graph
$$
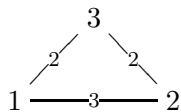
\xymatrix@-0.5pc{
& 3 \ar@{-}[dr]|2\\
1 \ar@{-}[rr]|3 \ar@{-}[ur]|2 && 2
}
$$
Here we use the notation $\xymatrix{i \ar@{-}[r]|a & j}$
if there are $a$ edges between $i$ and $j$.
Let $\ii := (i_{10},\ldots,i_1) := (2,3,2,1,2,1,3,1,2,1)$.
This is a reduced expression for a Weyl group element in $W_Q$.
The quiver $\GG_\ii$ looks as follows:
$$
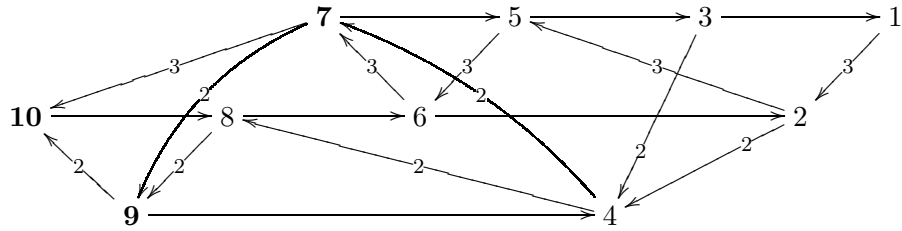
\xymatrix{
&&& {\bf 7} \ar@/_1pc/[ddll]|2 \ar[dlll]|3 \ar[rr] && 5 \ar[dl]|3 \ar[rr] && 
3 \ar[ddl]|>>>>>>>2 \ar[rr] && 1 \ar[dl]|3\\
{\bf 10} \ar[rr] && 8 \ar[dl]|2\ar[rr] && 6 \ar[ul]|3\ar[rrrr] &&&& 
2 \ar[ulll]|3\ar[dll]|<<<<<<2 \\
& {\bf 9} \ar[ul]|2\ar[rrrrr] &&&&& 4 \ar[ullll]|2\ar@/_0.8pc/[uulll]|2
}
$$
For the mutation sequence $\mu_\ii$ we get
\begin{align*}
\mu_\ii &= \overrightarrow{\mu_{10}} \circ \cdots \circ
\overrightarrow{\mu_2} \circ \overrightarrow{\mu_1}\\
&= (\id) \circ (\id) \circ
(\mu_2) \circ (\id) \circ (\mu_6\mu_2) \circ (\mu_1) \circ
(\mu_4) \circ (\mu_3\mu_1) \circ (\mu_8\mu_6\mu_2) \circ (\mu_5\mu_3\mu_1)
\end{align*}
Here are the quivers $\GG_\ii^k$:
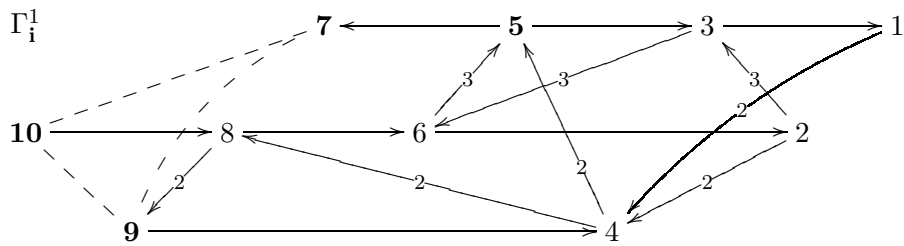
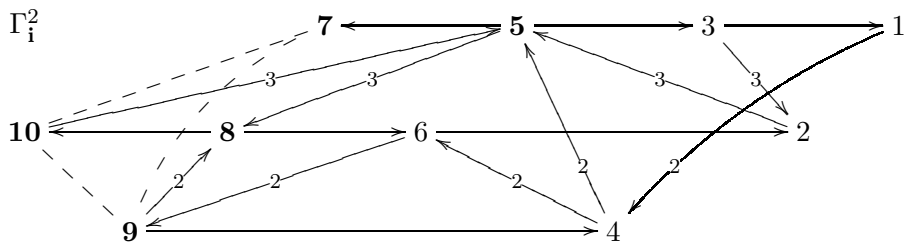
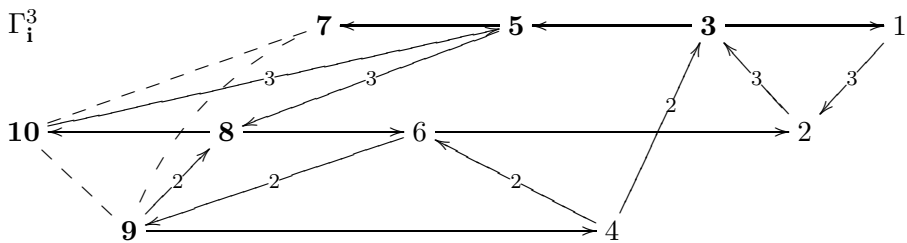
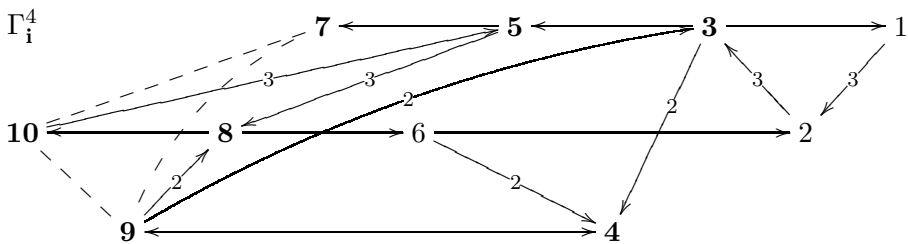
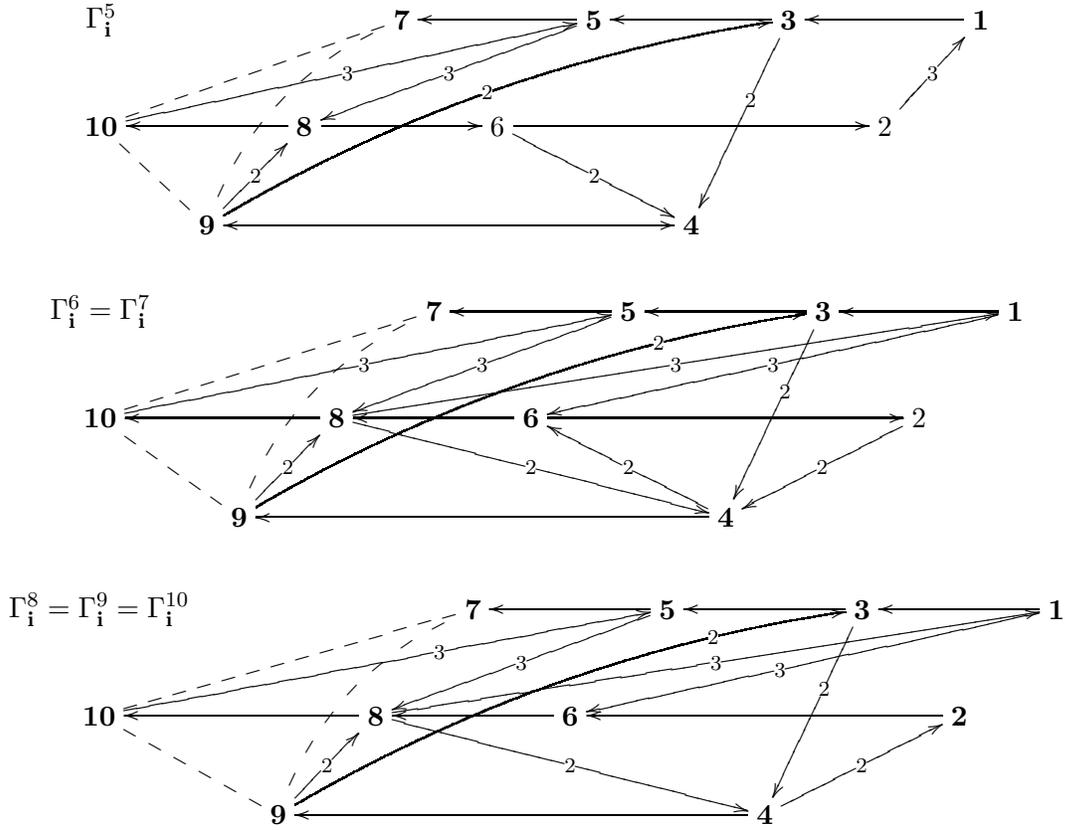
$$
\xymatrix{
\GG_\ii^1&&& {\bf 7} \ar@/_1pc/@{--}[ddll] 
\ar@{--}[dlll] && {\bf 5} \ar[ll] \ar[rr] 
&& 
3 \ar[dlll]|3 \ar[rr] && 1 \ar@/_0.7pc/[ddlll]|2\\
{\bf 10} \ar[rr] && 8 \ar[dl]|2 \ar[rr] && 6 \ar[ur]|3 \ar[rrrr] &&&& 
2 \ar[ul]|3 \ar[dll]|2 \\
& {\bf 9} \ar@{--}[ul] \ar[rrrrr] &&&&& 4 \ar[ullll]|2 \ar[uul]|<<<<<<<2
}
$$

$$
\xymatrix{
\GG_\ii^2&&& {\bf 7} \ar@/_1pc/@{--}[ddll] 
\ar@{--}[dlll] && {\bf 5} \ar[ll] \ar[rr]\ar[dlll]|3 
&&
3 \ar[dr]|3 \ar[rr] && 1 \ar@/_0.7pc/[ddlll]|>>>2\\
{\bf 10} \ar[urrrrr]|3 && {\bf 8} \ar[ll] \ar[rr] && 
6 \ar[dlll]|2 \ar[rrrr] &&&& 
2 \ar[ulll]|3 \\
& {\bf 9} \ar@{--}[ul]\ar[ur]|2 \ar[rrrrr] &&&&& 
4 \ar[ull]|2 \ar[uul]|<<<<<<<2
}
$$

$$
\xymatrix{
\GG_\ii^3 &&& {\bf 7} \ar@/_1pc/@{--}[ddll] 
\ar@{--}[dlll] && {\bf 5} \ar[ll] \ar[dlll]|3 
&&
{\bf 3} \ar[ll] \ar[rr] && 1 \ar[dl]|3\\
{\bf 10} \ar[urrrrr]|3 && {\bf 8} \ar[ll] \ar[rr] && 
6 \ar[dlll]|2 \ar[rrrr] &&&& 
2 \ar[ul]|3 \\
& {\bf 9} \ar@{--}[ul]\ar[ur]|2 \ar[rrrrr] &&&&& 
4 \ar[ull]|2 \ar[uur]|>>>>>>>>>2
}
$$

$$
\xymatrix{
\GG_\ii^4 &&& {\bf 7} \ar@/_1pc/@{--}[ddll] 
\ar@{--}[dlll] && {\bf 5} \ar[ll] \ar[dlll]|3 
&&
{\bf 3} \ar[ll] \ar[rr]\ar[ddl]|<<<<<<<<<2 && 1 \ar[dl]|3\\
{\bf 10} \ar[urrrrr]|3 && {\bf 8} \ar[ll] \ar[rr] && 
6 \ar[drr]|2 \ar[rrrr] &&&& 
2 \ar[ul]|3 \\
& {\bf 9} \ar@{--}[ul]\ar[ur]|2 \ar@/^1pc/[uurrrrrr]|2 \ar[rrrrr] &&&&& 
{\bf 4} \ar[lllll]
}
$$

$$
\xymatrix{
\GG_\ii^5 &&& {\bf 7} \ar@/_1pc/@{--}[ddll] 
\ar@{--}[dlll] && {\bf 5} \ar[ll] \ar[dlll]|3 
&&
{\bf 3} \ar[ll] \ar[ddl]|<<<<<<<<<2 && {\bf 1} \ar[ll]\\
{\bf 10} \ar[urrrrr]|3 && {\bf 8} \ar[ll] \ar[rr] && 
6 \ar[drr]|2 \ar[rrrr] &&&& 
2 \ar[ur]|3 \\
& {\bf 9} \ar@{--}[ul]\ar[ur]|2 \ar@/^1pc/[uurrrrrr]|2 \ar[rrrrr] &&&&& 
{\bf 4} \ar[lllll]
}
$$

$$
\xymatrix{
\GG_\ii^6 = \GG_\ii^7 &&& {\bf 7} \ar@/_1pc/@{--}[ddll] 
\ar@{--}[dlll] && {\bf 5} \ar[ll] \ar[dlll]|3 
&&
{\bf 3} \ar[ll] \ar[ddl]|<<<<<<<<<2 && {\bf 1} \ar[dlllll]|3 \ar[ll]\\
{\bf 10} \ar[urrrrr]|3 && {\bf 8} \ar[ll] \ar[urrrrrrr]|3 
\ar[drrrr]|2 && 
{\bf 6} \ar[ll] \ar[rrrr] &&&& 
2 \ar[dll]|2 \\
& {\bf 9} \ar@{--}[ul]\ar[ur]|2 \ar@/^1pc/[uurrrrrr]|>>>>>>>>>>2 &&&&& 
{\bf 4} \ar[lllll]\ar[ull]|2
}
$$

$$
\xymatrix{
\GG_\ii^8 = \GG_\ii^9 = \GG_\ii^{10} &&& {\bf 7} \ar@/_1pc/@{--}[ddll] 
\ar@{--}[dlll] && {\bf 5} \ar[ll] \ar[dlll]|3 
&&
{\bf 3} \ar[ll] \ar[ddl]|<<<<<<<<<2 && {\bf 1} 
\ar[dlllll]|>>>>>>>>>>>>>>>>>>>>>>>>3 \ar[ll]\\
{\bf 10} \ar[urrrrr]|>>>>>>>>>>>>>>>>>>>>>>>>>>3 
&& {\bf 8} \ar[ll] \ar[urrrrrrr]|3 
\ar[drrrr]|2 && 
{\bf 6} \ar[ll] &&&& 
{\bf 2} \ar[llll]\\
& {\bf 9} \ar@{--}[ul]\ar[ur]|2 \ar@/^1pc/[uurrrrrr]|>>>>>>>>>2 &&&&& 
{\bf 4} \ar[lllll]\ar[urr]|2
}
$$
Applying formula~(\ref{detform1new}) to 
$M[s,s_{\rm min}^{(k[i_s])}] := M[6,6_{\rm min}^{(2[i_6])}] = M[6,2]$ 
we get the following:
$$
M[6,2] \cdot M[8,6] = M[8,2] \cdot M[6,6] +
M[7,3]^3  \cdot M[4,4]^2 \;\;\;\; (s=6,k=2).
$$
Thus, we have
$$
M[8,2] = \frac{1}{M[6,6]} 
\left(M[6,2] \cdot M[8,6] - M[7,3]^3  \cdot M[4,4]^2  \right).
$$
Similarly, we obtain
\begin{align*}
M[2,2] \cdot M[6,6] &= M[6,2] +
M[5,3]^3  \cdot M[4,4]^2 & (s=2,k=2),\\
M[6,6] \cdot M[8,8] &= M[8,6] +
M[7,7]^3 & (s=6,k=6),\\
M[5,3] \cdot M[7,5] &= M[7,3] \cdot M[5,5] +
M[6,6]^3  \cdot M[4,4]^2 & (s=5,k=3),\\
M[3,3] \cdot M[5,5] &= M[5,3] +
M[4,4]^2 & (s=3,k=3),\\
M[5,5] \cdot M[7,7] &= M[7,5] +
M[6,6]^3 & (s=5,k=5).
\end{align*}
By our double induction (on $r$ and on the length of the intervals
$[b,a]$), in each of the above equations, we can write the functions
$M[6,2], M[8,6], M[7,3], M[5,3]$ and $M[7,5]$, respectively,
as a rational function of the functions
appearing in the same equation.
Now one can use these equations to express $\delta_{M[8,2]}$ as a rational 
function in $\delta_{M_1},\ldots,\delta_{M_8}$.
Remarkably, this rational function is a polynomial.

Finally, we display the dimension vectors of the modules $M_1,\ldots,M_8$:
\begin{align*}
\beta_\ii(1) &= {\bsm &0\\1&&0\esm} &
\beta_\ii(2) &= {\bsm &1\\3&&0\esm} &
\beta_\ii(3) &= {\bsm &3\\8&&0\esm} &
\beta_\ii(4) &= {\bsm &8\\24&&1\esm}\\
\beta_\ii(5) &= {\bsm &13\\40&&2\esm} &
\beta_\ii(6) &= {\bsm &63\\189&&8\esm} &
\beta_\ii(7) &= {\bsm &176\\527&&22\esm} &
\beta_\ii(8) &= {\bsm &465\\1392&&58\esm}
\end{align*}
As an exercise, the reader can compute $\beta_\ii(9)$ and $\beta_\ii(10)$.

Exchange equations are always homogeneous.
For example,
$$
M[5,3] \cdot M[7,5] = M[7,3] \cdot M[5,5] +
M[6,6]^3  \cdot M[4,4]^2
$$
is an equation of degree ${\bsm &205\\615&&26 \esm}$.

\subsection{The shift functor in $\stCC_w$ via mutations}\label{shiftcomb}
Fix a reduced expression $\ii = (i_r,\ldots,i_1)$ of some Weyl group
element $w$.
As before,
let $T_\ii := I_w \oplus \Omega_w^{-1}(V_\ii)$.
Define
$$
W_\ii := I_w \oplus \Omega_w(V_\ii).
$$
In Section~\ref{algorithm} we defined a sequence of mutations 
$$
\mu_\ii = \overrightarrow{\mu_r} \circ \cdots \circ \overrightarrow{\mu_1}
= \mu_{s_{r(\ii)}} \circ \cdots \circ \mu_{s_2} \circ \mu_{s_1}
$$
where $1 \le s_p \le r$ for all $p$, such that
$$
\mu_\ii(V_\ii) =\mu_{s_{r(\ii)}} \circ \cdots \circ \mu_{s_2}
\circ \mu_{s_1}(V_\ii)= T_\ii
\text{\;\;\; and \;\;\;}
\mu_\ii^{-1}(T_\ii) =\mu_{s_1} \circ \mu_{s_2} \circ \cdots
\circ \mu_{s_{r(\ii)}}(T_\ii)= V_\ii.
$$
Clearly, if $R$ is a basic $\CC_w$-maximal rigid module such that 
$R = \mu_{p_t} \circ \cdots \circ \mu_{p_1}(V_\ii)$,
then we have 
$R = \mu_{p_t} \circ \cdots \circ \mu_{p_1} \circ \mu_\ii^{-1}(T_\ii)$.

Now define an involution 
$$(-)^*\df 
\{ 1,\ldots,r \} \setminus \{ 1 \le k \le r \mid k^+ = r+1 \} \to 
\{ 1,\ldots,r \} \setminus \{ 1 \le k \le r \mid k^+ = r+1 \}
$$ 
by
$$
\left(k_{\rm min}^{(m)}\right)^* := k_{\rm min}^{(t_j-2-m)},
$$
where $j := i_{k_{\rm min}}$.
Observe that every $1 \le s \le r$ can be written as
$k_{\rm min}^{(m)}$ for some unique $k$ (namely $k=s$) and some
unique $0 \le m \le t_j-1$.
The following picture illustrates how $(-)^*$ permutes the vertices
of $\GG_\ii$:
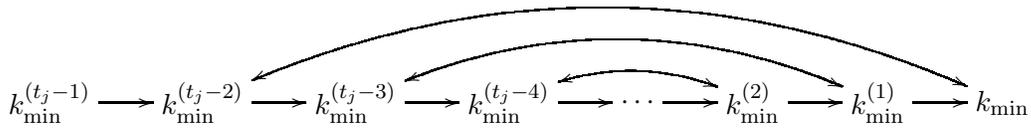
$$
\xymatrix@-0.3pc{
&&\\
k_{\rm min}^{(t_j-1)} \ar[r] & k_{\rm min}^{(t_j-2)} \ar[r] & 
k_{\rm min}^{(t_j-3)} \ar[r] & k_{\rm min}^{(t_j-4)} \ar[r] & 
\cdots \ar[r] & 
k_{\rm min}^{(2)} \ar[r]\ar@/_1pc/@{<->}[ll] & 
k_{\rm min}^{(1)} \ar[r]\ar@/_2pc/@{<->}[llll]& 
k_{\rm min} \ar@/_3pc/@{<->}[llllll]
}
$$
Set
$$
(\mu_\ii^{-1})^* := \mu_{s_1^*} \circ \mu_{s_2^*} \circ \cdots
\circ \mu_{s_{r(\ii)}^*}.
$$

\begin{Prop}\label{Omegamutation}
Let $R$ be a basic $\CC_w$-maximal rigid module which is mutation equivalent
to $V_\ii$.
Then $I_w \oplus \Omega_w^{-1}(R)$ and
$I_w \oplus \Omega_w(R)$ are mutation equivalent to $V_\ii$.
More precisely, let 
$$
R = \mu_{z_t} \circ \cdots \circ \mu_{z_1}(V_\ii)
\text{\;\;\; and \;\;\;}
R = \mu_{q_u} \circ \cdots \circ \mu_{q_1}(T_\ii).
$$
Then we have
$$
I_w \oplus \Omega_w^{-1}(R) = 
\mu_{z_t^*} \circ \cdots \circ \mu_{z_1^*}(T_\ii)
\text{\;\;\; and \;\;\;}
I_w \oplus \Omega_w(R) = 
\mu_{q_u^*} \circ \cdots \circ \mu_{q_1^*}(V_\ii). 
$$
\end{Prop}

Besides $\Omega_w^{-1}(R)$, we can also compute $\Omega_w(R)$
by just knowing a sequence of mutations from $V_\ii$ to $R$.
This works because $V_\ii$ and $T_\ii$ are connected via a known sequence of
mutations, namely $\mu_\ii$ if we start at $V_\ii$, and $\mu_\ii^{-1}$ if we
start at $T_\ii$.
The following picture illustrates the situation:
$$
\xymatrix{
W_\ii\ar@{~>}[dd]^{\mu_{z_t^*} \circ \cdots \circ \mu_{z_1^*}}  
&& V_\ii 
\ar@{~>}[ll]_{(\mu_\ii^{-1})^*}
\ar@{~>}[rr]^{\mu_{\ii}} \ar@{~>}[dd]^{\mu_{z_t} \circ \cdots
\circ \mu_{z_1}} 
&& T_\ii \ar@{~>}[dd]^{\mu_{z_t^*} \circ \cdots \circ \mu_{z_1^*}} \\
&&\\
I_w \oplus \Omega_w(R) && R && I_w \oplus \Omega_w^{-1}(R)
}
$$

\begin{proof}
As before, 
for $1 \le j \le n$ set $p_j := \min\{ 1 \le s \le r \mid i_s = j \}$.
Note that $p_j^{(t_j-1)} = {p_j}_{\rm max}$ and $p_j = {p_j}_{\rm min}$.
In the following pictures we display only the relevant horizontal arrows.
The quiver $\GG_{V_\ii}$ looks as follows:
$$
\xymatrix@-0.7pc{
M[p_1^{(t_1-1)},p_1] \ar@{.}[dd] \ar[r] & M[p_1^{(t_1-2)},p_1] 
\ar[r]\ar@{.}[dd] & 
\cdots \ar[r] & \ar@{.}[dd] M[p_1^{(1)},p_1]  
\ar[r] &\ar@{.}[dd] M[p_1,p_1]\\
&&&\\
M[p_n^{(t_n-1)},p_n] \ar[r] & M[p_n^{(t_n-2)},p_n] \ar[r] & \cdots  
\ar[r] & M[p_n^{(1)},p_n]
\ar[r] & M[p_n,p_n] 
}
$$
Next, we display the quiver $\GG_{T_\ii}$:
$$
\xymatrix@-0.7pc{
M[{p_1}_{\rm max},p_1] \ar@{.}[dd]& M[{p_1}_{\rm max},p_1^{(1)}]\ar[l] 
\ar@{.}[dd] &
M[{p_1}_{\rm max},p_1^{(2)}]\ar[l] \ar@{.}[dd]&
\ar[l] \cdots & 
\ar[l] \ar@{.}[dd] M[{p_1}_{\rm max},p_1^{(t_1-1)}]  \\
&&&\\
M[{p_n}_{\rm max},p_n] & M[{p_n}_{\rm max},p_n^{(1)}]\ar[l] &
M[{p_n}_{\rm max},p_n^{(2)}]\ar[l] &
\ar[l] \cdots & 
\ar[l] M[{p_n}_{\rm max},p_n^{(t_n-1)}]
}
$$
We know that $I_w \oplus \Omega_w^{-1}(V_\ii) = T_\ii$.
In particular, we have
$$ 
\Omega_w^{-1}(M[p_j^{(s-1)},p_j]) = M[{p_j}_{\rm max},p_j^{(s)}] 
$$
for all $1 \le j \le n$ and $1 \le s \le t_j-1$.
Thus $\GG_{T_\ii}$ looks like this:
$$
\xymatrix@-0.7pc{
M[p_1^{(t_1-1)},p_1] \ar@{.}[dd]& \Omega_w^{-1}(M[p_1,p_1])\ar[l] \ar@{.}[dd] &
 \Omega_w^{-1}(M[p_1^{(1)},p_1])\ar[l] \ar@{.}[dd]&
\ar[l] \cdots & 
\ar[l] \ar@{.}[dd]\Omega_w^{-1}(M[p_1^{(t_1-2)},p_1])  \\
&&&\\
M[p_n^{(t_n-1)},p_n] &\Omega_w^{-1}(M[p_n,p_n]) \ar[l] & 
\Omega_w^{-1}(M[p_n^{(1)},p_n]) \ar[l] &
\ar[l] \cdots & 
\ar[l] \Omega_w^{-1}(M[p_n^{(t_n-2)},p_n]) 
}
$$
The $n$ vertices of the form 
$M[p_j^{(t_j-1)},p_j]$ at the ``left'' of both quivers 
$\GG_{V_\ii}$ and $\GG_{T_\ii}$
are frozen vertices, to all other vertices we can apply the mutation
operation.

Now let 
$$
0 \to T_k \to T' \to T_k^* \to 0
$$ 
be an exchange sequence
associated to the cluster algebra $\RR(\CC_w,V_\ii)$.
This yields an {\it exchange triangle}
$T_k \to T' \to T_k^* \to T_k[1]$
in the stable category $\stCC_w$.
Note that $T_k[1] = \Omega_w^{-1}(T_k)$.
It follows that
$
T_k[1] \to T'[1] \to T_k^*[1] \to T_k[2]
$
is an exchange triangle as well.
There is an associated exchange sequence
$$
0 \to \Omega_w^{-1}(T_k) \to I \oplus \Omega_w^{-1}(T') \to
\Omega_w^{-1}(T_k^*) \to 0
$$
where $I$ is some module in $\add(I_w)$.
Thus, if we mutate the basic $\CC_w$-maximal rigid module
$I_w \oplus \Omega_w^{-1}(T)$ in direction $\Omega_w^{-1}(T_k)$,
we obtain $(\Omega_w^{-1}(T_k))^* = \Omega_w^{-1}(T_k^*)$.
We argue similarly to show that the mutation of $I_w \oplus \Omega_w(T)$
in direction $\Omega_w(T_k)$ gives
$(\Omega_w(T_k))^* = \Omega_w(T_k^*)$.
This finishes the proof.
\end{proof}

\begin{Cor}
If a $\LL$-module $R$ is $V_\ii$-reachable, then 
$\Omega_w^z(R)$ is $V_\ii$-reachable for all $z \in \Z$.
\end{Cor}

\subsection{Example}
Let $Q$ be a quiver with underlying graph 
$\xymatrix@-0.5pc{1 \ar@{-}[r] & 2 \ar@{-}[r] & 3 \ar@{-}[r] & 4}$ 
and let $\ii := (i_{10},\ldots,i_1) := (1,2,1,3,2,1,4,3,2,1)$.
Then we get
$$
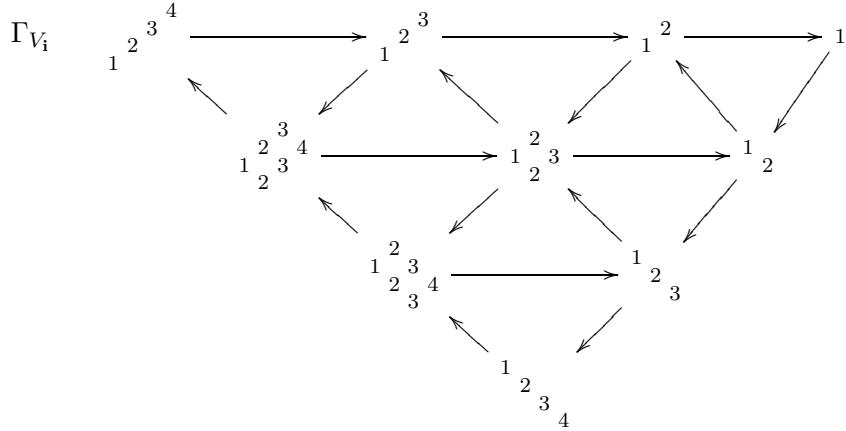
\xymatrix@-0.9pc{
\GG_{V_\ii} & {\bsm&&&4\\&&3\\&2\\1\esm} \ar[rr] &&
{\bsm&&3\\&2\\1\esm} \ar[rr]\ar[dl] && 
{\bsm&2\\1\esm} \ar[dl]\ar[rr] && {\bsm1\esm} \ar[dl]\\
&& {\bsm&&3\\&2&&4\\1&&3\\&2\esm} \ar[ul]\ar[rr]
&& {\bsm&2\\1&&3\\&2\esm} \ar[ul]\ar[rr]\ar[dl] && 
{\bsm1\\&2\esm} \ar[dl]\ar[ul]\\
&&& {\bsm&2\\1&&3\\&2&&4\\&&3\esm} \ar[ul]\ar[rr] && 
{\bsm1\\&2\\&&3\esm} \ar[dl]\ar[ul]\\
&&&& {\bsm1\\&2\\&&3\\&&&4\esm} \ar[ul]
}
$$
and
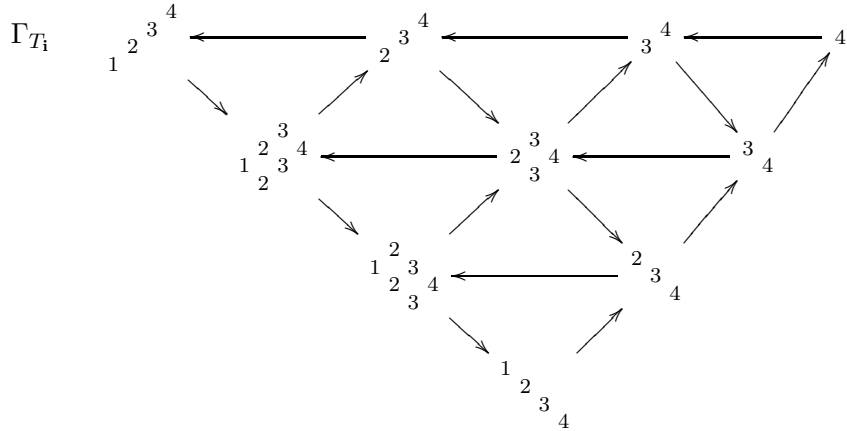
$$
\xymatrix@-0.9pc{
\GG_{T_\ii} & {\bsm&&&4\\&&3\\&2\\1\esm}\ar[dr] &&
{\bsm&&4\\&3\\2\esm} \ar[dr]\ar[ll] && 
{\bsm&4\\3\esm} \ar[dr]\ar[ll] && 
{\bsm4\esm} \ar[ll]\\
&& {\bsm&&3\\&2&&4\\1&&3\\&2\esm} \ar[dr]\ar[ur]
&& {\bsm&3\\2&&4\\&3\esm} \ar[ur]\ar[dr]\ar[ll] && 
{\bsm3\\&4\esm} \ar[ll]\ar[ur] \\
&&&  {\bsm&2\\1&&3\\&2&&4\\&&3\esm} \ar[dr]\ar[ur] 
&& {\bsm2\\&3\\&&4\esm} \ar[ur]\ar[ll]\\
&&&& {\bsm1\\&2\\&&3\\&&&4\esm} \ar[ur]
}
$$
Again, we identify the vertices of $\GG_{V_\ii}$ and $\GG_{T_\ii}$
with the indecomposable direct summands
of $V_\ii$ and $T_\ii$, respectively.
As before, we identify $\GG_{V_\ii}$ and the quiver 
$\GG_\ii$, which looks as follows:
$$
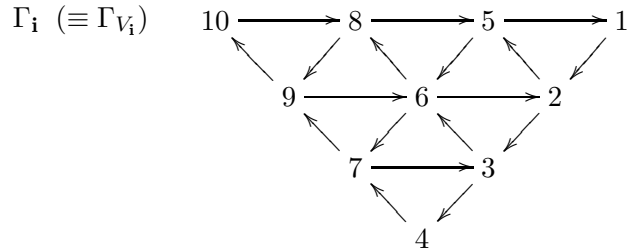
\xymatrix@-0.9pc{
\GG_\ii \;\;(\equiv \GG_{V_\ii}) & 10 \ar[rr] &&
8 \ar[rr]\ar[dl] && 
5 \ar[dl]\ar[rr] && 1 \ar[dl]\\
&& 9 \ar[ul]\ar[rr]
&& 6 \ar[ul]\ar[rr]\ar[dl] && 
2 \ar[dl]\ar[ul]\\
&&& 7 \ar[ul]\ar[rr] && 
3 \ar[dl]\ar[ul]\\
&&&& 4 \ar[ul]
}
$$
We have
\begin{align*}
\mu_\ii 
&= \overrightarrow{\mu_{10}} \circ \cdots 
\circ \overrightarrow{\mu_2} \circ \overrightarrow{\mu_1}\\
&= ({\rm id}) \circ
({\rm id}) \circ
(\mu_1) \circ
({\rm id}) \circ
(\mu_2) \circ 
(\mu_5 \circ \mu_1) \circ
({\rm id}) \circ 
(\mu_3) \circ 
(\mu_6 \circ \mu_2) \circ 
(\mu_8 \circ \mu_5 \circ \mu_1).
\end{align*}
Mutation of $V_\ii$ at $V_5$  yields the following quiver:
$$
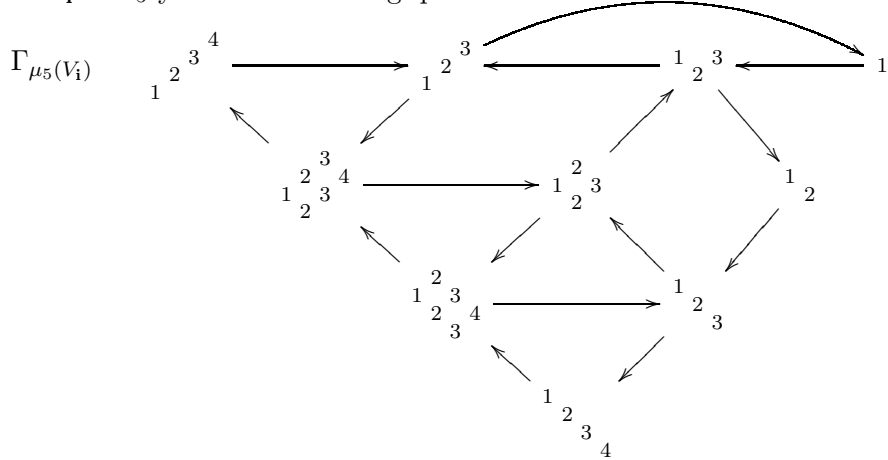
\xymatrix@-0.9pc{
\GG_{\mu_5(V_\ii)} & {\bsm&&&4\\&&3\\&2\\1\esm} \ar[rr] &&
{\bsm&&3\\&2\\1\esm} \ar[dl]\ar@/^2pc/[rrrr] && 
{\bsm1&&3\\&2\esm} \ar[dr]\ar[ll] && {\bsm1\esm} \ar[ll]\\
&& {\bsm&&3\\&2&&4\\1&&3\\&2\esm} \ar[ul]\ar[rr]
&& {\bsm&2\\1&&3\\&2\esm} \ar[dl]\ar[ur] && 
{\bsm1\\&2\esm} \ar[dl]\\
&&& {\bsm&2\\1&&3\\&2&&4\\&&3\esm} \ar[ul]\ar[rr] && 
{\bsm1\\&2\\&&3\esm} \ar[dl]\ar[ul]\\
&&&& {\bsm1\\&2\\&&3\\&&&4\esm} \ar[ul]
}
$$
The associated exchange sequences are
$$
0 \to  ({\bsm1&&3\\&2\esm}) \to {\bsm1\esm} \oplus ({\bsm&2\\1&&3\\&2\esm})  
\to ({\bsm&2\\1\esm}) \to 0
\text{\;\;\; and \;\;\;}
0 \to ({\bsm&2\\1\esm}) \to 
({\bsm1\\&2\esm}) \oplus ({\bsm&&3\\&2\\1\esm}) \to  ({\bsm1&&3\\&2\esm}) \to 0
$$
Next, we mutate at $V_6$.
The exchange sequences looks as follows:
\begin{multline*}
0 \to ({\bsm1&&3\\&2&&4\\&&3\esm}) \to 
({\bsm1&&3\\&2\esm}) \oplus ({\bsm&2\\1&&3\\&2&&4\\&&3\esm})
\to ({\bsm&2\\1&&3\\&2\esm}) \to 0\\
\text{\;\;\ and \;\;}
0 \to ({\bsm&2\\1&&3\\&2\esm}) \to 
({\bsm1\\&2\\&&3\esm})\oplus ({\bsm&&3\\&2&&4\\1&&3\\&2\esm}) 
\to ({\bsm1&&3\\&2&&4\\&&3\esm}) \to 0
\end{multline*}
Set $R := (\mu_6 \circ \mu_5)(V_\ii)$.
Thus we have $R = R_5 \oplus R_6 \oplus V_\ii/(V_5 \oplus V_6)$ 
with 
$$
R_5 = \bsm1&&3\\&2\esm 
\text{\;\;\; and \;\;\;}
R_6 = \bsm1&&3\\&2&&4\\&&3\esm.
$$

To calculate $\Omega_w^{-1}(R)$, we have to compute
$(\mu_{6^*} \circ \mu_{5^*})(T_\ii)$.
Mutation of $T_\ii$ at $5^* = 5$ yields the following quiver:
$$
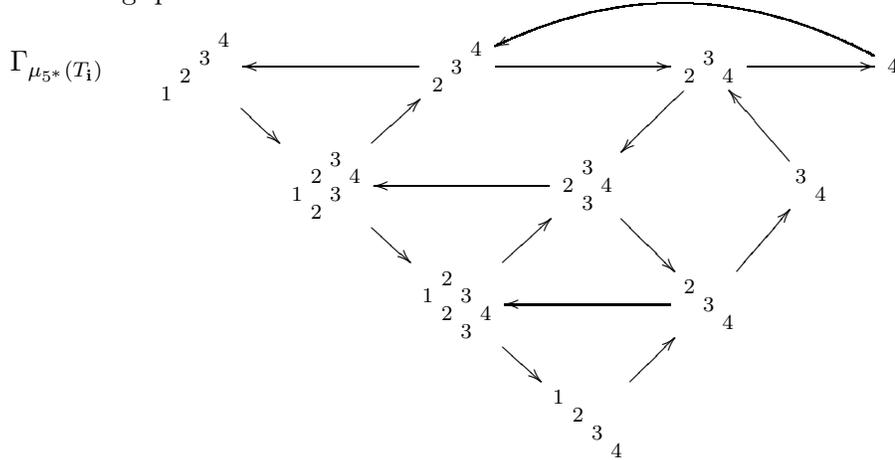
\xymatrix@-0.9pc{
\GG_{\mu_{5^*}(T_\ii)} & {\bsm&&&4\\&&3\\&2\\1\esm}\ar[dr] &&
{\bsm&&4\\&3\\2\esm} \ar[ll]\ar[rr] && 
{\bsm&3\\2&&4\esm} \ar[dl]\ar[rr] && 
{\bsm4\esm} \ar@/_2pc/[llll]\\
&& {\bsm&&3\\&2&&4\\1&&3\\&2\esm} \ar[dr]\ar[ur]
&& {\bsm&3\\2&&4\\&3\esm} \ar[dr]\ar[ll] && 
{\bsm3\\&4\esm} \ar[ul] \\
&&&  {\bsm&2\\1&&3\\&2&&4\\&&3\esm} \ar[dr]\ar[ur] 
&& {\bsm2\\&3\\&&4\esm} \ar[ur]\ar[ll]\\
&&&& {\bsm1\\&2\\&&3\\&&&4\esm} \ar[ur]
}
$$
The associated exchange sequences are
$$
0 \to {\bsm&3\\2&&4\esm} \to 
{\bsm3\\&4\esm} \oplus {\bsm&&4\\&3\\2\esm} \to {\bsm&4\\3\esm} \to 0
\text{\;\;\; and \;\;\;}
0 \to {\bsm&4\\3\esm} \to {\bsm4\esm} \oplus {\bsm&3\\2&&4\\&3\esm} 
\to {\bsm&3\\2&&4\esm} \to 0.
$$
Next, we mutate at $6^* = 2$.
The exchange sequences looks as follows:
$$
0 \to {\bsm2\esm} \to {\bsm&3\\2&&4\esm} \to {\bsm3\\&4\esm} \to 0
\text{\;\;\; and \;\;\;}
0 \to {\bsm3\\&4\esm} \to {\bsm2\\&3\\&&4\esm} \to {\bsm2\esm}\to 0.
$$
We get 
$\Omega_w^{-1}(R) := \Omega_w^{-1}(R_5) \oplus \Omega_w^{-1}(R_6) 
\oplus T_\ii/(T_5 \oplus T_6)$
with
$$
\Omega_w^{-1}(R_5) = {\bsm&3\\2&&4\esm}
\text{\;\;\; and \;\;\;} 
\Omega_w^{-1}(R_6) = {\bsm2\esm}.
$$


\section{Irreducible components associated to $\CC_w$}
\label{irredcomp}


\subsection{Module varieties}
Let $\GG := (\GG_0,\GG_1,s,t)$ be a finite quiver with vertex set 
$\GG_0 = \{ 1,\ldots, r\}$,
arrow set $\GG_1$ and  maps $s,t\df \GG_1 \ra \GG_0$ which map 
an arrow $a$ to its start vertex $s(a)$ and its terminal vertex $t(a)$, 
respectively.
In this section, we
interpret dimension vectors $\ff = (f_1,\ldots,f_r)$ for $\GG$
as maps $\ff \df \GG_0 \ra \N$. 
We consider the affine space
$$
\md(\C\GG,\ff) = \rep(\GG,\ff) =
\prod_{a \in \GG_1} \C^{\ff(t(a)) \times \ff(s(a))}
$$
of representations of $\GG$ with dimension vector $\ff$.
Here $\C^{p \times q}$ denotes the vector space of $(p \times q)$-matrices
with entries in $\C$.
This coincides with our definitions 
in Section~\ref{prenil}, except that we now work with
spaces $\C^{p \times q}$ of matrices rather than spaces 
$\Hom_\C(\C^q,\C^p)$ of linear maps.
So each element in $\md(\C\GG,\ff)$ is of the form $M = (M(a))_{a \in \GG_1}$
where $M(a)$ is a matrix of size $\ff(t(a)) \times \ff(s(a))$.

The group
$$
\GL_\ff := \prod_{i \in \GG_0} \GL_{\ff(i)}(\C)
$$
acts from the left by conjugation on $\md(\C\GG,\ff)$, \textit{i.e.}~for 
$M = (M(a))_{a \in \GG_1} \in \md(\C\GG,\ff)$ and
$g = (g(i))_{i \in \GG_0} \in \GL_\ff$ we have
$$
(g.M)(a) = g(t(a))M(a)g(s(a))^{-1} 
$$
for all $a \in \GG_1$.
The orbits of $\GL_\ff$ on $\md(\C\GG,\ff)$ correspond to the set of 
isomorphism classes of $\C\GG$-modules with dimension vector $\ff$.
Given a path $p = a_l \cdots a_2a_1$ in $\GG$ 
(\textit{i.e.} $a_1,\ldots,a_l$ are arrows with $s(a_{i+1}) = t(a_i)$ for 
$1 \le i \le l-1$) we define 
$$
M(p) := M(a_l)\cdots M(a_2)M(a_1)
$$ 
for any $M \in \md(\C\GG,\ff)$. 
More generally, for any element  
$\rho \in e_i\C\GG e_j$ we have 
$M(\rho) \in \C^{\ff(i) \times \ff(j)}$, since
$\rho$ is a linear combination of paths from $j$ to $i$. 
(For $k \in \GG_0$ we denote the associated path of length 0 by $e_k$.) 
Set $s(\rho) := j$ and $t(\rho) := i$.
If $I \subset \C\GG$
is a finitely generated ideal contained in the ideal generated by all
paths of length 2, we may assume that it is generated by elements
$\rho_1,\ldots,\rho_q$ with $\rho_k \in e_{t_k}\C\GG e_{s_k}$ for 
certain $s_k,t_k \in \GG_0$ where $1 \le k \le q$. 
Let $A := \C\GG/I$.
We consider the affine $\GL_\ff$-variety
$$
\md(A,\ff) := \{ M \in \md(\C\GG,\ff) \mid M(\rho_k) = 0 
\text{ for } 1 \le k \le q \}.
$$
Again, the $\GL_\ff$-orbits correspond to
the isomorphism classes of $A$-modules with dimension 
vector $\ff$.

Given $M \in \md(A,\ff)$ and $M' \in \md(A,\ff')$ we identify any homomorphism
$\vph \in \Hom_A(M,M')$ with a family of matrices
$$
(\vph(k))_{k \in \GG_0} \in \prod_{k \in \GG_0} \C^{\ff'(k) \times \ff(k)}
$$
such that
$$
\vph(t(b))M(b) = M'(b)\vph(s(b)) 
$$
for all $b \in \GG_1$.
In other words, the diagram
$$
\xymatrix{
\C^{\ff(s(b))} \ar[d]_{M(b)}\ar[r]^{\vph(s(b))} &
\C^{\ff'(s(b))} \ar[d]^{M'(b)} \\ 
\C^{\ff(t(b))}  \ar[r]^{\vph(t(b))} &\C^{\ff'(t(b))}
}
$$
commutes for all $b \in \GG_1$.

\subsection{A stratification of $\LL_d^w$}
Recall that for $X \in \nil(\LL)$ we have $X \in \CC_w$ if and only if 
there is a (unique) filtration 
$$
0 = X_0 \subseteq X_1 \subseteq \cdots \subseteq X_r = X
$$ 
by submodules such that $X_k/X_{k-1} \cong M_k^{a_k}$ for 
some $a_k \ge 0$ for all $1 \le k \le r$, see Proposition~\ref{M4}.
In this case, we have
$$
\dimv_{B_\ii}(\Hom_\LL(V_\ii,X)) = 
\sum_{k=1}^r a_k\, \dimv_{B_\ii}(\Delta_k),
$$
\textit{i.e.} $\ba := (a_1,\ldots,a_r)$ is the $\Delta$-dimension vector of 
$\Hom_\LL(V_\ii,X)$.
Thus, with
$$
\mu(\ba) := \sum_{k=1}^r a_k\, \dimv_\LL(M_k)
$$
we may consider
\begin{multline*}
\LL^\ba := \{ X \in \LL_{\mu(\ba)} \mid X \text{ has a filtration }
0 = X_0 \subseteq X_1 \subseteq \cdots \subseteq X_r = X\\ 
\text{ with } X_k/X_{k-1}\cong M_k^{a_k},\, 1 \le k \le r \}. 
\end{multline*}
In other words, 
$\LL^\ba = \{ X \in \LL_{\mu(\ba)} \mid X \in \CC_{M_\ii,\ba} \}$.
Define
$$
\LL^w_\bd := \{ X \in \LL_\bd \mid X \in \CC_w \}.
$$
We get a finite decomposition
$$
\LL^w_\bd \;\; =
\bigcup_{\ba \in \N^r,\;\mu(\ba) = \bd} \LL^\ba
$$
into disjoint subsets.

\begin{Lem}\label{stratlem}
$\LL^\ba$ is an irreducible constructible subset of $\LL_{\mu(\ba)}$.
\end{Lem}

\begin{proof}
We know from Proposition~\ref{M6} that $X \in \LL^\ba$ 
if and only if there exists a short exact sequence
$$
0 \ra \bigoplus_{k=1}^r V_{k^-}^{a_k} \ra 
\bigoplus_{k=1}^r V_k^{a_k} \ra X \ra 0
$$
with $V_{k^-} = 0$ if $k^- = 0$. 
Now the result follows from
\cite[Section~2.1]{Bo}.
\end{proof}

\begin{Rem}\label{rem:ssec-strat}
It is not hard to see that for $X \in \nil(\LL)$ the following
are equivalent
\begin{itemize}

\item[(i)]
$X \in \CC_w$;

\item[(ii)] 
$\Hom_\LL(D(J_w),X) = 0 = \Ext^1_\LL(D(\LL/J_w),X)$.

\end{itemize}
Here $J_i$ is by definition the ideal of $\LL$ which is as a $\C$-vector space
generated by all paths $p$ in $\overline{Q}$ with $p \not= e_i$, and we
set $J_w := J_{i_r}\cdots J_{i_1}$.
It follows that $\LL_\bd^w$ is an open subset in $\LL_\bd$
and it follows that $\LL^\ba$ is a locally closed subset of $\LL_{\mu(\ba)}$.
However, we will not need this fact.
\end{Rem}

\subsection{Review of Bongartz's bundle  construction}\label{reviewbong}
Following
Bongartz~\cite[Section~4]{Bo},
we apply the above definitions and conventions in order to relate
the varieties $\LL^\ba$ and $\md(B_\ii,\ff)$.
Assume that
$$
\ff = \sum_{k=1}^r a_k\, \dimv_{B_\ii}(\Delta_k).
$$ 
Recall that this implies
$\dm \Hom_\LL(V_k,X) = \ff(k)$ for all $X \in \LL^\ba$. 
It follows, that 
$$
\left\{ (X,\vph) \mid X \in \LL^\ba \text{ and } \vph \in \Hom_\LL(V_k,X) 
\right\}
$$
is a (usually non-trivial) algebraic vector bundle of rank $\ff(k)$ over
$\LL^\ba$. 
Thus, setting
$$
I(\ff) := \left\{ (i,j) \in \N_1^2 \mid 1 \le i \le r 
\text{ and } 1 \le j \le \ff(i) \right\}
$$
we consider
\begin{multline*}
H^\ba := \{ (X,(\vph^{(i)}_j)_{(i,j) \in I(\ff)}) \mid X \in \LL^\ba
\text{ and }\\
(\vph^{(k)}_1,\ldots,\vph^{(k)}_{\ff(k)}) \text{ is a basis of }
\Hom_\LL(V_k,X),\, 1 \le k \le r \}
\end{multline*}
equipped with a left $\GL_\bd$-action given by
$$
g.(X,(\vph^{(k)}_j)_{(k,j)\in I(\ff)}) = 
(g.X, (g\circ\vph^{(k)}_j)_{(k,j)\in I(\ff)}),
$$
and with a right $\GL_\ff$-action given by
$$
(X,(\vph^{(k)}_j)_{(k,j)\in I(\ff)}).h=
(X,(\sum_{t=1}^{\ff(k)}\vph^{(k)}_t h_{t,j}(k))_{(k,j)\in I(\ff)}).
$$
Here $h_{t,j}(k)$ denotes the entry in row $t$ and column $j$ of the matrix
$h_k$.
Clearly, the map
$$
\pi_1\df H^\ba \ra \LL^\ba
$$
defined by
$$
(X,(\vph^{(k)}_j)_{(k,j)\in I(\ff)}) \mapsto X
$$
is a $\GL_\bd$-equivariant $\GL_\ff$-principal bundle.

In order to define a map $\pi_2\df H^\ba \ra \md(B_\ii,\ff)$ we write
$B_\ii = \C\GG/I$ for an admissible ideal I, and we identify the vertices
$\GG_0 = \{ 1,2,\ldots,r \}$ with the summands $V_1,\ldots,V_r$ of $V_\ii$.
Recall that $\GG \equiv \GG_\ii$. 
Thus we may think of each arrow $b\df i \ra j$ in $\GG_1$ as a certain element
$b \in \Hom_\LL(V_j,V_i)$. 
With these identifications
$$
\pi_2(X,(\vph^{(k)}_j)_{(k,j)\in I(\ff)}) = (M(b))_{b \in \GG_1}
$$
is determined  by
$$
\vph_j^{(s(b))} \circ b = \sum_{u=1}^{\ff(t(b))} \vph_u^{(t(b))}M_{u,j}(b).
$$
Here $M_{u,j}(b)$ denotes the entry in row $u$ and column $j$ of the matrix
$M(b)$.
It is easy to verify that $\pi_2$ is a $\GL_\bd$-invariant 
$\GL_\ff$-equivariant morphism,
if we view $\md(B_\ii,\ff)$ with the {\it right} $\GL_\ff$-action induced from
the usual left action via the anti-automorphism $h \mapsto h^{-1}$ of 
$\GL_\ff$. 
Moreover, by construction
$$
\pi_2(X,(\vph^{(k)}_j)_{(k,j)\in I(\ff)}) \cong \Hom_\LL(V_\ii,X)
$$
as a $B_\ii$-module. 
Thus, in fact $\Ima(\pi_2) = \F(\Delta,\ff)$, where
$\F(\Delta,\ff)$ is the subset of $\md(B_\ii,\ff)$ consisting
of the $\Delta$-filtered $B_\ii$-modules with dimension vector $\ff$.
It is shown in \cite[Corollary~1.5]{cbsr01} that $\F(\Delta,\ff)$ is 
open in $\md(B_\ii,\ff)$. Since $\pi_1$ is a  $\GL_\ff$-principal
bundle it follows from Lemma~\ref{stratlem} that 
$\F(\Delta,\ff)=\pi_2(\pi_1^{-1}(\LL^\ba))$ is also irreducible.
In particular, $\overline{\F(\Delta,\ff)}$ is an irreducible component of
$\md(B_\ii,\ff)$.

Finally, for $\pi_2(X,(\vph^{(k)}_j)_{(k,j)\in I(\ff)}) = M$ we have
\begin{align*}
\dm \GL_\bd . X 
&= \dm \GL_\bd - \dm \End_\LL(X)\\
\dm \pi_1^{-1}(\GL_\bd . X) 
&= \dm \GL_\bd -\dm \End_\LL(X) + \dm \GL_\ff\\
\dm(\GL_\ff.M)   &= \dm \GL_\ff - \dm \End_\LL(X).
\end{align*}
The last equation holds, since the functor 
$F_\ii\df \CC_w \ra \F(\Delta)$ which maps $X$ to $\Hom_\LL(V_\ii,X)$
is an equivalence of additive categories. 
By the same token
$\pi_2^{-1}(M.\GL_\ff) = \pi_1^{-1}(\GL_\bd.X)$. 
We conclude $\dm \pi_2^{-1}(M) = \dm \GL_\bd$.
Thus we proved the following:

\begin{Lem}
For $\ba \in \N^r$ and 
$$
\ff = \sum_{k=1}^r a_k\, \dimv_{B_\ii}(\Delta_k)
$$ 
there exists a 
variety $H^\ba$ with a $\GL_\bd$-$\GL_\ff$-action
together with two surjective morphisms
$$
\xymatrix@!@-1.8pc{ 
& \ar[dl]_{\pi_1} H^\ba \ar[dr]^{\pi_2}\\
\LL^\ba && \F(\Delta,\ff)
}
$$
such that $\pi_1$ is a $\GL_\bd$-equivariant $\GL_\ff$-principal bundle,
and $\pi_2$ is a $\GL_\ff$-equivariant and $\GL_\bd$-invariant morphism.
Moreover, $\dm \pi_2^{-1}(M) = \dm \GL_\bd$ for all $M \in \F(\Delta,\ff)$.
\end{Lem}

Since $\CC_w = \Gen(V_\ii)$, it is easy to see that
for $g \in \GL_\bd$ and $h \in H^\ba$ with $g.h=h$ we have $g = 1_{\GL_\bd}$.

\begin{Rem}
It seems plausible that with a dual bundle construction, 
as in \cite[4.3]{Bo}, one can show that $\pi_2$ is a $\GL_\bd$-principal
bundle. 
\end{Rem}

\subsection{Parametrization of components}\label{paracomp}

\begin{Lem}\label{dimdelta}
For $\ba = (a_1,\ldots,a_r)$, $\bd = \mu(\ba)$ and 
$\ff = \sum_{k=1}^r a_k\, \dimv_{B_\ii}(\Delta_k)$
we have
$$
\dm \F(\Delta,\ff) = \dm \GL_\ff - \bil{\bd,\bd}_Q.
$$
\end{Lem}

\begin{proof}
For any $N \in \F(\Delta,\ff)$ we have $\pdim_{B_\ii}(N) \le 1$, thus
$\Ext^2_{B_\ii}(N,N) = 0$, which implies that $N$ is a smooth 
point~\cite[3.7]{geiss96} of the {\it scheme} $\md(B_\ii,\ff)$. 
Recall that $\mu(\ba) = \sum_{k=1}^r a_k\, \dimv(M_k)$.
Now Voigt's Lemma~\cite[1.3]{gabr74} and our Lemma~\ref{M8} allow 
the calculation
\begin{align*}
\dm \F(\Delta,\ff) 
&= \dm \GL_\ff.N +\dm \Ext^1_{B_\ii}(N,N)\\
&= \dm \GL_\ff - \bil{\ff,\ff}_{B_\ii}\\
&= \dm \GL_\ff - \left(\sum_{k=1}^r a_k^2 
\bil{\Delta_k,\Delta_k}_{B_\ii}
+ \sum_{1 \le s<k \le r}a_ka_s 
\bil{\Delta_k,\Delta_s}_{B_\ii}\right)\\
&= \dm \GL_\ff - \left(\sum_{k=1}^r a_k^2 \bil{M_k,M_k}_Q
+ \sum_{1 \le s<k \le r} a_ka_s (M_k,M_s)_Q\right)\\
&= \dm \GL_\ff - \bil{\bd,\bd}_Q.
\end{align*}          
For the fourth equality we used Lemma~\ref{M8} and the fact that
$$
\bil{\Delta_k,\Delta_k}_{B_\ii} = 1 = \bil{M_k,M_k}_Q.
$$
This finishes the proof.
\end{proof}

\begin{Prop}\label{dimdelta2}
The (Zariski-) closure $Z^\ba$ of $\LL^\ba$ 
is an irreducible component of $\LL_{\mu(\ba)}$.
In particular, $Z^\ba$ is the unique irreducible component of $\LL_{\mu(\ba)}$
which contains a dense open subset which belongs to $\LL^\ba$.
\end{Prop}

\begin{proof}
We know from Lemma~\ref{stratlem} that $\LL^\ba$ 
is an irreducible constructible subset of $\LL_{\mu(\ba)}$.
Thus, $Z^\ba$ is an irreducible subvariety of $\LL_{\mu(\ba)}$.
Since $\LL_{\mu(\ba)}$ is equi-dimensional, it remains to show that 
$$
\dm \LL^\ba = \dm \GL_\bd - \bil{\bd,\bd}_Q = \dm \LL_\bd.
$$
Recall that $\bd = \mu(\ba)$ and
$$
\ff = \sum_{k=1}^r a_k\, \dimv_{B_\ii}(\Delta_k).
$$ 
As before, $\F(\Delta,\ff)$ denotes the irreducible open subset of 
$\Delta$-good modules in the affine $\GL_\ff$-variety $\md(B_\ii,\ff)$ of 
$B_\ii$-modules with dimension vector $\ff$.
By Lemma~\ref{dimdelta} we know that
$$
\dm \F(\Delta,\ff) = \dm \GL_\ff - \bil{\bd,\bd}_Q.
$$
In Section~\ref{reviewbong} we constructed a
$\GL_\bd$-$\GL_\ff$-variety $H^\ba$ together with
surjective morphisms
$$
\xymatrix@!@-1.8pc{
& \ar[dl]_{\pi_1} H^\ba \ar[dr]^{\pi_2} \\
\LL^\ba && \F(\Delta,\ff)
}
$$
with $\pi_1$ a $\GL_\bd$-equivariant $\GL_\ff$-principal bundle, and
$\pi_2$ a $\GL_\ff$-equivariant morphism with all fibers having the
same dimension as $\GL_\bd$.
Our claim about the dimension of $\LL^\ba$ follows.
\end{proof}

Let $M = M_1^{a_1} \oplus \cdots \oplus M_r^{a_r}$ for some
$\ba = (a_1,\ldots,a_r) \in \N^r$.
We just proved that
$Z^\ba$ is an irreducible component of $\LL_{\mu(\ba)}$.
Let us denote the corresponding dual semicanonical basis
vector $\rho_{Z^\ba}$ by $s_M$.
Thus there is a dense open subset $U^\ba \subseteq Z^\ba$ such that
$s_M = \delta_X$ for all $X \in U^\ba$.


\section{A dual PBW-basis and a dual semicanonical basis for $\cA(\CC_w)$}
\label{PBWsection}


In this section we prove Theorem~\ref{main7} and Theorem~\ref{basesthm}.
We also deduce from these results the existence of semicanonical bases
for the cluster algebras $\widetilde{\RR}(\CC_w,T)$ and
$\underline{\RR}(\CC_w,T)$ obtained by inverting and specializing
coefficients, respectively.

\subsection{Proof of Theorem~\ref{main7}}\label{clustercattoclusteralg}
By the definition of the cluster algebra $\cA(\CC_w,T)$, 
its initial seed
is $(\yy, B(T)^\circ)$ where $\yy = (y_1,\ldots,y_r)$.
In particular, $\cA(\CC_w,T)$ is a subalgebra of
$\F := \C(y_1,\ldots,y_r)$.
Since $T$ is rigid, by Theorem~\ref{2.6i} and \cite[Theorem 1.1]{GLSSemi1}
every monomial in the $\delta_{T_k}$ belongs to
the dual semicanonical basis $\cS^*$, hence the $\delta_{T_k}$ 
are algebraically
independent, and $(\delta_{T_1}, \ldots , \delta_{T_r})$ is a transcendence
basis of the subfield $\G$ it generates inside the fraction field
of $U(\n)^*_{\text{gr}}$.
Let $\iota\df \F\to\G$ be the field isomorphism defined by
$\iota(y_k) = \delta_{T_k}$ where $1\le k \le r$.
Combining Theorems~\ref{2.6ii} and \ref{main4} we see that
the cluster variable $z$ of $\cA(\CC_w,T)$ obtained from the
initial seed $(\yy, B(T)^\circ)$ through a sequence of seed mutations
in successive directions $k_1,\ldots,k_s$ will be mapped by $\iota$
to $\delta_X$, where $X \in \CC_w$ is the indecomposable rigid module obtained
by the same sequence of mutations of rigid modules.
It follows that $\iota$ restricts to an isomorphism from
$\cA(\CC_w,T)$ to $\RR(\CC_w,T)$.
This isomorphism is completely determined by the images of the 
elements $y_k$, hence the unicity. 
The cluster monomials are mapped to elements $\delta_R$ where 
$R$ is a (not necessarily $\CC_w$-maximal or basic) 
rigid module in $\CC_w$, hence an element of $\cS^*$.
More precisely, the cluster monomials in $\RR(\CC_w,T)$ are
the elements $\delta_R$, where $R$ runs through the set of all 
$T$-reachable modules (see Section~\ref{reach} for the definition
of $T$-reachable).
This finishes the proof of Theorem~\ref{main7}.

\subsection{Proof of Theorem~\ref{basesthm}}
Let $M_\ii = M_1 \oplus \cdots \oplus M_r$ be as before.
For $1 \le k \le r$ we proved that
$\dimv(M_k) = \beta_\ii(k)$.
Set $\beta(k) := \beta_\ii(k)$.

We have
$$
\C[\delta_{M_1},\ldots, \delta_{M_r}] \subseteq 
\RR(\CC_w,V_\ii) \subseteq
\Span_\C\ebrace{\delta_X \mid X \in \CC_w},
$$
where the first inclusion follows from the observation that each of the
$\LL$-modules $M_k$ for $1 \le k \le r$ is the direct summand of a
$\CC_w$-maximal rigid module on the mutation path from $V_\ii$ to $T_\ii$, 
see Section~\ref{section17}. 
The second inclusion follows from the observation
that $\Span_\C\ebrace{\delta_X \mid X \in \CC_w}$ is an algebra. 
This follows
from the fact that $\CC_w$ is an additive category together with
Theorem~\ref{2.6i}.

For each $M \in \add(M_\ii)$ 
we constructed a dual semicanonical basis vector
$s_M$, see the explanation at the end of Section~\ref{paracomp}.
If $M = M_k$ is an indecomposable direct summand of $M_\ii$,
then $s_M = \delta_{M_k}$.
(For every rigid $\LL$-module $R \in \nil(\LL)$, the function
$\delta_R$ belongs to the dual semicanonical basis.
The modules $M_k$ are rigid by Corollary~\ref{M2}.)

The following theorem is a slightly more explicit statement
of Theorem~\ref{basesthm}:

\begin{Thm}\label{proofmain1}
Let $w$ be a Weyl group element, and let
$\ii = (i_r,\ldots,i_1)$ be a reduced expression of $w$.
Then the following hold:
\begin{itemize}

\item[(i)]
We have
$$
\RR(\CC_w,V_\ii) = \C[\delta_{M_1},\ldots,\delta_{M_r}]
= \Span_\C\ebrace{\delta_X \mid X \in \CC_w};
$$

\item[(ii)]
The set
$$
\left\{ \delta_M \mid M \in \add(M_\ii) \right\}
$$
is a $\C$-basis of $\RR(\CC_w,V_\ii)$;

\item[(iii)]
The subset
$$
\cS_w^* := \left\{ s_M \mid M \in \add(M_\ii) \right\}
$$
of the dual semicanonical basis
is a $\C$-basis of $\RR(\CC_w,V_\ii)$, and all
cluster monomials of $\RR(\CC_w,V_\ii)$ belong
to $\cS_w^*$.

\end{itemize}
\end{Thm}

The basis $\left\{ \delta_M \mid M \in \add(M_\ii) \right\}$ 
will be called 
{\it dual PBW-basis} 
of
$\RR(\CC_w,V_\ii)$,
and $\cS_w^*$ the 
{\it dual semicanonical basis} 
of
$\RR(\CC_w,V_\ii)$.
The proof of this theorem will be given after a series of lemmas.

Let
$$
\n = \bigoplus_{d \in \Delta^+} \n_d
$$
be the root space decomposition of $\n$.
We consider $\n$ as a subspace of the universal enveloping
algebra $U(\n)$.
Since we identify $U(\n)$ and $\MM$,
we can think of an element $f$ in $\n_d$ as a constructible function
$f\df \LL_d \to \C$ in $\MM_d$.

\begin{Lem}\label{prooflemma2}
Let $f \in \n_d$.
If 
$d \not\in \{ \beta(k) \mid 1 \le k \le r \}$,
then 
$$
f(X) = 0 \text{ for all } X \in \CC_w.
$$
\end{Lem}

\begin{proof}
Let $X \in \CC_w$, and let
$f \in \n_d$ with $f(X) \not= 0$.
In particular, $f \in \MM_d$, and we have
$d = \dimv(X) \in \Delta^+$.
We know that $X \in \CC_{M_\ii,\aaa}$ for some
$\aaa = (a_1,\ldots,a_r)$.
Thus 
$$
\dimv(X) = \sum_{k=1}^r a_k\, \dimv(M_k).
$$
By  Lemma~\ref{bracketclosed1},
$\Delta_w^+ = \{ \beta(k) \mid 1 \le k \le r \}$ 
is a bracket closed subset of $\Delta^+$.
Thus $d = \beta(s)$ for some $1 \le s \le r$.
This finishes the proof.
\end{proof}

As before, let
$\ii = (i_r,\ldots,i_1)$ be a reduced expression of a Weyl group element 
$w$, and let 
$$
\cP = \left\{ p_{\bf m} \mid {\bf m} \in \N^{(J)}_{\phantom{0}} \right\}
$$ 
be an $\ii$-compatible PBW-basis of $U(\n)$, see 
Sections~\ref{universal} and \ref{nw}.

\begin{Lem}\label{prooflemma3}
Let $p_\mm \in \cP$ where $\mm = (m_j)_{j\in J}$.
If $m_j > 0$ for some $j > r$, then
$$
p_\mm(X) = 0 \text{ for all } X \in \CC_w.
$$
Equivalently, $\delta_X(p_\mm) = 0$ for all $X \in \CC_w$.
\end{Lem}

\begin{proof}
We regard $p_\mm$ as an element of $\MM$, hence as a convolution
product
\[
p_\mm = p_1^{(m_1)}\star p_2^{(m_2)} \star \cdots \star p_s^{(m_s)}.
\] 
Let us assume that $s > r$ and $m_s > 0$.
It follows that 
$p_\mm = p \star p_s$ where 
$$
p := \frac{1}{m_s} \left(p_1^{(m_1)} 
\star p_{2}^{(m_{2})} \star \cdots \star 
p_{s-1}^{(m_{s-1})} \star p_s^{(m_s-1)}\right).
$$
Now let $X \in \CC_w$.
Then 
$$
p_{\mathbf m}(X) = (p \star p_s)(X) = 
\sum_{m \in \C} m \, \chi_{\rm c}(\{ U \subseteq X \mid p(U)p_s(X/U) = m \}).
$$
Since $\CC_w$ is closed under factor modules, we get $X/U \in \CC_w$
for all submodules $U$ of $X$.
Now Lemma~\ref{prooflemma2} yields $p_s(X/U) = 0$ for all such $U$.
Thus we proved that $p_\mm(X) = 0$ for all $X \in \CC_w$.
\end{proof}

Recall from Section~\ref{nw} that
$$
\cP_\ii^* = \left\{ (p_1^*)^{m_1} \cdots (p_r^*)^{m_r} \mid 
m_k \ge 0 \text{ for all } 1 \le k \le r \right\}
$$
is a subset of the dual PBW-basis $\cP^*$ of $U(\n)_{\rm gr}^*$.
The following lemma is of central importance:

\begin{Lem}\label{prooflemma6}
For $1 \le k \le r$ we have $p_k^* = \delta_{M_k}$ (up to rescaling of
$p_k$).
\end{Lem}

\begin{proof}
For each $1 \le k \le r$
there exists some $\mm = (m_i)_{i \ge 1}$ such that
$p_\mm(M_k) \not = 0$, since
$\delta_{M_k} \in \MM^* \equiv U(\n)_{\rm gr}^*$ 
is a linear combination of elements in $\cP^*$.
Let $s$ be the natural number with $m_s \ge 1$, but
$m_j = 0$ for all $j > s$.

By Lemma~\ref{prooflemma3},
if $s > r$, then $p_\mm(X) = 0$ for all modules $X \in \CC_w$,
a contradiction.
Thus, we know that $s \le r$.
We even know that $s \le k$, since $M_k$ is an object of the subcategory
$\CC_u$ of $\CC_w$, where $u = s_{i_k} \cdots s_{i_2}s_{i_1}$.

If $s = k$, then for dimension reasons $m_1 = \cdots = m_{k-1} = 0$ and
$m_k = 1$.
So we get $p_\mm = p_k$.

Finally, assume $s < k$.
Since $p_\mm(M_k) \not= 0$, we know that
$M_k$ has a filtration
$$
0 = U_{1,0} \subseteq U_{1,1} \subseteq \cdots \subseteq U_{1,m_1} =
U_{2,0} \subseteq U_{2,1} \subseteq \cdots \subseteq U_{2,m_2} =
U_{3,0} \subseteq \cdots \subseteq U_{s,m_s} = M_k
$$
such that $p_i(U_{i,j}/U_{i,j-1}) \not= 0$ for all $1 \le i \le s$ and
$1 \le j \le m_i$.
But we know that $p_i$ lies in $\MM_{\beta(i)}$.
In other words, we have
$\dimv(U_{i,j}/U_{i,j-1}) = \beta(i)$.
This implies that $\beta(k)$, the dimension vector of $M_k$, is 
a positive integer linear combination of $\beta(i)$'s with
$i < k$.
More precisely, 
$$
\beta(k) = m_1\beta(1) + \cdots + m_s\beta(s).
$$
But $\beta(1),\ldots,\beta(s)$ belong to the bracket closed set
$\Delta_v^+$ where $v := s_{i_s} \cdots s_{i_2}s_{i_1}$.
Thus $\beta(k)$ is also in $\Delta_v^+$, which is a contradiction,
since $s < k$.

Summarizing, we proved that $p_\mm(M_k) \not= 0$ if and only
if $p_\mm = p_k$.
Now we can rescale our PBW-basis elements
$p_k$, and we obtain without loss of generality that $p_k(M_k) = 1$.
Thus we proved that
$$
\delta_{M_k}(p_\mm) := p_\mm(M_k) =
\begin{cases}
1 & \text{if $p_\mm = p_k$},\\
0 & \text{otherwise}.
\end{cases}
$$
In other words, $\delta_{M_k} = p_k^*$.
\end{proof}

\begin{Cor}\label{prooflemma7}
Under the identification $U(\n)_{\rm gr}^* \equiv \MM^*$ we have 
$$
\cP_\ii^* = \left\{ \delta_M \mid M \in \add(M_\ii) \right\}.
$$
\end{Cor}

\begin{proof}
By definition
$$
\cP_\ii^* = \left\{ (p_1^*)^{m_1} \cdots (p_r^*)^{m_r} \mid 
m_k \ge 0 \text{ for all } 1 \le k \le r \right\} \subseteq \cP^*.
$$
This implies the result, since
$p_k^* = \delta_{M_k}$ and 
$\delta_X \cdot \delta_Y = \delta_{X \oplus Y}$
for all nilpotent $\LL$-modules $X$ and $Y$.
\end{proof}

\begin{proof}[Proof of Theorem~\ref{proofmain1}]
Let $X \in \CC_w$.
By Lemma~\ref{prooflemma3} and Corollary~\ref{prooflemma7}, 
$\delta_X$ is a linear combination of
dual PBW-basis vectors of the form $\delta_M$ with
$M \in \add(M_\ii)$.
Hence
$\delta_X \in \C[\delta_{M_1},\ldots,\delta_{M_r}]$,
and
$$
\Span_\C\ebrace{\delta_X \mid X \in \CC_w}
\subseteq
\C[\delta_{M_1},\ldots,\delta_{M_r}]
\subseteq
\RR(\CC_w,V_\ii). 
$$
Using the known reverse inclusions we
get (i) and (ii) of Theorem~\ref{proofmain1}.

Next, let $M = M_1^{a_1} \oplus \cdots \oplus M_r^{a_r}$
be a module in $\add(M_\ii)$.
Set $\ba := (a_1,\ldots,a_r)$.
Then $s_M = \delta_X$ for some module $X$
in $Z^\ba$. 
In particular, $X$ is contained in $\CC_w$.
Thus, by what we proved up to now we get
$$
s_M = \delta_X \in \RR(\CC_w,V_\ii).
$$
For dimension reasons this implies that
$$
\cS_w^* := 
\left\{ s_M \mid M \in \add(M_\ii) \right\}
= \cS^* \cap \RR(\CC_w,V_\ii)
$$
is a $\C$-basis of $\RR(\CC_w,V_\ii)$.
By what we proved before, the set of
cluster monomials of $\RR(\CC_w,V_\ii)$ are a subset
of $\cS_w^*$.
This finishes the proof of Theorem~\ref{proofmain1}.
\end{proof}

By Theorem~\ref{proofmain1}, we know that 
$$
\RR(\CC_w,V_\ii) = \C[p_1^*,\ldots,p_r^*].
$$
Thus Proposition~\ref{prop17} yields the following result:

\begin{Prop}\label{caN(w)}
Under the identification $U(\n)^*_{\rm gr}\equiv \C[N]$
the cluster algebra $\RR(\CC_w,V_\ii)$ gets identified to the
ring of invariants $\C[N]^{N'(w)}$, which is isomorphic to
$\C[N(w)]$.
\end{Prop}

\begin{Cor}
Let $\ii = (i_r,\ldots,i_1)$ 
be a reduced expression of $w$.
For $X \in \CC_w$, the function
$\varphi_X \in \C[N]$ is determined by its values
on 
$
\{ x_\ii(t) \mid t = (t_r,\ldots,t_1) \in (\C^*)^r \}
$
where $x_\ii(t) := x_{i_r}(t_r) \cdots x_{i_2}(t_2)x_{i_1}(t_1)$.
\end{Cor}

\begin{proof}
Let $\varphi,\psi \in \C[N]^{N'(w)}$.
Then $\varphi = \psi$ if and only if
$\varphi(x_\ii(t)) = \psi(x_\ii(t))$ 
for all $t \in (\C^*)^r$:
Recall that
each $x \in N$ can be written as $x = yy'$ for 
a unique $(y,y') \in N(w) \times N'(w)$.
For $x \in N^w$ we have $\pi_w(x) = y$.
Furthermore, 
the image of $\pi_w$ is dense in $N(w)$, see Proposition~\ref{coordNw}.
It is well known that the set $\{ x_\ii(t) \mid t \in (\C^*)^r \}$
contains a dense open subset of $N^w$.
For $x = x_\ii(t)$ we get
$$
\varphi(\pi_w(x)) = \varphi(y) = \varphi(yy') = \varphi(x).
$$
For the second equality we used that $\varphi$ is $N'(w)$-invariant.
Since $\varphi$ is a regular map, its values on the whole of
$N(w)$ are already determined by its values on
$\pi_w(x_\ii(t))$, $t \in (\C^*)^r$.
\end{proof}

\subsection{Proof of Theorem~\ref{wident}}
By Proposition~\ref{coordNw}, we know that $\C[N^w]$ is
the localization of the ring $\C[N(w)]$ with respect to the
minors $D_{\vpi_i,w^{-1}(\vpi_i)}$. 
By Proposition~\ref{caN(w)}, $\C[N(w)]$ is equal to the cluster
algebra $\RR(\CC_w,V_\ii)$. 
By Proposition~\ref{minors}, the minors $D_{\vpi_i,w^{-1}(\vpi_i)}$
coincide with the functions $\vph_X$ where $X$ runs through the 
set of indecomposable $\CC_w$-projective-injectives.
In other words, the $D_{\vpi_i,w^{-1}(\vpi_i)}$ coincide with
the generators of the coefficient ring of $\RR(\CC_w,V_\ii)$.
Hence $\C[N^w]$ is equal to the cluster algebra
$\widetilde{\RR}(\CC_w,V_\ii)$.

\subsection{Example}
Let us discuss an example of base change between 
$\cP_\ii^*$ and $\cS_w^*$.
Let $Q$ be a quiver with underlying graph
$\xymatrix@-0.5pc{1 \ar@{-}[r] & 2 \ar@{-}[r] & 3}$ and
let $\ii := (i_6,\ldots,i_1) := (2,3,1,2,3,1)$, 
which is a reduced expression of the
Weyl group element $w := s_2s_3s_1s_2s_3s_1$.
As before, let $V_\ii = V_1 \oplus \cdots \oplus V_6$ and
$M_\ii = M_1 \oplus \cdots \oplus M_6$, where as always $M_k = M[k,k]$.
The $\LL$-modules $V_k$ are the following:
\begin{align*}
V_1 &= M_1 = {\bsm 1 \esm} & 
V_2 &= M_2 = {\bsm 3 \esm} &
V_3 &= M_3 = {\bsm 1&&3\\&2 \esm} \\
V_4 &= M[4,1] = {\bsm &&3\\&2\\1\esm} & 
V_5 &= M[5,2] = {\bsm 1\\&2\\&&3\esm} &
V_6 &= M[6,3] = {\bsm &2\\1&&3\\&2\esm}.
\end{align*}
The initial cluster of our cluster algebra $\RR(\CC_w,V_\ii)$ looks
as follows:
$$
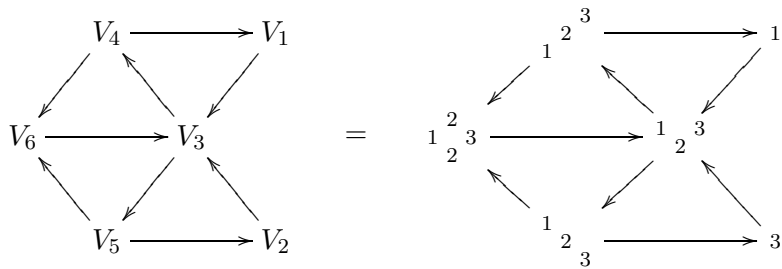
\xymatrix@=0.5cm{
& V_4 \ar[rr]\ar[dl] && V_1 \ar[dl] &&
& {\bsm &&3\\&2\\1\esm} \ar[rr]\ar[dl] && {\bsm 1 \esm} \ar[dl] \\
V_6 \ar[rr] && V_3 \ar[ul]\ar[dl] &&=&
{\bsm &2\\1&&3\\&2\esm} \ar[rr] && {\bsm 1&&3\\&2\esm} \ar[ul]\ar[dl]\\
& V_5 \ar[ul]\ar[rr] && V_2 \ar[ul] &&
& {\bsm 1\\&2\\&&3\esm} \ar[ul]\ar[rr] && {\bsm 3 \esm} \ar[ul]
}
$$
We have
\begin{align*}
M_4 &= {\bsm&3\\2\esm} & 
M_5 &= {\bsm1\\&2 \esm} &
M_6 &= {\bsm 2 \esm}.
\end{align*}
There are only three more indecomposable $\LL$-modules, namely
\begin{align*}
W_1 &= {\bsm2\\&3\esm} & 
W_2 &= {\bsm&2\\1\esm} &
W_3 &= {\bsm&2\\1&&3\esm}.
\end{align*}
Observe that $\Omega(V_k) = W_k$ for $1 \le k \le 3$.

The functions $\delta_{M_k}$ can be computed easily.
Indeed, for all $\jj$ and $k$, the variety
$\F_{\jj,M_k}$ is either empty or
a single point, so $\chi_{\rm c}(\F_{\jj,M_k})$ is either 0 or 1.
Using Theorem~\ref{detthmnew} we get
\begin{align*}
\delta_{V_4} &= \delta_{M_1}\cdot\delta_{M_4} 
- \delta_{M_3},\\
\delta_{V_5} &= \delta_{M_2}\cdot\delta_{M_5}
-\delta_{M_3},\\
\delta_{V_6} &= \delta_{M_3}\cdot\delta_{M_6}
- \delta_{M_4} \cdot \delta_{M_5} .
\end{align*}
Some further exchange relations are
\begin{align*}
\delta_{V_3} \delta_{W_3} &= \delta_{V_4}\cdot\delta_{V_5} 
+ \delta_{V_1}\delta_{V_2}\delta_{V_6},\\
\delta_{V_2}\delta_{W_2} &= \delta_{W_3}+\delta_{V_4},\\
\delta_{V_1}\delta_{W_1} &= \delta_{W_3}+
\delta_{V_5}.
\end{align*}
The cluster variables in $\RR(\CC_w,V_\ii)$
are
$$
\left\{
\delta_{M_k},\, \delta_{V_s},\, \delta_{W_t} \mid
1 \le k \le 6,\, 4 \le s \le 6 \text{ and } 1 \le t \le 3 \right\}.
$$
(Here we consider the three coefficients $\delta_{V_k}$ with $4 \le k \le 6$
also as cluster variables.)
Using the above formulas we get
\begin{align*}
\delta_{W_3} &= \delta_{M_1}\delta_{M_2}\delta_{M_6} 
- \delta_{M_1}\delta_{M_4} - \delta_{M_2}\delta_{M_5} + \delta_{M_3},\\
\delta_{W_2} &= \delta_{M_1}\delta_{M_6} - \delta_{M_5},\\
\delta_{W_1} &= \delta_{M_2}\delta_{M_6} - \delta_{M_4}.
\end{align*}
So we wrote all cluster variables as linear combinations of
dual PBW-basis vectors.

\subsection{Generalities on bases of algebras}\label{generalbase}
A  $\LL$-module $M=\oplus_{k=1}^r M_k^{a_k}$ in $\add(M_\ii)$ 
{\it has gaps} if for each $1\leq j\leq n$ there is some $1\leq k\leq r$ with
$i_k=j$ and $a_k=0$. In other words, $M$ has gaps if and only if $M$ has
no direct summand of the form 
\[
M_\ii(I_{\ii,j}):=
M_{k_\text{max}}\oplus\cdots\oplus M_{k^+_{\text{min}}}\oplus  M_{k_{\text{min}}}
\]
where $i_k=j$.

\begin{Lem}\label{multiproj}
Let
$M = M' \oplus M''$ be in $\add(M_\ii)$ such that
$$
M'' \cong M_\ii(I_{\ii,j}) 
$$
for some $1 \le j \le n$.
Then we have
$s_M = s_{M'} \cdot s_{M''}$.
\end{Lem}

\begin{proof}
We have $s_{M''} = \delta_{I_{\ii,j}}$, and
$I_{\ii,j}$ is $\CC_w$-projective-injective.
The claim follows now easily from
\cite[Theorem 1.1]{GLSSemi1} in combination with
the explanations in \cite[Section 2.6]{GLSSemi1}.
\end{proof}

Let ${\rm B} := \{ b_i \mid i \ge 1 \}$ be a $K$-basis of a
commutative $K$-algebra $A$.
For some fixed $n \ge 1$ let ${\rm C} := \{ b_1,\ldots,b_n \}$.
A basis vector $b \in {\rm B}$ is called 
${\rm C}$-{\it free}
if $b \notin b_i{\rm B}$ for some $b_i \in {\rm C}$.
Assume that the following hold:
\begin{itemize}

\item[(i)]
For all $b_i \in {\rm C}$ we have
$b_i{\rm B} \subseteq {\rm B}$;

\item[(ii)]
If $b_1^{z_1}\cdots b_n^{z_n}b = b_1^{z_1'}\cdots b_n^{z_n'}b'$
for some $z_i,z_i' \ge 0$ and some ${\rm C}$-free elements 
$b,b' \in {\rm B}$, then
$b=b'$ and $z_i = z_i'$ for all $i$.

\end{itemize}
It follows that
${\rm B} = \left\{ b_1^{z_1} \cdots b_n^{z_n}b \mid b \in {\rm B} 
\text{ is {\rm C}-free}, z_i \ge 0 \right\}$.
Define
$$
\underline{A} := A/(b_1-1,\ldots,b_n-1).
$$
For $a \in A$, let $\underline{a}$ be the residue class
of $a$ in $\underline{A}$.
Furthermore, 
let $A_{b_1,\ldots,b_n}$ be the localization of
$A$ at $b_1,\ldots,b_n$.
The following lemma is easy to show:

\begin{Lem}\label{basislemma}
With the notation above, the following hold:
\begin{itemize}

\item[(1)]
The set
$
\underline{\rm B} := \{ \underline{b} \mid b 
\text{ is {\rm C}-free} \}
$
is a $K$-basis of $\underline{A}$;

\item[(2)]
The set
$
{\rm B}_{b_1,\ldots,b_n} := \left\{ b_1^{z_1} \cdots b_n^{z_n}b \mid 
b \in {\rm B}
\text{ is {\rm C}-free}, z_i \in \Z \right\}
$
is a $K$-basis of $A_{b_1,\ldots,b_n}$.

\end{itemize}
\end{Lem}

\subsection{Inverting and specializing coefficients}\label{sectinvertspecial}
One can rewrite the basis $\cS_w^*$ appearing in 
Theorem~\ref{basesthm} as
$$
\cS_w^* = 
\left\{ 
(\delta_{I_{\ii,1}})^{z_1} \cdots (\delta_{I_{\ii,n}})^{z_n}
s_M 
\mid 
M \in \add(M_\ii), M \text{ has gaps, } z_i \ge 0 \right\}.
$$
The next two theorems deal with the situation of invertible
coefficients and specialized coefficients.

\begin{Thm}[Invertible coefficients]\label{basesthm2}
The set
$$
\widetilde\cS_w^* := 
\left\{ 
(\delta_{I_{\ii,1}})^{z_1} \cdots (\delta_{I_{\ii,n}})^{z_n}
s_M 
\mid 
M \in \add(M_\ii), M \text{ has gaps,} z_i \in \Z \right\}
$$ 
is a $\C$-basis of $\widetilde{\RR}(\CC_w,V_\ii)$,
and
$\widetilde\cS_w^*$
contains all cluster monomials of the cluster algebra
$\widetilde{\RR}(\CC_w,V_\ii)$.
\end{Thm}

Next, we specialize all $n$ coefficients $\delta_{I_{\ii,j}}$
of the cluster algebra $\RR(\CC_w,V_\ii)$ to $1$.
We obtain a new cluster algebra 
$\underline{\RR}(\CC_w,V_\ii)$
which does not have any coefficients.
The residue class of $\delta_X \in \RR(\CC_w,V_\ii)$ 
is denoted by $\underline{\delta}_X$.
The residue class of a dual semicanonical basis vector
$s_M$ is denoted by $\underline{s}_M$.

\begin{Thm}[No coefficients]\label{basesthm3}
The set
$$
\underline\cS_w^* := \left\{ \underline{s}_M \mid 
M \in \add(M_\ii), M \text{ has  gaps} \right\}
$$ 
is a $\C$-basis of $\underline{\RR}(\CC_w,V_\ii)$,
and
$\underline\cS_w^*$
contains all cluster monomials of the cluster algebra
$\underline{\RR}(\CC_w,V_\ii)$.
\end{Thm}

\begin{proof}[Proof of Theorem~\ref{basesthm2} and Theorem~\ref{basesthm3}]
Let ${\rm B} := \{ b_i \mid i \ge 1 \} := \cS_w^*$
be the dual semicanonical basis of $\RR(\CC_w,V_\ii)$.
We can label the $b_i$ such that
$$
\{ b_1,\ldots,b_n \} = \left\{ \delta_{I_{\ii,1}},\ldots,
\delta_{I_{\ii,n}} \right\}.
$$
Using Lemma~\ref{multiproj} it is easy to check that
the elements $b_i$ satisfy the properties (i) and (ii) mentioned
in Section~\ref{generalbase}.
Then apply Lemma~\ref{basislemma}.
\end{proof}


\section{Acyclic cluster algebras}\label{acycliccase}


In this section we will study the case of acyclic cluster algebras, 
which is of special interest.
As before, let $Q$ be an acyclic quiver with vertices $1,\ldots,n$.
Without loss of generality we assume that $i < j$ whenever there is an
arrow $a\df i \to j$ in $Q$.
We define two Weyl group elements $c := s_n \cdots s_2s_1$ and $w := c^2$.
For simplicity we assume that $Q$ is not a linearly oriented quiver of type
$\A_n$.
This implies that $\ii := (n,\ldots,2,1,n,\ldots,2,1)$ is
a reduced expression of $w$.
Define $V_\ii = V_1 \oplus \cdots \oplus V_{2n}$ and
$M_\ii = M_1 \oplus \cdots \oplus M_{2n}$ as before.

It follows that for $1 \le j \le n$ we have $M_j = I_j$ and
$M_{n+j} = \tau_Q(I_j)$.
Here $I_j$ denotes the indecomposable injective $KQ$-module with
socle $S_j$, and $\tau_Q$ is the Auslander-Reiten translation in
$\md(KQ)$.

Observe that $\RR(\CC_w,V_\ii)$ is an acyclic cluster algebra
associated to $Q$ having $n$ non-invertible coefficients,
whereas $\underline{\RR}(\CC_w,V_\ii)$ is the acyclic
cluster algebra associated to $Q$ having no coefficients.

\begin{Thm}\label{thm16.1}
With $w$ and $\ii$ as above, the following hold:
\begin{itemize}

\item[(i)]
$\RR(\CC_w,V_\ii) = \C[\delta_{M_1},\ldots,\delta_{M_{2n}}]
= \Span_\C\ebrace{\delta_X \mid X \in \CC_w}$;

\item[(ii)]
$\left\{ \delta_M \mid M \in \add(M_\ii) \right\}$ and
$\left\{ s_M \mid M \in \add(M_\ii) \right\}$ 
are both a $\C$-basis 
of $\RR(\CC_w,V_\ii)$;

\item[(iii)]
$\left\{ \underline{s}_M \mid M \in \add(M_\ii), 
M \text{ has gaps} \right\}$ 
is a $\C$-basis of $\underline{\RR}(\CC_w,V_\ii)$;

\item[(iv)]
$\left\{ \underline{\delta}_M \mid M \in \add(M_\ii), 
M \text{ has  gaps} \right\}$ 
is a $\C$-basis of $\underline{\RR}(\CC_w,V_\ii)$;

\item[(v)]
There is an isomorphism of cluster algebras
$\underline{\RR}(\CC_w,V_\ii) \cong \cA_Q$,
where $\cA_Q$ is the coefficient-free acyclic 
cluster algebra associated to $Q$.

\end{itemize}
\end{Thm}

\begin{proof}
Parts (i), (ii) and (iii) were already proved before for arbitrary 
reduced expressions of arbitrary Weyl group elements.
Part (v) is clear from our description of the initial seed
(labeled by $V_\ii$) for the cluster algebra $\RR(\CC_w,V_\ii)$.
It remains to prove (iv):
We have
$$
\RR(\CC_w,V_\ii) = \bigoplus_{d \in \N^n} \RR_d
$$
where $\RR_d$ is the $\C$-vector space with basis
$\left\{ s_M \mid M \in \add(M_\ii) \cap \rep(Q,d) \right\}$.
We know that
$\left\{ \delta_M \mid M \in \add(M_\ii) \cap \rep(Q,d) \right\}$
is a basis of $\RR_d$ as well.
After specializing the coefficients $\delta_{I_{\ii,j}}$, $1 \le j \le n$ 
to 1,
we get
$$
\underline{\RR}(\CC_w,V_\ii) = \bigoplus_{d \in \N^n} 
\underline{\RR}_d
$$
where $\underline{\RR}_d$ is the $\C$-vector space with
basis
$$
\left\{ \underline{s}_M \mid M \in \add(M_\ii) \cap \rep(Q,d),
M \text{ has  gaps} \right\}.
$$
Now one can use the formula
$$
\delta_{I_{\ii,i}} = \delta_{M_{n+i}} \cdot \delta_{M_i} -
\prod_{j \to i} \delta_{M_{n+j}} \cdot
\prod_{i \to k} \delta_{M_k}
$$
(where the products are taken over all arrows of $Q$ which
start and end in $i$, respectively) and 
an induction on the number of vertices
of $Q$ to show that for every  $M \in \add(M_\ii)$ which has  gaps,
the vector $\underline{s}_M$ is a linear combination
of elements of the form $\underline{\delta}_{M'}$ where $M'$  
in $\add(M_\ii)$ has gaps and $|\dimv(M')| \leq |\dimv(M)|$.
For dimension reasons we get that the vectors $\underline{\delta}_{M'}$
with $M'$ having  gaps form a linearly independent set.
This implies (iv).
\end{proof}

It is interesting to compare Theorem~\ref{thm16.1},(iv) 
to Berenstein, Fomin and
Zelevinsky's construction of a basis for 
the acyclic cluster algebra $\cA_Q$.
Let ${\mathbf y} := (y_1,\ldots,y_n)$ be the initial
cluster whose exchange matrix $B_Q$ is encoded by $Q$, as in
Section~\ref{clustintro}.
Let $y_1^*,\ldots,y_n^*$ be the $n$ cluster variables
obtained from $\mathbf{y}$ by one mutation in each of the $n$ possible
directions.  
Thus the $n$ sets $\{ y_1,\ldots,y_n \} \setminus \{ y_k \} \cup\{ y_k^* \}$
form 
the neighboring clusters of our initial cluster
$\mathbf{y}$.
Using a simple Gr\"obner basis argument, the following is shown
in \cite{BFZ}:

\begin{Thm}[Berenstein, Fomin, Zelevinsky]\label{BFZbasis}
The monomials
$$
\left\{ y_1^{p_1}(y_1^*)^{q_1} \cdots y_n^{p_n}(y_n^*)^{q_n} \mid
p_i,q_i \ge 0,\; p_iq_i = 0 \right\}
$$
form a $\C$-basis of the acyclic cluster algebra $\cA_Q$.
\end{Thm}

Starting with the initial seed $(\yy,B_Q)$, which corresponds to 
$\GG_\ii \equiv \GG_{V_\ii}$,
we perform the sequence of  mutations $\mu_n \cdots \mu_2\mu_1$.
In each step we obtain a new cluster variable which we denote by
$y_k^\dag$.
Note that $y_1^\dag = y_1^*$, but already $y_2^\dag$ and $y_2^*$
may be different.
Observe that $\mu_n \cdots \mu_2\mu_1(B_Q) = B_Q$.
We get that
$$
((y_1^\dag,\ldots,y_n^\dag),B_Q)
$$
is a seed of the cluster algebra $\cA_Q$
where 
$$
\{ y_1,\ldots,y_n \} \cap \{ y_1^\dag,\ldots,y_n^\dag \} = \varnothing.
$$
Our version of Theorem~\ref{BFZbasis} looks then as follows:

\begin{Thm}
The monomials
$$
\left\{ 
y_1^{p_1}(y_1^\dag)^{q_1} \cdots y_n^{p_n}(y_n^\dag)^{q_n} \mid
p_i,q_i \ge 0, \; p_iq_i = 0 
\right\}
$$
form a $\C$-basis of the acyclic cluster algebra $\cA_Q$.
\end{Thm}

Note that
the initial cluster $(y_1,\ldots,y_n)$ comes from $V_\ii$
and the cluster 
$(y_1^\dag,\ldots,y_n^\dag)$ 
comes from
$T_\ii$.


\section{Coordinate rings of unipotent radicals}


In this section, we assume that $Q$ is of finite Dynkin type $\A, \D, \E$.
We first recall some standard notation (we refer the reader to 
\cite{GLSFlag} for more details).
The group $G$ is now a simple complex algebraic group of the same type as $Q$.
Let $J$ be a subset of the set $I$ of vertices of $Q$, and let $K$ be
the complementary subset. To $K$ one can attach a standard parabolic
subgroup $B_K$ containing the Borel subgroup $B=HN$. We denote by
$N_K$ the unipotent radical of $B_K$. This is a subgroup of~$N$.
Let $W_K$ be the subgroup of the Weyl group $W$ generated by the 
reflections $s_k$ with $k \in K$. 
This is a finite Coxeter group and
we denote its longest element by $w_0^K$. 
The longest element of $W$ is denoted by $w_0$.

In finite type, the preprojective algebra $\LL$ is finite-dimensional
and selfinjective.
In agreement with \cite{GLSFlag}, we shall denote by $P_i$ the
indecomposable projective $\LL$-module with top $S_i$ and by
$Q_i$ the indecomposable injective module with socle $S_i$.
We write 
$$
Q_J = \bigoplus_{j\in J} Q_j 
\text{\;\;\; and \;\;\;}
P_J = \bigoplus_{j\in J} P_j.
$$

In \cite{GLSFlag} we have shown that $\C[N_K]$ is naturally
isomorphic to the subalgebra
\[
R_K := \Span_\C\ebrace{\vph_X \mid X\in\Cogen(Q_J)}
\]
of $\C[N]$.
As before, $\Cogen(Q_J)$ is the full subcategory of $\md(\LL)$ whose objects
are submodules of direct sums of finitely many copies of $Q_J$.
This allowed us to introduce a cluster algebra $\cA_J \subseteq R_K$,
whose cluster monomials are of the form $\vph_X$ with $X$ a rigid
module in $\Cogen(Q_J)$. 
We conjectured that in fact $\cA_J = R_K$, see \cite[Conjecture 9.6]{GLSFlag}.

We are going to prove that this conjecture follows from the results
of this paper.
Let $w := w_0w_o^K$, and let $\ii$ be a reduced expression for $w$.

\begin{Lem}
We have $N_K = N'(w_0^K) = N(w_0w_0^K)$.
\end{Lem}

\begin{proof}
We know that $N'(w_0^K)$ is the subgroup of $N$ 
generated by the one-parameter subgroups $N(\alpha)$ with $\alpha>0$
and $w_0^K(\alpha)>0$. These are exactly the one-parameter subgroups
of $N$ which do not belong to the Levi subgroup of $B_K$, hence the
first equality follows.
Now, since $N=w_0N_-w_0$, we have 
$$
 N'(w_0^K) = N \cap \left(w_0^KNw_0^K\right) = 
N \cap \left(w_0^Kw_0N_{-}w_0w_0^K\right) = N(w_0w_0^K).
$$
\end{proof}

As before, let 
$\Gen(P_J)$ be
the subcategory of $\md(\LL)$ whose objects
are factor modules of direct sums of finitely many copies of $P_J$.

\begin{Lem}
We have $\CC_{w_0w_0^K} = \Gen(P_J)$.
\end{Lem}

\begin{proof}
By Proposition~\ref{minors}, we know that the indecomposable 
$\CC_w$-projective-injective object $I_{\ii,i}$ with socle $S_i$
satisfies 
$$
\vph_{I_{\ii,i}} = D_{\vpi_i,w_0^Kw_0(\vpi_i)},\qquad (i\in I).
$$
By \cite[\S 6.2]{GLSFlag}, it follows that 
$I_{\ii,i} = {\mathcal E}_{w_0^K}Q_i$,
where ${\mathcal E}_{w_0^K}$ is the functor defined in \cite[\S 5.2]{GLSFlag}.
It readily follows that $I_{\ii,i}$ is the 
indecomposable projective-injective 
object of $\Gen(P_J)$ with simple socle $S_i$. 
Hence $\CC_{w_0w_0^K}$ and $\Gen(P_J)$ have the same projective-injective
generator.
\end{proof}

Let $S$ denote the self-duality of $\md(\LL)$ induced 
by the involution $a \mapsto a^*$
mapping an arrow $a$ of $\overline{Q}$ to its opposite arrow $a^*$,
see \cite[\S 1.7]{GLSAus}.
This restricts to an
anti-equivalence of categories 
$\Gen(P_J) \to \Cogen(Q_J)$, that we shall again denote by $S$.

\begin{Lem}\label{lem:inverse}
For every $X \in \nil(\LL)$ and every $n \in N$ we have
\[
\vph_{X}(n^{-1})=(-1)^{\dm X}\vph_{S(X)}(n).
\]
\end{Lem}

\begin{proof}
We know that $N$ is generated by the one-parameter subgroups $x_i(t)$
attached to the simple positive roots.
By Proposition~\ref{phi_form} we have
\[
\vph_X(x_{i_1}(t_1)\cdots x_{i_k}(t_k))=
\sum_{\aaa=(a_1,\ldots,a_k)\in\N^k}\chi_{\rm c}(\F_{\ii^\aaa,X})
\frac{t_1^{a_1}\cdots t_k^{a_k}}{a_1!\cdots a_k!}.
\] 
Now, if $n=x_{i_1}(t_1)\cdots x_{i_k}(t_k)$, we have
$n^{-1}=x_{i_k}(-t_k)\cdots x_{i_1}(-t_1)$ and the result 
follows from the fact that 
$\F_{\ii^\aaa,X} \cong \F_{\ii_\op^{\aaa_\op},S(X)}$,
where $\ii_\op$ and $\aaa_\op$ denote the sequences 
obtained by reading $\ii$ and $\aaa$ from right to left.
\end{proof}

We can now prove the following:

\begin{Thm}
Conjecture  9.6 of \cite{GLSFlag} holds.
\end{Thm}

\begin{proof}
As before, let $w := w_0w_0^K$, and let
$\ii$ be a reduced expression of $w$.
The cluster algebra $\RR(\CC_w) = \RR(\Gen(P_J))$ is 
isomorphic to $\cA_J$ via the map $\vph_X \mapsto \vph_{S(X)}$.
This comes from the fact that $S\df \Gen(P_J) \to \Cogen(Q_J)$
is an anti-equivalence which maps the $\CC_w$-maximal rigid module
$V_\ii$ used to define the initial seed of $\RR(\CC_w)$
to the maximal rigid module
$U_\jj$ of \cite[\S 9.2]{GLSFlag} used to define the initial seed of $\cA_J$. 
(Here we assume that $\jj$ is the reduced expression of $w_0^Kw_0$
obtained by reading the reduced expression $\ii$ of $w_0w_0^K$
from right to left.)
In particular the cluster variables $\vph_{M_k}$ which,
by Theorem~\ref{proofmain1}, generate
$\RR(\Gen(P_J)) \equiv \C[N(w_0w_0^K)]$ are mapped to cluster variables
$\vph_{S(M_k)}$ of $\cA_J$. 
They also form a system of generators of the 
polynomial algebra $\C[N(w_0w_0^K)] = \C[N_K]$ by Lemma~\ref{lem:inverse},
because the map $n\mapsto n^{-1}$ is a biregular automorphism of $N_K$.
Hence $\cA_J = \C[N_K]$.  
\end{proof}

\begin{Rem}\label{finrem}
The previous discussion shows that we obtain two different cluster
algebra structures on $\C[N_K]$, coming from the two different subcategories
$\Gen(P_J)$ and $\Cogen(Q_J)$.
When using $\Gen(P_J) = \CC_{w_0w_0^K}$, we regard $\C[N_K]$ as the subring
of $N'(w_0w_0^K)$-invariant functions of $\C[N]$ for the action of
$N'(w_0w_0^K)$ on $N$ by \emph{right} translations, 
see Section~\ref{invariant_ring}.  
This allows us to relate the first cluster structure to the cluster structure
of the unipotent cell $\C[N^{w_0w_0^K}]$, see Proposition~\ref{coordNw}.
When using $\Cogen(Q_J)$, we regard $\C[N_K]$ as the subring
of $N'(w_0w_0^K)$-invariant functions of $\C[N]$ for the action of
$N'(w_0w_0^K)=N(w_0^K)$ on $N$ by \emph{left} translations.
These functions can then be ``lifted'' to $B_K^{-}$-invariant functions
on $G$ for the action of $B_K^{-}$ on $G$ by left translations.
This allows us to ``lift'' the second cluster structure
to a cluster structure on $\C[B_K^-\backslash G]$, see \cite[\S 10]{GLSFlag}.
\end{Rem}


\section{Remarks and open problems}
\label{openproblems}


\subsection{Calculation of $M_\ii(R)$}
Let $\ii$ be a reduced expression of a Weyl group element
$w$, and let $R$ be a $V_\ii$-reachable
$\LL$-module, see Section~\ref{reach}.
Based on Theorem~\ref{main7} we
can combine Corollary~\ref{QTexact} and Proposition~\ref{dimcount}
to determine algorithmically $M_\ii(R)$.
(For the definition of $M_\ii(R)$ see Section~\ref{Mstrata}.)
Recall that the $V_\ii$-reachable modules $R$ 
are in 1-1 correspondence with the cluster monomials
$\delta_R$ in $\RR(\CC_w)$.

\subsection{Calculation of Euler characteristics}
\label{Euler2}
Let $\ii$ be a reduced expression of a Weyl group element
$w$,
and let $R$ be a $V_\ii$-reachable $\LL$-module,
and let $\jj = (j_1,\ldots,j_p)$. 
By Proposition~\ref{phi_form}
the Euler characteristic
$\chi_{\rm c}(\F_{\jj,R})$ is equal to the coefficient of $t_1 \cdots t_p$ in
$\vph_R(x_{j_1}(t_1)\cdots x_{j_p}(t_p))$.
Using mutations, we can express algorithmically $\vph_R$
as a Laurent polynomial in the functions $\vph_{V_1},\ldots,\vph_{V_r}$.
Now we can use the calculations from Section~\ref{Euler1} to compute
all the Euler characteristics 
$\chi_{\rm c}(\F_{\jj,R})$.

\subsection{Open orbit conjecture}
It is known that the (specialized) dual canonical basis ${\mathcal B}^*$ 
and the dual semicanonical
basis $\cS^*$ of $\MM^* \equiv U(\n)^*_{\rm gr}$ 
do not coincide, see \cite[Section 1.5]{GLSSemi1}.
But one might at least hope that both bases have some interesting
elements in common:

\begin{Conj}[Open Orbit Conjecture]\label{finalconj2}
Let $Z$ be an irreducible component of $\LL_d$,
and let $b_Z$ and $\rho_Z$ be the associated dual canonical
and dual semicanonical basis vectors of $\MM^*$.
If $Z$ contains an open $\GL_d$-orbit, then
$b_Z = \rho_Z$.
\end{Conj}

We know that each cluster monomial of the cluster algebra 
$\cA(\CC_w)$ is of the form $\rho_Z$, where
$Z$ contains an open $\GL_d$-orbit.
So if the conjecture is true, then all cluster monomials belong
to the dual canonical basis.

\subsection{Example}
Finally, we would like to ask the following question.
Is it possible to realize every element of
the dual canonical basis of $\MM^*$ as 
a $\delta$-function?
We know several examples of elements $b$ of ${\mathcal B}^*$ which
do not belong to $\cS^*$. In all these examples,
$b$ is however equal to $\delta_X$ for a non-generic $\LL$-module $X$.
(We say that $X \in \nil(\LL)$ is {\it generic} if $\delta_X \in \cS^*$.)

Let us look at an example.
Let $Q$ be the quiver 
$\xymatrix@-0.5pc{1 & 2 \ar@<0.5ex>[l]\ar@<-0.5ex>[l]}$
and let $\LL$ be the associated preprojective algebra.
For $\lambda \in \C^*$ we define representations $M(\lambda,1)$ and
$M(\lambda,2)$ of $Q$ as follows:
$$
M(\lambda,1) := \;
\xymatrix@+1.0pc{
\C & \C \ar@<0.5ex>[l]^{\left(\bsm \lambda \esm\right)}
\ar@<-0.5ex>[l]_{\left(\bsm 1 \esm\right)}
}
\text{\;\; and \;\;}
M(\lambda,2) := \;
\xymatrix@+1.0pc{
\C^2 & \C^2 \ar@<0.5ex>[l]^{\left(\bsm \lambda & 1\\0&\lambda\esm\right)}
\ar@<-0.5ex>[l]_{\left(\bsm 1 &0\\0&1\esm\right)}
}
$$
Let $\iota\df \rep(Q,(2,2)) \to \LL_{(2,2)}$ be the obvious canonical
embedding.
Clearly, the image of $\iota$ is an irreducible component
of $\LL_{(2,2)}$, which we denote by $Z_Q$.
It is not difficult to check that the set
$$
\left\{ M(\lambda,1) \oplus M(\mu,1) \mid \lambda,\mu \in \C^* \right\}
$$
is a dense subset of $Z_Q$.
It follows that 
$$
\delta_{M(\lambda,1) \oplus M(\mu,1)} = \rho_{Z_Q}
$$
is an element of the dual semicanonical basis $\cS^*$.
It is easy to check that 
$$
\delta_{M(\lambda,2)}\not = \delta_{M(\lambda,1) \oplus M(\mu,1)}.
$$
Indeed, the variety $\F_{\jj,X}$ of composition series of type 
$\jj = (1,2,1,2)$ 
has Euler characteristic $2$ for $X = M(\lambda,1) \oplus M(\mu,1)$
and Euler characteristic $1$ for $X = M(\lambda,2)$.
Furthermore, one can show that
$$
\delta_{M(\lambda,2)} = b_{Z_Q}
$$ belongs to the dual canonical basis 
${\mathcal B}^*$ of $\MM^*$.

Note that the functions $\delta_{M(\lambda,1) \oplus M(\mu,1)}$
and $\delta_{M(\lambda,2)}$ do not belong to any of the
subalgebras $\RR(\CC_w)$, since 
$M(\lambda,1)$ and $M(\lambda,2)$ are regular representations
of the quiver $Q$ for all $\lambda$.


\bigskip
{\parindent0cm \bf Acknowledgements.}\,
We are grateful to {\O}yvind Solberg for answering our
questions on relative homology theory.
We thank Shrawan Kumar for his kind help concerning
Kac-Moody groups.
It is a pleasure to thank the Mathematisches Forschungsinstitut
Oberwolfach (MFO) for two weeks of hospitality  
in July/August 2006, where this work was started.
Furthermore, the first and second authors like to thank the 
Max-Planck Institute for Mathematics in Bonn for a research stay in
September - December 2006 and October 2006, respectively.
We also thank the Sonderforschungsbereich/Transregio SFB 45 for financial 
support.
The three authors are grateful to the Mathematical
Sciences Research Institute in Berkeley (MSRI) for 
an invitation in Spring 2008 during which substantial parts of 
this work were written. 
The first author was partially supported by PAPIIT grant IN103507-2
and CONACYT grant 81948.
All three authors thank the Hausdorff Center for Mathematics
in Bonn for financial support.



\begin{thebibliography}{999}


\bibitem[Ab]{Ab}
{\it E. Abe},
Hopf algebras.
Cambridge Tracts in Mathematics 74.
Cambridge University Press, 1980.

\bibitem[ASS]{ASS}
{\it I. Assem, D. Simson, A. Skowro{\'n}ski},
Elements of the representation theory of associative algebras. 
Vol. 1. Techniques of representation theory. 
London Mathematical Society Student Texts, 65. 
Cambridge University Press, Cambridge, 2006. x+458 pp.

\bibitem[Au]{A}
{\it M. Auslander},
Representation theory of artin algebras. I. 
Comm. Alg. 1 (1974), 177--268. 

\bibitem[APR]{APR}
{\it M. Auslander, M. Platzeck, I. Reiten},
Coxeter functors without diagrams.
Trans. Amer. Math. Soc. 250 (1979), 1--46.

\bibitem[ARS]{ARS}
{\it M. Auslander, I. Reiten, S. Smal\o}, 
Representation theory of Artin algebras. 
Corrected reprint of the 1995 original. 
Cambridge Studies in Advanced Mathematics, 36. Cambridge University 
Press, Cambridge, 1997. xiv+425 pp. 

\bibitem[AS1]{AS1}
{\it M. Auslander, \O. Solberg},
Relative homology and representation theory. I.
Relative homology and homologically finite subcategories.
Comm. Algebra 21 (1993), no. 9, 2995--3031.

\bibitem[AS2]{AS2}
{\it M. Auslander, \O. Solberg},
Relative homology and representation theory. II.
Relative cotilting theory.
Comm. Algebra 21 (1993), no. 9, 3033--3079.


\bibitem[BZ]{BZ}
{\it A. Berenstein, A. Zelevinsky},
Total positivity in Schubert varieties.
Comment. Math. Helv. 72 (1997), 128--166.

\bibitem[BFZ]{BFZ}
{\it A. Berenstein, S. Fomin, A. Zelevinsky},
Cluster algebras III: Upper bounds and double Bruhat
cells.
Duke Math. J. 126 (2005), no. 1, 1--52.

\bibitem[Bo]{Bo}
{\it K. Bongartz}, 
Minimal singularities for representations of Dynkin quivers. 
Comment. Math. Helv. 69 (1994), no.4, 575--611.

\bibitem[BMRRT]{BMRRT}
{\it A. Buan, R. Marsh, M. Reineke, I. Reiten, G. Todorov},
Tilting theory and cluster combinatorics.  
Adv. Math. 204 (2006), no. 2, 572--618.

\bibitem[BM]{BM}
{\it A. Buan, R. Marsh},
Cluster-tilting theory.
In: Trends in Representation Theory of Algebras and Related Topics, 
Contemp.~Math. 406 (2006), 1--30.

\bibitem[BIRS]{BIRS}
{\it A. Buan, O. Iyama, I. Reiten, J. Scott},
Cluster structures for 2-Calabi-Yau categories and
unipotent groups.
Compositio Math. 145 (2009), 1035--1079.

\bibitem[CK1]{CK}
{\it P. Caldero, B. Keller},
From triangulated categories to cluster algebras.
Invent.~Math. 172 (2008), 169--211.

\bibitem[CK2]{CK2}
{\it P. Caldero, B. Keller},
From triangulated categories to cluster algebras II.
Ann.~Sci. {\'E}cole Norm.~Sup. (4) 39 (2006), no. 6, 983--1009. 

\bibitem[CPS]{CPS}
{\it E. Cline, B. Parshall, L. Scott},
Finite-dimensional algebras and highest weight categories.  
J. Reine Angew. Math. 391 (1988), 85--99. 

\bibitem[CB]{CB}
{\it W. Crawley-Boevey}, 
On the exceptional fibres of Kleinian singularities.
Amer.~J. Math. 122 (2000), 1027--1037.

\bibitem[CBS]{cbsr01}
{\it W. Crawley-Boevey, J.~Schr{\"o}er}, 
Irreducible components of varieties of modules. 
J. Reine Angew. Math. 553 (2002), 201--220.

\bibitem[DWZ1]{DWZ1}
{\it H. Derksen, J. Weyman, A. Zelevinsky},
Quivers with potentials and their representations I: Mutations.
Selecta Math. (N.S.) 14 (2008), no. 1, 59--119.

\bibitem[DWZ2]{DWZ2}
{\it H. Derksen, J. Weyman, A. Zelevinsky},
Quivers with potentials and their representations II: 
Applications to cluster algebras.
J. Amer.~Math.~Soc. 23 (2010), no. 3, 749--790.

\bibitem[DFK]{DFK}
{\it P. Di Francesco, R. Kedem},
Q-systems as cluster algebras II: Cartan matrix of finite type and 
the polynomial property.
Lett.~Math.~Phys. 89 (2009), no.~3, 183--216.

\bibitem[FZ1]{FZ5}
{\it S. Fomin, A. Zelevinsky}, 
Double Bruhat cells and total positivity.
J. Amer.~Math.~Soc. 12 (1999), no. 2, 335--380.

\bibitem[FZ2]{FZ1}
{\it S. Fomin, A. Zelevinsky}, 
Cluster algebras. I. Foundations. 
J. Amer.~Math.~Soc. 15 (2002), no. 2, 497--529.

\bibitem[FZ3]{FZ2}
{\it S. Fomin, A. Zelevinsky}, 
Cluster algebras. II. Finite type classification.
Invent.~Math. 154 (2003), no. 1, 63--121. 

\bibitem[FZ4]{FZSurv}
{\it S. Fomin, A. Zelevinsky}, 
Cluster algebras: notes for the CDM-03 conference.  
Current developments in mathematics, 2003,  1--34, Int. Press,
Somerville, MA, 2003.

\bibitem[FZ5]{FZ3}
{\it S. Fomin, A. Zelevinsky}, 
Cluster algebras. IV. Coefficients.
Compositio Math. 143 (2007), 112--164.

\bibitem[FK]{FK}
{\it C. Fu, B. Keller},
On cluster algebras with coefficients and 2-Calabi-Yau categories.
Trans.~Amer. Math.~Soc. 362 (2010), 859--895. 

\bibitem[Ga]{gabr74}
{\it P. Gabriel}, 
Finite representation type is open. 
Representations of algebras, 
(Proc. Internat. Conf., Carleton Univ., Ottawa, Ont., 1974, Springer-Verlag, 
Berlin, 1975, Lecture Notes in Math., 488, pp.~132--155.

\bibitem[GR]{GR}
{\it P. Gabriel, A.V. Roiter},
Representations of finite-dimensional algebras.
Translated from the Russian.
With a chapter by B. Keller.
Reprint of the 1992 English translation.
Springer-Verlag, Berlin, 1997. iv+177 pp.

\bibitem[Ge]{geiss96}
{\it C. Gei{\ss}}, 
Geometric methods in representation theory of finite-dimensional algebras. 
Representation theory of algebras and related topics
(Mexico City, 1994), CMS Conf.~Proc., vol.~19, Amer. Math. Soc.,
Providence, RI, 1996, pp.~53--63.

\bibitem[GLS1]{GLSSemi1}
{\it C. Gei{\ss}, B. Leclerc, J. Schr\"oer},
Semicanonical bases and preprojective algebras.
Ann.~Sci. {\'E}cole Norm.~Sup. (4)  38  (2005),  no. 2, 193--253.

\bibitem[GLS2]{GLSAus}
{\it C. Gei{\ss}, B. Leclerc, J. Schr\"oer},
Auslander algebras and initial seeds for cluster algebras.
J. London Math. Soc. (2) 75 (2007), 718--740.

\bibitem[GLS3]{GLSVerma}
{\it C. Gei{\ss}, B. Leclerc, J. Schr\"oer},
Verma modules and preprojective algebras.
Nagoya Math.~J. 182 (2006), 241--258.

\bibitem[GLS4]{GLSSemi2}
{\it C. Gei{\ss}, B. Leclerc, J. Schr\"oer},
Semicanonical bases and preprojective algebras II: A multiplication
formula.
Compositio Math. 143 (2007), 1313--1334.

\bibitem[GLS5]{GLSRigid}
{\it C. Gei{\ss}, B. Leclerc, J. Schr\"oer},
Rigid modules over preprojective algebras.
Invent.~Math. 165 (2006), no. 3, 589--632.

\bibitem[GLS6]{GLSFlag}
{\it C. Gei{\ss}, B. Leclerc, J. Schr\"oer},
Partial flag varieties and preprojective algebras.
Ann.~Inst. Fourier (Grenoble) 58 (2008), 825--876.

\bibitem[GLS7]{GLSUni1}
{\it C. Gei{\ss}, B. Leclerc, J. Schr\"oer},
Cluster algebra structures and semicanonical bases for unipotent
groups. 
Unpublished (2007), 121 pp. arXiv:math/0703039

\bibitem[Ha1]{H2}
{\it D. Happel},
Triangulated categories in the representation theory of 
finite-dimensional algebras.
London Mathematical Society Lecture Note Series, 119. 
Cambridge University Press, Cambridge, 1988. x+208 pp.

\bibitem[Ha2]{H3}
{\it D. Happel},
Partial tilting modules and recollement.
Proceedings of the International Conference on Algebra, Part 2 
(Novosibirsk, 1989), 345--361, Contemp.~Math., 131, Part 2, Amer. Math. Soc.,
Providence, RI, 1992.

\bibitem[Iy]{Iy}
{\it O. Iyama},
Auslander correspondence.
Adv. Math. 210 (2007), no. 1, 51--82.

\bibitem[IR]{IR}
{\it O. Iyama, I. Reiten}, 
2-Auslander algebras associated with reduced words in Coxeter
groups.
Preprint (2010), 14pp. arXiv:1002.3247.

\bibitem[Jo]{J}
{\it A. Joseph},
Quantum groups and their primitive ideals.
Ergebnisse der Mathematik und ihrer Grenzgebiete (3), 29. 
Springer-Verlag, Berlin, 1995. x+383pp.

\bibitem[Ka]{Ka}
{\it V. Kac},
Infinite-dimensional Lie algebras.
Third edition.
Cambridge University Press, Cambridge, 1990.
xxii+400pp.

\bibitem[KP]{KP1}
{\it V. Kac, D. Peterson},
Regular functions on certain infinite-dimensional groups.
Arithmetic and geometry, Vol II, 141--166, Progr. Math., 36,
Birkh\"auser Boston, Boston, MA, 1983.

\bibitem[KS]{KSa}
{\it M. Kashiwara, Y. Saito},
Geometric construction of crystal bases.
Duke Math.~J. 89 (1997), 9--36.

\bibitem[Ke]{Ke}
{\it B. Keller},
On triangulated orbit categories.  
Doc. Math. 10 (2005), 551--581 (electronic).

\bibitem[KR]{KR}
{\it B. Keller, I. Reiten},
Cluster tilted algebras are Gorenstein and stably Calabi-Yau.
Adv.~Math. 211 (2007), 123--151.

\bibitem[Ku]{Ku}
{\it S. Kumar},
Kac-Moody groups, their flag varieties and representation theory. 
Progress in Mathematics, 204. Birkh\"auser Boston, Inc., Boston, MA, 2002. 

\bibitem[Le]{LeMZ}
{\it B. Leclerc},
Dual canonical bases, quantum shuffles and $q$-characters.  
Math.~Z.  246  (2004), 691--732. 

\bibitem[Lu1]{Lu1}
{\it G. Lusztig},
Quivers, perverse sheaves, and quantized enveloping algebras. 
J. Amer.~Math.~Soc. 4 (1991), no. 2, 365--421.

\bibitem[Lu2]{Lu2}
{\it G. Lusztig},
Semicanonical bases arising from enveloping algebras. 
Adv.~Math. 151 (2000), no. 2, 129--139.

\bibitem[Pa]{P}
{\it Y. Palu},
Cluster characters for 2-Calabi-Yau triangulated categories.  
Ann. Inst. Fourier (Grenoble)  58  (2008),  no. 6, 2221–2248.

\bibitem[Re]{Reu}
{\it C. Reutenauer}, 
Free Lie algebras. 
London Mathematical Society Monographs.
New Series, 7. Oxford Science Publications. The Clarendon Press, Oxford
University Press, New York, 1993. xviii+269 pp.

\bibitem[Rm]{R}
{\it N. Richmond},
A stratification for varieties of modules.
Bull.~London Math.~Soc. 33 (2001), no. 5, 565--577.

\bibitem[Ri1]{Ri1}
{\it C.M. Ringel},
Tame algebras and integral quadratic forms. 
Lecture Notes in Mathematics, 1099. 
Springer-Verlag, Berlin, 1984. xiii+376 pp. 

\bibitem[Ri2]{Ri5}
{\it C.M. Ringel},
The category of modules with good filtrations over a quasi-hereditary 
algebra has almost split sequences.  
Math.~Z.  208  (1991),  no. 2, 209--223.

\bibitem[Ri3]{Ri4}
{\it C.M. Ringel},
The category of good modules over a quasi-hereditary algebra.  
Proceedings of the Sixth International Conference on Representations of 
Algebras (Ottawa, ON, 1992),  17 pp., 
Carleton-Ottawa Math. Lecture Note Ser., 14, Carleton Univ., Ottawa, ON, 1992. 

\bibitem[Ri4]{Ri6}
{\it C.M. Ringel},
PBW-bases of quantum groups.
J.~Reine Angew. Math. 470 (1996), 51--88.

\bibitem[Ri5]{Ri7}
{\it C.M. Ringel},
Iyama's finiteness theorem via strongly quasi-hereditary algebras.
Preprint (2009), 5~pp.  arXiv:0912.5001.

\end{thebibliography}
\end{document}